\magnification=\magstep1
\font\sc=cmcsc10
 1
\font\Bf=cmbx10 scaled \magstep 1
\def\UN{{\bf 1}}
\def\cD{{\cal D}}
\def\Sl{\leftharpoonup}
\def\Sr{\rightharpoonup}

\def\qed{\ \vrule height 5pt width 5pt depth 0pt}
\def\m@th{\mathsurround=0pt }
\def\leftharpoonupfill{$\m@th \mathord\leftharpoonup \mkern-8mu
    \cleaders\hbox{$\mkern-2mu \mathord- \mkern-2mu$}\hfill
    \mkern-8mu \mathord-$}
\def\Cl#1{\vbox{\ialign{##\crcr
    \leftharpoonupfill\crcr\noalign{\kern-1pt\nointerlineskip}
    $\hfil\displaystyle{#1}\hfil$\crcr}}}         
\def\CL#1{\vbox{\ialign{##\crcr
    \leftharpoonupfill\crcr\noalign{\kern-1pt\nointerlineskip}
    $\hfil\scriptstyle{#1}\hfil$\crcr}}}         
\def\rightharpoonupfill{$\m@th \mathord- \mkern-8mu
    \cleaders\hbox{$\mkern-2mu \mathord- \mkern-2mu$}\hfill
    \mkern-8mu \mathord\rightharpoonup$}
\def\Cr#1{\vbox{\ialign{##\crcr
    \rightharpoonupfill\crcr\noalign{\kern-1pt\nointerlineskip}
    $\hfil\displaystyle{#1}\hfil$\crcr}}}         
\def\CR#1{\vbox{\ialign{##\crcr
    \rightharpoonupfill\crcr\noalign{\kern-1pt\nointerlineskip}
    $\hfil\scriptstyle{#1}\hfil$\crcr}}}         
\def\cros{\raise1.9pt\hbox{$\scriptscriptstyle
          >$}\!\raise1.5pt\hbox{$\scriptstyle\triangleleft\,$}}

\centerline{\Bf Larson--Sweedler Theorem}\bigskip
 
\centerline{\Bf and the Role of Grouplike Elements}\bigskip

\centerline{\Bf in Weak Hopf Algebras}\bigskip
\vskip 1truecm

\centerline{{\sc Peter Vecserny\'es}\footnote{}{
E-mail: vecser@rmki.kfki.hu\hfill\break\indent
Supported by the Hungarian Research Fund, OTKA -- T 034 512}}\bigskip
\bigskip

\centerline{Research Institute for Particle and Nuclear Physics, Budapest}

\bigskip
\centerline{H-1525 Budapest 114, P.O.B. 49, Hungary}
  
\vskip 3truecm
\centerline{\bf Abstract}\bigskip
 
We extend the Larson--Sweedler theorem [10] to weak Hopf algebras by proving 
that a finite dimensional weak bialgebra is a weak Hopf algebra iff 
it possesses a non-degenerate left integral. We show that the  category of 
modules over a weak Hopf algebra is autonomous monoidal with semisimple  
unit and invertible modules. We also reveal the connection of invertible 
modules to left and right grouplike elements in the dual weak Hopf algebra. 
Defining distinguished left and right grouplike elements we 
derive the Radford formula [15] for the fourth power of the antipode in a weak 
Hopf algebra and prove that the order of the antipode is finite up to an 
inner automorphism by a grouplike element in the trivial subalgebra $A^T$ 
of the underlying weak Hopf algebra $A$. 

\vfill\eject

\noindent
{\bf 0. Introduction}\medskip

Weak Hopf algebras have been proposed recently [1, 2, 18] as a 
generalization of Hopf algebras by weakening the compatibility conditions 
between the algebra and coalgebra structures of Hopf algebras.
Comultiplication is allowed to be non-unital, $\Delta(\UN)\equiv
\UN^{(1)}\otimes\UN^{(2)}\not= \UN\otimes \UN$, just like in weak quasi 
Hopf algebras [11] and in rational Hopf algebras [19, 8], but the 
comultiplication is coassociative. In exchange for coassociativity, the 
multiplicativity of the counit is replaced by a weaker condition: 
$\varepsilon(ab)=\varepsilon(a\UN^{(1)})\varepsilon(\UN^{(2)}b)$, 
implying that the unit representation is not necessarily one-dimensional and 
irreducible. Like weak quasi and rational Hopf algebras, they
can possess non-integral (quantum) dimensions even in the finite dimensional 
and semisimple cases, which is necessary if we want to recover them
as global symmetries of low-dimensional quantum field theories. In situations
where only the representation category matters, these two concepts are 
equivalent. Nevertheless, just like finite dimensional Hopf algebras, 
finite dimensional weak Hopf algebras (WHA) obey the mathematical beauty of 
giving rise to a self-dual notion: the dual space of a WHA can be 
canonically endowed with a WHA structure. For a recent review, see [12].

Here we continue the study [2] of the structural properties of finite 
dimensional weak Hopf algebras over a field $k$. The main results of 
this paper are:
\smallskip

\item{1.} The generalization of the Larson--Sweedler theorem [10] to WHAs
claiming that a finite dimensional weak bialgebra is a weak Hopf algebra if 
and only if it possesses a non-degenerate left integral.

\item{2.} The characterization of inequivalent invertible modules of WHAs
through left/right grouplike elements in the dual WHA and the proof of 
the semisimplicity of invertible modules, which include the unit module 
serving as a monoidal unit in the monoidal category of left (right) modules.

\item{3.} A finiteness claim about the order of the antipode (up to an
inner automorphism by a grouplike element in the trivial subalgebra) and 
the derivation of the Radford formula [15] in a weak Hopf algebra $A$:
$S^4(a)=\sigma\Sr s^{-1}as\Sl\hat S^{-1}(\sigma),\ a\in A$, 
where $S$ ($\hat S$) is the antipode in $A$ ($\hat A$), and $s$ and 
$\sigma$ are distinguished left grouplike elements in $A$ and in the dual WHA 
$\hat A$, respectively. 

\smallskip
The existence of a non-degenerate left integral $l\in B$ in a finite 
dimensional bialgebra $B$ implies the existence of a non-degenerate left 
integral $\lambda\in\hat B$ in the dual bialgebra $\hat B$ with the property
$\lambda\Sr l=\UN$. Then the formula $S(a):= (\lambda\Sl a)\Sr l, a\in B$ 
gives rise to the antipode for $B$ proving one direction of the Larson-Sweedler
theorem [10]. The proof of the opposite direction [10] involves the structure 
theorem for Hopf modules, which are one-sided $H$-modules and $H$-comodules 
of the Hopf algebra $H$ together with a compatibility condition. The structure 
theorem for a finite dimensional Hopf module $M^H_H$ claims that $M^H_H\simeq
C(M)\otimes H^H_H$ as Hopf modules, where $C(M)$ is the space of coinvariants
in $M^H_H$ and $H^H_H$ is the canonical Hopf module. Proving that the dual 
Hopf algebra $\hat H$ carries a Hopf module structure, $\hat H^H_H\simeq
C(\hat H)\otimes H^H_H$ follows and observing that by dimensionality argument
$C(\hat H)$ is one dimensional, a non-degenerate left integral in $\hat H$
emerges in the space of coinvariants $C(\hat H)$.

The proof of the corresponding statement (Theorem 4.1) in the case of finite 
dimensional weak bialgebras is in the same spirit. The existence of a 
non-degenerate left integral in a finite dimensional WBA implies the 
existence of a non-degenerate left integral in the dual WBA and the previous
classical formula leads to the antipode. The proof of the opposite direction
is more involved: besides weak Hopf modules one has to introduce multiple
weak Hopf modules, in which bimodule or bicomodule structures are also present
together with compatibility conditions between the module and comodule 
structures. Then the structure theorem (Theorem 3.2) for a multiple weak Hopf 
module ${_AM^A_A}$ of a WHA $A$ claims that 
${_AM^A_A}\simeq {_A(C(M)\times A^A_A)}$, i.e., the right weak Hopf module 
structure of $M$ is given by the canonical weak Hopf module $A^A_A$; while
as a left $A$-module, $M$ is isomorphic to the product module of the 
coinvariants ${_AC(M)}$ and the left regular module ${_AA}$. The left 
$A$-module structure of the coinvariants arises from the bimodule structure 
of $M$ (Lemma 3.1 iii).
In particular, the dual WHA $\hat A$ is a multiple weak Hopf module 
${_A{\hat A}^A_A}$ and its coinvariants $C(\hat A)$ are the left integrals
$\hat I^L\subset\hat A$ (Theorem 3.2). Moreover, $\hat I^L$ becomes a free 
left $A^R$- and $A^L$-module with a single generator by restricting the 
left $A$-module structure of ${_A{\hat I}^L}$ to the canonical coideal
subalgebras $A^R$ and $A^L$ of $A$, respectively (Corollary 3.5). It is the 
latter result that replaces the dimensionality argument of the classical 
Hopf case and together with the isomorphism 
${_A{\hat A}^A_A}\simeq {_A(\hat I^L\times A^A_A)}$ of multiple
weak Hopf modules leads to the existence of a non-degenerate left integral in 
$\hat I^L\subset\hat A$.

The modules of a WHA that are invertible with respect to their monoidal product
are important in low dimensional quantum field theories. Hence, it is worth
characterising them in purely (weak) Hopf algebraic terms. Although a 
WHA $A$ is not a semisimple algebra in general, its unit and invertible modules 
are semisimple (Theorem 2.4, resp. Prop. 5.4 ii). The origin of this property 
is that the trivial subWHA $A^T$, which is generated by the canonical coideal 
subalgebras $A^L$ and $A^R$ of a WHA $A$, is in the coradical of $A$ 
(Lemma 2.3). We derive two other equivalent characterizations of invertible 
modules: they are precisely the modules that become free rank one $A^L$- and 
$A^R$-modules by restricting the $A$-module structure to these 
coideal subalgebras (Prop. 5.4 i). For example, the invertible left $A$-module
structure of right integrals $I^R\subset A$ and left integrals $\hat I^L\subset\hat A$ follows in this way. 
The second equivalent characterization of invertible $A$-modules 
involves left or right grouplike elements (Def. 5.1) in the dual WHA: 
an $A$-module is invertible iff it is isomorphic to a cyclic
submodule in the second regular $A$-module ${_A\hat A}$ generated by a 
left (right) grouplike element in $\hat A$ (Prop. 5.7). Moreover,
the isomorphism classes of invertible $A$-modules are given by the (finite)
factor group $G_L(\hat A)/G_L^T(\hat A)$ (or by $G_R(\hat A)/G_R^T(\hat A)$) 
(Prop. 5.7),
where $G_L^T(\hat A)$ is the intersection of the (in general infinite) set of 
left grouplike elements $G_L(\hat A)$ and the trivial subWHA ${\hat A}^T$ 
in $\hat A$. 

If $l\in A$ and $\lambda\in\hat A$ are dual left integrals, i.e., if they are 
non-degenerate and satisfy $\lambda\Sr l=\UN$, then 
$s:=l\Sl\lambda$ and $\sigma:=\lambda\Sl l$ will define (distinguished) 
left grouplike elements (Def. 6.1 and discussion before) like in the Hopf case
[15]. $\sigma$ falls into a central element of the factor group 
$G_L(\hat A)/G_L^T(\hat A)$ and determines the unimodularity of $A$, that is 
the possible existence of a two-sided non-degenerate integral in $A$ (Corollary
6.3). The Nakayama automorphism $\theta_\lambda\colon A\to A$ corresponding 
to a non-degenerate left integral $\lambda\in\hat A$ can be given in terms
of distinguished left grouplike elements in two different ways, which contain
the square or the inverse square of the antipode. Hence, these 
expressions lead to the generalization of the Radford formula [15] to WHAs
(Theorem 6.4). Since the factor groups 
$G_L(A)/G_L^T(A)$ and $G_L(\hat A)/G_L^T(\hat A)$ are finite and since even 
powers of the antipode are WHA automorphisms, the iteration of the Radford 
formula leads to the claim that the order of the antipode is finite up to a 
conjugation by an element in $G_L^T(A)\cap G_R^T(A)$ (Theorem 6.4). 
The explicit form of the Nakayama automorphism $\theta_\lambda$, like in 
the Hopf case [16], can be used to prove the unimodularity of the double 
of a WHA (Corollary 6.5).  

We note that it was established in [2] that WHAs are quasi-Frobenius algebras.
Result 1 implies that they are Frobenius algebras. Grouplike elements in a WHA, 
which are just the intersection of left and right grouplike elements in our 
formulation, were introduced in [2]. The modules associated with them 
were studied in [13]. However, this notion of grouplike elements is too 
restrictive: 
for characterization of isomorphism classes of invertible modules (Result 2) 
one has to introduce the less restrictive notion of left (right) grouplike 
elements, because the factor group $G(\hat A)/G^T(\hat A)$ of grouplike and 
trivial grouplike elements is, in general, smaller than the corresponding 
factor group $G_L(A)/G_L^T(A)$ of left grouplike elements (Prop 5.8). 
Result 3 was proved in [13] in the case when the square 
of the antipode is the identity mapping on the coideal subalgebra $A^L$ of 
the WHA $A$.
  
\smallskip
The organization of the paper is as follows. In Section 1 we review the axioms
and the main properties of weak bialgebras (WBA) and weak Hopf algebras. 
Here and throughout the paper they are considered to be finite dimensional. 
Section 2 is devoted to the autonomous monoidal category of modules of a 
WHA and to properties of the unit module including semisimplicity. We 
derive also a lower bound for the $k$-dimension of an $A$-module in terms
of the $k$-dimensions of the simple submodules of the unit $A$-module. This
estimation leads to a sufficient condition for an $A$-module to become a 
free rank one $A^L$- and $A^R$-module.  
In Section 3 we prove a structure theorem for multiple weak Hopf modules 
and show that the left $A$-modules spanned by right integrals in $A$ and left
integrals in $\hat A$ become free rank one $A^L$- and $A^R$-modules.  
Section 4 contains the generalization of the Larson--Sweedler theorem 
to the weak Hopf case. 
In Section 5 we reveal the connection between invertible 
modules of a WHA $A$ and left (right) grouplike elements in the dual WHA 
$\hat A$ and prove that invertible modules are semisimple. 
Section 6 contains the definition and some basic properties of distinguished 
left and right grouplike elements, the derivations of the form of 
the Nakayama automorphism $\theta_\lambda\colon A\to A$ corresponding to 
a non-degenerate left integral $\lambda\in\hat A$ and the Radford formula. 
In addition, we prove the claim about the order of the 
antipode and unimodularity of the double of a WHA. 
In Appendix A we give a simple example of a WHA 
in which the order of the antipode is not finite. Finally, Appendix B
contains the generalization of the cyclic category module [4] to weak Hopf 
algebras containing a modular pair of grouplike elements in involution.
\bigskip

\vbox{\noindent
{\bf 1. Preliminaries}\medskip

Here we give a quick survey of weak bialgebras and weak Hopf algebras 
[2]. We restrict ourselves to their main properties, however,
some useful identities we use later on are also given.}  

\smallskip\noindent
{\bf 1.1 The axioms}\smallskip

A {\it weak bialgebra} $(A;u,\mu;\varepsilon,\Delta)$ is defined by 
the properties i--iii):

\item{i)} $A$ is a finite dimensional associative algebra over a field $k$ 
with multiplication $\mu\colon A\otimes A\to A$ and unit $u\colon k\to A$, 
which are $k$-linear maps.
\item{ii)} $A$ is a coalgebra over $k$ with comultiplication $\Delta\colon
A\to A\otimes A$ and counit $\varepsilon\colon A\to k$, which are $k$-linear 
maps.
\item{iii)} The algebra and coalgebra structures obey the compatibility
conditions
$$
  \Delta(ab)=\Delta(a)\Delta(b),\qquad a,b\in A\eqno(1.1a)
$$$$
  \varepsilon(ab^{(1)})\varepsilon(b^{(2)}c)
 =\varepsilon(abc)
 =\varepsilon(ab^{(2)})\varepsilon(b^{(1)}c), \qquad a,b,c\in A\eqno(1.1b)
$$$$
  \UN^{(1)}\otimes\UN^{(2)}\UN^{(1')}\otimes\UN^{(2')}
 =\UN^{(1)}\otimes\UN^{(2)}\otimes\UN^{(3)}
 =\UN^{(1)}\otimes\UN^{(1')}\UN^{(2)}\otimes\UN^{(2')},\eqno(1.1c)
$$
where (and later on) $ab\equiv\mu(a,b), \UN:=u(1)$ and we used Sweedler 
notation [17] for iterated coproducts omitting summation indices and a 
summation symbol.

A {\it weak Hopf algebra} $(A;u,\mu;\varepsilon,\Delta;S)$ is a WBA together 
with property iv):
\item{iv)} There exists a $k$-linear map $S\colon A\to A$, called the 
antipode, satisfying
$$\eqalignno{
  a^{(1)}S(a^{(2)})&=\varepsilon(\UN^{(1)}a)\UN^{(2)}, &(1.2a)\cr
  S(a^{(1)})a^{(2)}&=\UN^{(1)}\varepsilon(a\UN^{(2)}),\qquad a\in A. &(1.2b)\cr
  S(a^{(1)})a^{(2)}S(a^{(3)})&=S(a), &(1.2c)\cr}
$$

WBAs and WHAs are self-dual notions, the dual space $\hat A:={\rm Hom}_k(A,k)$
of a WBA (WHA) equipped with structure maps $\hat u,\hat\mu,\hat\varepsilon,
\hat\Delta, (\hat S)$ defined by transposing the structure maps of $A$ by
means of the canonical pairing $\langle\ ,\ \rangle\colon\hat A\times A\to k$
gives rise to a WBA (WHA).
\vbox{
\medskip\noindent
{\bf 1.2 Properties of WBAs}\smallskip

Let $A$ be a WBA. The images $A^{L/R}=\Pi^{L/R}(A)=\bar\Pi^{L/R}(A)$ of 
the projections $\Pi^{L/R}\colon A\to A$ and $\bar\Pi^{L/R}\colon A\to A$
defined by}
$$\eqalign{
  \Pi^L(a)&:=\varepsilon(\UN^{(1)}a)\UN^{(2)},\cr
  \bar\Pi^L(a)&:=\varepsilon(a\UN^{(1)})\UN^{(2)},\cr}\qquad\eqalign{
  \Pi^R(a)&:=\UN^{(1)}\varepsilon(a\UN^{(2)}),\cr
  \bar\Pi^R(a)&:=\UN^{(1)}\varepsilon(\UN^{(2)}a),\cr}\qquad a\in A\eqno(1.3)
$$
are unital subalgebras (i.e. containing $\UN$) of $A$ that commute with each 
other. $A^L$ and $A^R$ are called {\it left} and {\it right subalgebras},
respectively. The image $\Delta(\UN)$ of the unit is in $A^R\otimes A^L$ 
and the coproduct on $A^{L/R}$ reads as:
$$ 
  \Delta(x^L)=\UN^{(1)}x^L\otimes\UN^{(2)},\quad x^L\in A^L,\qquad 
  \Delta(x^R)=\UN^{(1)}\otimes x^R\UN^{(2)},\quad x^R\in A^R.\eqno(1.4)
$$
Hence, $A^L$ and $A^R$ are left and right coideals, respectively, and the 
{\it trivial subalgebra} $A^T:=A^L\vee A^R\subset A$ generated by the coideal
subalgebras $A^L$ and $A^R$ is a subWBA of $A$. 

The maps $\kappa_L\colon A^L\to\hat A^R$ and $\kappa_R\colon A^R\to\hat A^L$ 
given by Sweedler arrows 
$$
  \kappa_L(x^L):=x^L\Sr\hat\UN,\qquad 
  \kappa_R(x^R):=\hat\UN\Sl x^R,\quad x^{L/R}\in A^{L/R}\eqno(1.5)
$$  
are algebra isomorphisms with inverses $\hat\kappa_R$ and $\hat\kappa_L$, 
respectively. Moreover, 
$$\eqalign{
  x^L\Sr\varphi&=(x^L\Sr\hat\UN)\varphi,\cr
  x^R\Sr\varphi&=\varphi(x^R\Sr\hat\UN),\cr}\qquad
 \eqalign{
  \varphi\Sl x^L&=(\hat\UN\Sl x^L)\varphi,\cr
  \varphi\Sl x^R&=\varphi(\hat\UN\Sl x^R),\cr}
  \qquad\varphi\in\hat A, x^{L/R}\in A^{L/R}.\eqno(1.6)
$$
Defining $Z^{L/R}:=A^{L/R}\cap {\rm Center}\, A$ and $Z:=A^L\cap A^R$ 
the restrictions of $\kappa_{L/R}$ to $Z^{L/R}$ and $Z$ lead to the 
algebra isomorphisms $Z^{L/R}\to\hat Z$ and $Z\to\hat Z^{R/L}$, respectively.
Hence, the {\it hypercenter} $H:=Z\cap {\rm Center}\, A=Z^L\cap Z=Z^R\cap Z$ 
of $A$ is isomorphic to the hypercenter $\hat H$ of $\hat A$ via the 
restriction of $\kappa_L$ or $\kappa_R$ to $H$.

The restrictions of the canonical pairing to $\hat A^{L/R}\times A^{L/R}$
(four possibilities) are non-degenerate. The maps $\Pi^{L/R}$ and 
$\hat\Pi^{L/R}$ ($\bar\Pi^{L/R}$ and $\hat{\bar\Pi}^{R/L}$) are transposed 
to each other
$$\eqalign{
  \langle\varphi ,\Pi^{L/R}(a)\rangle
 &=\langle\hat\Pi^{L/R}(\varphi),a\rangle,\cr
  \langle\varphi ,\bar\Pi^{L/R}(a)\rangle
 &=\langle\hat{\bar\Pi}^{R/L}(\varphi),a\rangle,\cr}\qquad 
  a\in A,\varphi\in\hat A,\eqno(1.7)
$$
(note the switch of $L$ and $R$ in the second equation) 
and obey the identities
$$
  \eqalign{\Pi^R(a^{(1)})\otimes a^{(2)}&=\UN^{(1)}\otimes a\UN^{(2)},\cr
           a^{(1)}\otimes\Pi^L(a^{(2)})&=\UN^{(1)}a\otimes\UN^{(2)},\cr}
  \quad
  \eqalign{\bar\Pi^R(a^{(1)})\otimes a^{(2)}&=\UN^{(1)}\otimes \UN^{(2)}a,\cr
           a^{(1)}\otimes\bar\Pi^L(a^{(2)})&=a\UN^{(1)}\otimes\UN^{(2)}.\cr}
  \qquad a\in A.\eqno(1.8)
$$
due to (1.1c). The space of {\it left/right integrals} $I^{L/R}$ in $A$ is 
defined by
$$
  I^L:=\{ l\in A\,\vert\, al=\Pi^L(a)l, a\in A\},\qquad
  I^R:=\{ r\in A\,\vert\, ra=r\Pi^R(a), a\in A\}.\eqno(1.9)
$$

\smallskip\noindent
{\bf 1.3 Properties of WHAs}\smallskip

Let $A$ be a WHA. The antipode $S$, as in the case of Hopf algebras, 
turns out to be invertible, antimultiplicative, anticomultiplicative and 
leaves the counit invariant: $\varepsilon=\varepsilon\circ S$.
The restriction of the antipode to $A^L$ leads to algebra antiisomorphism
$S\colon A^L\to A^R$, therefore $A^T$ is a subWHA of $A$, moreover,
$$
  \Pi^L\circ S=\Pi^L\circ\Pi^R=S\circ\Pi^R,\qquad   
  \Pi^R\circ S=\Pi^R\circ\Pi^L=S\circ\Pi^L.\eqno(1.10)
$$
The projections (1.3) to left and right subalgebras can be expressed as
$$\eqalign{
  \Pi^L(a)&=a^{(1)}S(a^{(2)}),\cr  
  \bar\Pi^L(a)&=S^{-1}(\Pi^R(a)),\cr}\qquad\eqalign{  
  \Pi^R(a)&=S(a^{(1)})a^{(2)},\cr  
  \bar\Pi^R(a)&=S^{-1}(\Pi^L(a)),\cr}\qquad a\in A.\eqno(1.11)
$$
The first two equations follow from the antipode axioms (1.2a and b). 
The other two can be seen using the aforementioned properties 
of the antipode and the WBA identity $\varepsilon(abc)=
\varepsilon(\Pi^R(a)b\Pi^L(c))$ following from (1.1b) and (1.3). 
The left and right subalgebras become separable $k$-algebras with separating 
idempotents [14, p.182] $q^L=S(\UN^{(1)})\otimes\UN^{(2)}\in A^L\otimes A^L$ 
and $q^R=\UN^{(1)}\otimes S(\UN^{(2)})\in A^R\otimes A^R$, respectively, that
obey
$$\eqalign{
  x^LS(\UN^{(1)})\otimes\UN^{(2)}&=S(\UN^{(1)})\otimes\UN^{(2)}x^L,\qquad
  x^L\in A^L,\cr
  x^R\UN^{(1)}\otimes S(\UN^{(2)})&=\UN^{(1)}\otimes S(\UN^{(2)})x^R,\qquad
  x^R\in A^R\cr}\eqno(1.12)
$$
by definition. The product $q^Lq^R\in A^T\otimes A^T$ is a separating 
idempotent for $A^T$, thus the trivial subalgebra is a separable 
$k$-algebra, too. The separating idempotent $q^{L/R}$ serves as a 
{\it quasibasis} [20, p.6] for the counit: 
$$\eqalign{
  S(\UN^{(1)})\varepsilon(\UN^{(2)}x^L)&=x^L
 =\varepsilon(x^LS(\UN^{(1)}))\UN^{(2)},\qquad x^L\in A^L,\cr
  \UN^{(1)}\varepsilon(S(\UN^{(2)})x^R&=x^R
 =\varepsilon(x^R\UN^{(1)})S(\UN^{(2)}),\qquad
  x^R\in A^R,\cr}\eqno(1.13)
$$
thus the counit is a non-degenerate functional on $A^{L/R}$. The properties 
$S(\UN^{(1)})\UN^{(2)}=\UN$ and $\UN^{(1)}S(\UN^{(2)})=\UN$ of separating
idempotents $q_L$ and $q_R$ ensure that the counit $\varepsilon$ 
is an index $\UN$ functional [20, p.7] on $A^L$ and on $A^R$, respectively. 
Due to the identities (1.5),(1.7),(1.10) and (1.12) the corresponding 
Nakayama automorphisms $\theta_{L/R}\colon A^{L/R}\to A^{L/R}$, which 
are defined by
$$
  \varepsilon(y^{L/R}\theta_{L/R}(x^{L/R})):=\varepsilon(x^{L/R}y^{L/R}),
  \qquad x^{L/R},y^{L/R}\in A^{L/R},\eqno(1.14)
$$
can be given as
$$\eqalign{
  \theta_L(x^L)&=\UN\Sl\hat S^{-1}(\hat\UN\Sl x^L)=S^2(x^L),
                 \qquad x^L\in A^L,\cr
  \theta_R(x^R)&=\hat S(\hat\UN\Sl x^R)\Sr\UN=S^{-2}(x^R),
                 \qquad x^R\in A^R.\cr}\eqno(1.15)  
$$
Hence, $\theta_L$ ($\theta_R$) is the restriction of the square of the 
(inverse of the) antipode to $A^L$ ($A^R$). Since any separable algebra 
admits a non-degenerate (reduced) trace [6, p.165], the counit, being a 
non-degenerate functional on $A^{L/R}$, can be given with the help of the 
corresponding trace as 
$\varepsilon(\cdot)={\rm tr}_{L/R}(t_{L/R}\cdot )$ 
with $t_{L/R}\in A^{L/R}$ invertible. Therefore, the Nakayama automorphisms 
$\theta_{L/R}$ are given by ${\rm ad}\, t_{L/R}$ and $S^2$ is inner on 
$A^{L/R}$, hence, on $A^T$, too.
 
In a WHA a left integral $l\in I^L$ and a right integral $r\in I^R$ obey the
identities 
$$\eqalign{
  l^{(1)}\otimes al^{(2)}&=S(a)l^{(1)}\otimes l^{(2)},\cr
  r^{(1)}a\otimes r^{(2)}&=r^{(1)}\otimes r^{(2)}S(a),\cr}\qquad a\in A,
  \eqno(1.16)
$$
respectively. Moreover, there exist projections $L/R\colon A\to I^{L/R}$ and 
$\bar L/\bar R\colon A\to I^{L/R}$:
$$\eqalign{
           L(a)&:=\hat S^2(\beta_i)\Sr (b_ia),\cr
      \bar L(a)&:=(b_ia)\Sl\hat S^{-2}(\beta_i),\cr}\qquad\eqalign{
           R(a)&:=(ab_i)\Sl\hat S^2(\beta_i),\cr
      \bar R(a)&:=\hat S^{-2}(\beta_i)\Sr (ab_i),\cr}\qquad a\in A
  \eqno(1.17)
$$
where $\{b_i\}\subset A$ and $\{\beta_i\}\subset\hat A$ are dual $k$-bases
with respect to the canonical pairing. They obey the properties
$$\eqalign{
  \langle\hat L/\hat R(\varphi),a\rangle
               &=\langle\varphi,R/L(a)\rangle,\cr
  \langle \hat{\bar L}/\hat{\bar R}(\varphi),a\rangle
               &=\langle\varphi,\bar L/\bar R(a)\rangle,\cr}\qquad 
 a\in A,\varphi\in\hat A,\eqno(1.18)
$$
therefore the restrictions of the canonical pairing to $\hat I^{L/R}\times 
I^{L/R}$ (four possibilities) are non-degenerate.

\medskip\noindent
{\bf 2. Properties of the unit module}\medskip

In this chapter $A$ denotes a WHA over a field $k$.

A {\it left (right) $A$-module} ${_AM}\equiv (M,\mu_L)$ 
(${M_A}\equiv (M,\mu_R)$)
is a $k$-linear space together with the $k$-linear map 
$\mu_L\colon A\otimes M\to M$ ($\mu_R\colon M\otimes A\to M$) satisfying
$$\eqalign{                
            a\cdot(b\cdot m)&=(ab)\cdot m,\cr
            \UN\cdot m&=m,\cr}\qquad\eqalign{
            (m\cdot a)\cdot b&=m\cdot (ab),\cr
             m\cdot\UN&=m,\cr}\qquad m\in M, a,b\in A,
$$
where (and later on) $\mu_L(a\otimes m)\equiv a\cdot m$ and
$\mu_R(m\otimes a)\equiv m\cdot a$. The role of the unit module will be played
by the trivial representation [2, p.400] of $A$:
 
\noindent
{\bf Definition 2.1} {\sl The unit left 
(right) $A$-module ${_AA^L}$ ($A^R_A$) is defined by
$$\eqalign{
  a\cdot x^L&:=\Pi^L(ax^L)=a^{(1)}x^LS(a^{(2)}),\qquad x^L\in A^L, a\in A\cr
  x^R\cdot a&:=\Pi^R(x^Ra)=S(a^{(1)})x^Ra^{(2)},\qquad x^R\in A^R, a\in A.\cr}
  \eqno(2.1)
$$}

We note that these modules need not be one-dimensional as in the case of 
Hopf algebras, they are not even simple in general. Nevertheless, they play 
the role of the unit object in the monoidal category of finite dimensional left (right) $A$-modules. We deal with only the category of left $A$-modules 
since the one-to-one correspondence between left and right $A$-modules 
induced by the antipode, $m\cdot a:=S(a)\cdot m, a\in A,m\in {_AM}$, extends to a categorical isomorphism. 

\smallskip\noindent
{\bf Proposition 2.2} {\sl The category ${\cal L}$ consisting of finite 
dimensional left $A$-modules of a WHA $A$ as objects and left $A$-module
maps as arrows can be endowed with an autonomous (relaxed) monoidal structure:
$({\cal L};\times, {_AA^L}, \{ 1_{K\times M\times N}\},\{ X^L_M\},\{ X^R_M\};
\Cl{\phantom{M}},\Cr{\phantom{M}})$, where $\times$ is the monoidal product,
${_AA^L}$ is the monoidal unit, $\{ 1_{K\times M\times N}\},
\{ X^L_M\},\{ X^R_M\}$ are natural equivalences satisfying the pentagon and 
the triangle identities, while $\Cl{\phantom{M}}$ and $\Cr{\phantom{M}}$ are 
the functors of left and right conjugations, respectively.}

\noindent{\it Proof.} Let us define first the monoidal product: $\times$.    
The product module $_A(M\times N)$ of the modules ${_AM}$ and ${_AN}$ as 
a $k$-linear space is defined to be
$$
  M\times N:=\UN^{(1)}\cdot M\otimes \UN^{(2)}\cdot N\eqno(2.2a)
$$
and the left $A$-module structure on $M\times N$ is given by
$$
  a\cdot\left(m\otimes n\right) 
  :=a^{(1)}\cdot m\otimes a^{(2)}\cdot n, 
  \qquad a\in A,\ m\otimes n\in M\times N,\eqno(2.2b)
$$
where (and later on) we have suppressed possible or necessary summation 
for tensor product elements in product modules. The product on the 
arrows $T_\alpha\colon M_\alpha\to N_\alpha, \alpha=1,2$ is
defined by $T_1\times T_2:=(T_1\otimes T_2)\circ \Delta(\UN)$, i.e. by the 
restriction of the tensor product of the linear maps $T_1$ and $T_2$ to
$M_1\times M_2$. One can easily check that $T_1\times T_2\colon M_1\times M_2
\to N_1\times N_2$ is a left $A$-module map. The given monoidal product is 
associative due to the associativity of the coproduct and property (1.1c) 
of the unit, hence the components $M_1\times (M_2\times M_3)\to 
(M_1\times M_2)\times M_3$ of the natural equivalence responsible for 
associativity in a monoidal category are the identity mappings 
$1_{M_1\times M_2\times M_3}$ in our case. 

The monoidal unit property of the left $A$-module $A^L$ can be seen
by verifying that for any object $M$ the $k$-linear invertible maps 
$X^L_M\colon M\to A^L\times M$ and $X^R_M\colon M\to M\times A^L$
defined by
$$\eqalign{
  X^L_M(m)&:=S(\UN^{(1)})\otimes\UN^{(2)}\cdot m,\cr 
  (X^L_M)^{-1}(x^L\otimes m) 
 &:=x^L\cdot m,\cr}\qquad
  \eqalign{ 
  X^R_M(m)&:=\UN^{(1)}\cdot m\otimes\UN^{(2)},\cr 
  (X^R_M)^{-1}(m\otimes x^L)
 &:=S^{-1}(x^L)\cdot m,\cr}\eqno(2.3)
$$
are left $A$-module maps and the identities
$$\eqalign{
  X^L_NT&=(1_{A^L}\times T)X^L_M\cr
  X^R_NT&=(T\times 1_{A^L})X^R_M\cr}\qquad M,N\in {\rm Obj}\,
  {\cal L},\quad T\colon M\to N,\eqno(2.4)
$$$$
  (X^R_M\times 1_N)(1_M\times (X^L_N)^{-1})=1_{M\times A^L\times N}
  \eqno(2.5)
$$
hold, i.e. $X^L=\{ X^L_M\}$ and $X^R=\{ X^R_M\}$ are natural equivalences 
satisfying the triangle identity.

An autonomous category [21] contains both left and right conjugation functors
by definition. The left conjugate $\Cl{M}$ of an object $M$ in 
${\cal L}$ is given by the $k$-dual $\hat M:={\rm Hom}_k(M,k)$ as a $k$-linear
space. The left $A$-module structure $\Cl{M}\equiv(\hat M,\Cl{\mu}_L)$ 
is defined to be 
$$
  \langle a\cdot\hat m,m\rangle := \langle \hat m,S(a)\cdot m\rangle,
  \qquad a\in A, m\in M,\hat m\in\hat M,\eqno(2.6)
$$
where $\langle\, ,\, \rangle$ is the $k$-valued canonical bilinear paring 
on the cartesian product of $\hat M$ and $M$. Dual bases with respect to this 
pairing will be denoted by $\{\hat m_i\}_i\subset\hat M$ and 
$\{ m_i\}_i\subset M$.
Due to the definition (2.6) of the left $A$-module $\Cl{M}$ we have 
$$\eqalignno{
  m_i\otimes a\cdot \hat m_i &=S(a)\cdot m_i\otimes\hat m_i, 
     \quad a\in A&(2.7a)\cr
  m_i\otimes \hat m_i&=\UN\cdot m_i\otimes\hat m_i 
    =\UN^{(1)}S(\UN^{(2)})\cdot m_i\otimes\hat m_i\cr 
     &=\UN^{(1)}\cdot m_i\otimes \UN^{(2)}\cdot\hat m_i\in  
     M\times \Cl{M},&(2.7b)\cr}
$$
where (and later on) we omit summation symbol for the sum of tensor product 
of dual basis elements.

The arrow family of left evaluation and coevaluation maps 
$E^l_M\colon \Cl{M}\times M\to A^L$ and  
$C^l_M\colon A^L\to M\times\Cl{M}$, respectively, are defined 
to be
$$\eqalign{
  E^l_M(\hat m\otimes m)&:=\UN^{(2)}\langle \hat m,\UN^{(1)}\cdot m\rangle,
  \qquad\hat m\otimes m\in\Cl{M}\times M\cr
  C^l_M(x^L)&:=x^L\cdot m_i\otimes\hat m_i,\qquad x^L\in A^L.\cr}\eqno(2.8)
$$
They are left $A$-module maps
$$\eqalign{
   E^l_M(a\cdot(\hat m\otimes m))
  &=\UN^{(2)}\langle a^{(1)}\cdot\hat m,\UN^{(1)}a^{(2)}\cdot m\rangle
   =\UN^{(2)}\langle \hat m,S(a^{(1)})\UN^{(1)}a^{(2)}\cdot m\rangle\cr
  &=\Pi^L(a^{(3)})\langle \hat m,S(a^{(1)})a^{(2)}\cdot m\rangle
   =\Pi^L(a^{(2)})\langle \hat m,\Pi^R(a^{(1)})\cdot m\rangle\cr
  &=\Pi^L(a\UN^{(2)})\langle \hat m,\UN^{(1)}\cdot m\rangle
   =a\cdot E^l_M(\hat m\otimes m),\cr
   C^l_M(a\cdot x^L)&=a^{(1)}x^LS(a^{(2)})\cdot m_i\otimes \hat m_i
   =a^{(1)}x^L\cdot m_i\otimes a^{(2)}\cdot\hat m_i\cr 
  &=a\cdot C^l_M(x^L)\cr}\eqno(2.9)
$$
due to the identities (1.8) and (2.7a) and they satisfy the left rigidity 
identities [21]
$$\eqalignno{
  (X^R_M)^{-1}(1_M\times E_M^l)(C_M^l\times 1_M)X^L_M(m)
 &:=S^{-1}(\UN^{(2')})S(\UN^{(1)})\cdot m_i
             \langle\hat m_i,\UN^{(1')}\UN^{(2)}\cdot m\rangle\cr
 =S^{-1}(\UN^{(2')})S(\UN^{(1)})\UN^{(1')}\UN^{(2)}\cdot m&=m,\qquad m\in M,
  &(2.10a)\cr
  (X^L_{\CL{M}})^{-1}(E_M^l\times 1_{\CL{M}})(1_{\CL{M}}\times C_M^l)
   X^R_{\CL{M}}(\hat m)
 &:=\langle\UN^{(1)}\cdot\hat m,\UN^{(1')}\UN^{(2)}\cdot m_i\rangle     
     \UN^{(2')}\cdot\hat m_i\cr
  =\UN^{(2')}S^{-1}(\UN^{(1')}\UN^{(2)})\UN^{(1)}\cdot\hat m
 &=\hat m,\qquad \hat m\in\Cl{M}&(2.10b)\cr}
$$
for any $M\in {\rm Obj}\, {\cal L}$. Thus defining the left conjugated 
arrow $\Cl{T}\colon\Cl{N}\to\Cl{M}$ of $T\colon M\to N$ by
$$
  \Cl{T}:=(X^L_{\CL{M}})^{-1}(E_N^l\times 1_{\CL{M}})
         (1_{\CL{N}}\times T\times 1_{\CL{M}})
         (1_{\CL{N}}\times C_M^l)X^R_{\CL{N}}\eqno(2.11)
$$
one arrives at the antimonoidal contravariant left conjugation functor
$\Cl{\phantom{M}}\colon{\cal L}\to{\cal L}$ [21].

Similarly, the right conjugate $\Cr{M}$ of an object $M$ in 
${\cal L}$ is the $k$-linear space $\hat M$ equipped with the
left $A$-module structure $\Cr{M}\equiv(\hat M,\Cr{\mu}_L)$ 
$$
  \langle a\cdot\hat m,m\rangle := \langle\hat m,S^{-1}(a)\cdot m\rangle,
  \qquad a\in A, m\in M,\hat m\in\hat M\eqno(2.12)
$$
implying 
$$\eqalignno{
  \hat m_i\otimes a\cdot m_i &=S(a)\cdot\hat m_i\otimes m_i, 
     \quad a\in A&(2.13a)\cr
  \hat m_i\otimes m_i &=\UN\cdot\hat m_i\otimes m_i 
    =\UN^{(1)}S(\UN^{(2)})\cdot\hat m_i\otimes m_i\cr 
    &=\UN^{(1)}\cdot\hat m_i\otimes \UN^{(2)}\cdot m_i\in  
     \Cr{M}\times M.&(2.13b)\cr}
$$
The arrow family of right evaluation and coevaluation maps 
$E^r_M\colon M\times \Cr{M}\to A^L$ and  
$C^r_M\colon A^L\to \Cr{M}\times M$, respectively, are defined 
to be
$$\eqalign{
  E^r_M(m\otimes\hat m)&:=\UN^{(2)}\langle\UN^{(1)}\cdot\hat m, m\rangle,
  \qquad m\otimes\hat m\in M\times \Cr{M},\cr
  C^r_M(x^L)&:=x^L\cdot\hat m_i\otimes m_i,\qquad x^L\in A^L.\cr}\eqno(2.14)
$$
As in the previous case, one proves that they are left $A$-module maps
satisfying the right rigidity identities [21]
$$\eqalign{
  (X^L_M)^{-1}(E_M^r\times 1_M)(1_M\times C_M^r)X^R_M&=1_M,\cr
  (X^R_{\CR{M}})^{-1}(1_{\CR{M}}\times E_M^r)(C_M^r\times 1_{\CR{M}})
   X^L_{\CR{M}}&=1_{\CR{M}}.\cr}\eqno(2.15)
$$
Hence, defining the right conjugated arrow $\Cr{T}\colon\Cr{N}\to\Cr{M}$ of
$T\colon M\to N$ by
$$
  \Cr{T}:=(X^R_{\CR{M}})^{-1}(1_{\CR{M}}\times E_N^r)
         (1_{\CR{M}}\times T\times 1_{\CR{N}})
         (C_M^r\times 1_{\CR{N}})X^L_{\CR{N}}\eqno(2.16)
$$
one arrives at the antimonoidal contravariant right conjugation functor
$\Cr{\phantom{M}}\colon{\cal L}\to{\cal L}$.\qed

\smallskip
In order to prove semisimplicity of the unit module, we show that the trivial 
subWHA is not only semisimple but also cosemisimple:

\medskip\noindent
{\bf Lemma 2.3} {\sl The trivial weak Hopf subalgebra $A^T\subset A$ is a 
sum of simple subcoalgebras, i.e. $A^T$ is contained in the coradical 
$C_0$ of $A$.}

\smallskip\noindent{\it Proof.} First we decompose the WHA $A^T$ into 
a direct sum of subWHAs. 

The intersection $Z:=A^L\cap A^R$ is in the center of the separable algebra 
$A^T$, because the unital coideal subalgebras $A^L$ and $A^R$ that generate 
$A^T$ commute with each other. The WHA identity (1.10) implies $z=S(z)$ for all 
$z\in Z$. Hence, $Z$ is a unital, pointwise $S$-invariant subalgebra of the 
$k$-algebra ${\rm Center}\, A^T$ and one can write $A^T$ 
as a tensor product algebra $A^T\simeq A^L\otimes_Z A^R$. 
Let $\{z_\alpha\}_\alpha$ be the set of primitive orthogonal
idempotents in $Z$. They are central idempotents in $A^T$; thus, 
$A^T=\oplus_\alpha A^T_\alpha,\ A^T_\alpha:=A^Tz_\alpha$ is an
ideal decomposition of the algebra $A^T$. It is also a WHA decomposition:
first, $S(A^T_\alpha)=A^T_\alpha$ since $Z$ is pointwise 
$S$-invariant, and second, 
$\Delta(A^T_\alpha)\subset A^T_\alpha\otimes A^T_\alpha$, because
$$
  \Delta(xz_\alpha)=\Delta(xz_\alpha z_\alpha)
  =\Delta(x)(z_\alpha\otimes\UN)(\UN\otimes z_\alpha)
  =\Delta(x)(z_\alpha\otimes z_\alpha),\quad x\in A^T,\eqno(2.17)
$$
due to $z_\alpha\in A^L\cap A^R$ and due to coproduct property (1.4) of 
elements in $A^L$ and in $A^R$. 

This WHA decomposition implies that $(A^T_\alpha)^X=A^X_\alpha$ with
$X=L,R,T$ and that the WHA $A^T_\alpha$ has the tensor product 
algebra structure $A^T_\alpha\simeq A^L_\alpha\otimes_{Z_\alpha} A^R_\alpha$ 
with unit $\UN_\alpha=z_\alpha$. The algebra $Z_\alpha:=Zz_\alpha=A^L_\alpha
\cap A^R_\alpha$ is an Abelian division algebra over the ground 
field $k$, that is $Z_\alpha$ is a subfield in the center of the separable 
algebra $A^T_\alpha$, hence 
$Z_\alpha$ is a finite separable field extension of $k$ [14, p.191].

Now we prove that $\widehat{A^T_\alpha}$, the dual of the WHA $A^T_\alpha$, is 
isomorphic to the simple $k$-algebra $M_{n_\alpha}(Z_\alpha)$, where 
$n_\alpha={\rm dim}_{Z_\alpha}\, A^R_\alpha$, i.e. $A^T_\alpha$
is simple as a $k$-coalgebra. We stress that the inclusion 
$(\widehat{A^T_\alpha})^T\subset\widehat{A^T_\alpha}$ is proper in general.
Therefore, simplicity of $\widehat{A^T_\alpha}$ as an algebra is
a `non-trivial' property in the sense that it goes for a WHA which is not 
trivial, i.e., not generated by the canonical coideal subalgebras
$(\widehat{A^T_\alpha})^L$ and $(\widehat{A^T_\alpha})^R$.

Consider the cyclic left $\widehat{A^T_\alpha}$-module $(A^T_\alpha)^R=
A^R_\alpha$ with the Sweedler action $\varphi\cdot x^R:=\varphi\Sr x^R;
\varphi\in\widehat{A^T_\alpha},x^R\in A^R_\alpha$. It is just the trivial 
representation [2, p.401] of the WHA $\widehat{A^T_\alpha}$; hence,
its endomorphism ring ${\rm End}_{\widehat{A^T_\alpha}} A^R_\alpha$ is given 
by $\hat Z^R_\alpha\Sr$ [2, p.402], where $\hat Z^R_\alpha:=
{\rm Center}\widehat{A^T_\alpha}\cap (\widehat{A^T_\alpha})^R$.
However, the maps in $\hat Z^R_\alpha\Sr$ are just  
multiplications by elements of $Z_\alpha$ due to the statements after (1.6), 
i.e., ${\rm End}_{\widehat{A^T_\alpha}} A^R_\alpha =Z_\alpha$. 
In addition, $A^R_\alpha$ is a faithful module in this case, because 
$\varphi\Sr A^R_\alpha=0$ implies $\varphi=0$ due to
$$\eqalign{
  0&=\varepsilon((\varphi\Sr A^R_\alpha)A^L_\alpha)
    =\varepsilon(\UN_\alpha^{(1)}A^L_\alpha)
     \langle\varphi,A^R_\alpha\UN_\alpha^{(2)}\rangle
    =\langle\varphi,A^R_\alpha\Pi^L(A^L_\alpha)\rangle
    =\langle\varphi, A^R_\alpha A^L_\alpha\rangle\cr
    &=\langle\varphi,A^T_\alpha\rangle,\cr}\eqno(2.18)
$$
where we used (1.4) and (1.3). Hence, $\widehat{A^T_\alpha}$ is isomorphic 
to the unital subalgebra $\phi(\widehat{A^T_\alpha}):=\widehat{A^T_\alpha}\Sr$ 
of ${\rm End}_k A^R_\alpha$. Therefore, the relations
$\widehat{A^T_\alpha}\simeq\phi(\widehat{A^T_\alpha})
\subset {\rm BiEnd}_{\widehat{A^T_\alpha}} A^R_\alpha
={\rm End}_{Z_\alpha} A^R_\alpha\simeq M_{n_\alpha}(Z_\alpha)$
together with the equality 
$$\eqalign{
    {\rm dim}_k \widehat{A^T_\alpha}&={\rm dim}_k A^T_\alpha
   ={\rm dim}_k (A^L_\alpha\otimes_{Z_\alpha} A^R_\alpha)
   =({\rm dim}_{Z_\alpha} A^R_\alpha)^2 {\rm dim}_k Z_\alpha\cr
  &=n_\alpha^2{\rm dim}_k Z_\alpha={\rm dim}_k M_{n_\alpha}(Z_\alpha)\cr}
   \eqno(2.19) 
$$
concerning the $k$-dimensions imply the $k$-algebra isomorphism
$\widehat{A^T_\alpha}\simeq M_{n_\alpha}(Z_\alpha)$.\qed   

\medskip\noindent
{\bf Theorem 2.4} {\sl The unit left $A$-module ${_AA^L}$ 
is the direct sum of simple submodules
$$
  {_AA^L}=\bigoplus_p {_AA^L_p},\qquad A^L_p:=A^Lz^L_p,\eqno(2.20)
$$
where $\{ z_p^L\}_p$ is the set of primitive orthogonal idempotents in 
$Z^L:=A^L\cap {\rm Center}\, A$.}

\smallskip
\noindent{\it Proof.} Let $N$ be the radical of $A$. Since $N$ is an ideal
in $A$, we have $N\cdot A^L:=\Pi^L(NA^L)\subset \Pi^L(N)$. Due to the 
identity $N=(\hat C_0)^\perp:=\{ a\in A\,\vert\, \langle \hat C_0, 
a\rangle=0\}$ [17, p.183], where $\hat C_0$ is the coradical of 
the dual weak Hopf algebra $\hat A$, the previous Lemma leads to the 
containment
$N=(\hat C_0)^\perp\subset(\hat A^T)^\perp\subset(\hat A^L)^\perp$>. Hence,
using (1.7) the canonical pairing gives rise to
$$
  \langle\hat A, \Pi^L(N)\rangle=\langle \hat\Pi^L(\hat A), N\rangle
 =\langle\hat A^L,N\rangle=0,\eqno(2.21) 
$$
i.e. $\Pi^L(N)=0$. Therefore, the radical of $A$ is in the annihilator ideal
of the left module ${_AA^L}$, that is ${_AA^L}$ is semisimple [6, p.118].

The endomorphism ring for the unit module is given by 
${\rm End}\, {_AA^L}=Z^L\cdot$ [2, p.402], that is by the restriction
of the $A$-action to the subalgebra $Z^L$. Since the unit module is
a free, hence faithful $A^L$-module, it is also faithful as a $Z^L$-module.
Now, the direct sum decomposition (2.20) is clear and 
${\rm End}\, {_AA^L_p}=Z^Lz_p^L\cdot$. But $Z^Lz_p^L$ is a field, therefore 
${_AA^L_p}$ is indecomposable [6, p.121]. Together 
with semisimplicity this leads to simplicity of the direct summands 
${_AA^L_p}$.\qed

The analogous result holds for the unit right $A$-module:

\smallskip\noindent
{\bf Remark 2.5} {\sl The unit right $A$-module 
$A^R_A$ given in Def. 2.1 is the direct sum of simple submodules
$$
  {A^R_A}=\bigoplus_p {A^R_{pA}},\qquad A^R_p:=A^Rz_p^R\eqno(2.22)
$$
where $\{ z_p^R:=S(z_p^L)\}_p$ is the set of primitive orthogonal 
idempotents in $Z^R=S(Z^L)$.}\smallskip

We have seen that the simple submodules of the unit left (right) $A$-module 
are labelled by primitive idempotents in $Z^L$ ($Z^R$). Although a generic 
$A$-module does not need to be semisimple, it is always a direct sum of 
submodules labelled by pairs of primitive orthogonal idempotents in the 
cartesian product $Z^L\times Z^R$. Indeed, the product of primitive orthogonal 
idempotents in $Z^L$ and $Z^R$ gives rise to a decomposition of the unit 
$\UN=\sum_{p,q}z_p^Lz_q^R\equiv\sum_{p,q}z_p^LS(z_q^L)$ in 
$Z^L\vee Z^R\subset A^T\cap{\rm Center}\, A$.\footnote{${}^1$}{Note that 
$S^2(z_p^{L/R})=z_p^{L/R}$, because $S^2$ is inner on $A^{L/R}$ and the 
idempotents are central.} 
Since $Z^L\vee Z^R\simeq Z^L\otimes_H Z^R$ certain products $z_p^Lz_q^R$ 
can be identically zero due to the presence of the hypercenter 
$H:=Z^L\cap Z^R=A^L\cap A^R\cap {\rm Center}\, A$ in $A^T$. 
If $z_p^Lz_q^R\not= 0$ we refer to $(p,q)$ as an 
{\it admissible pair}. Hence, the non-zero summands are labelled by 
admissible pairs in the decomposition of the unit, which induces a direct sum 
decomposition of every $A$-module ${_AM}$ into submodules as 
$M=\oplus_{(p,q)}M_{(p,q)}$, where $M_{(p,q)}:=z_p^Lz_q^R\cdot M$. 
We will call the left $A$-module ${_AM}$ a {\it member
of the $(p,q)$ class} and write ${_AM}_{(p,q)}$ if $(p,q)$ is an admissible 
pair and 
$$
  z_{p'}^Lz_{q'}^R\cdot M=\delta_{pp'}\delta_{qq'}M,\eqno(2.23) 
$$
for all product idempotents $z_{p'}^Lz_{q'}^R\in Z^L\vee Z^R$.
Clearly, the simple submodule ${_AA}^L_p$ of the 
unit module ${_AA}^L$ is in the `diagonal' class $(p,p)$, because
$$
  z_q^Lz_r^R\cdot xz_p^L:=z_q^Lxz_p^LS(z_r^R)=z_q^Lxz_p^Lz_r^L=
  \delta_{p,q}\delta_{p,r}xz_p^L, \qquad x\in A^L.\eqno(2.24)
$$
The next Lemma shows that the simple submodules of the unit module
${_AA}^L$ obey a kind of minimality property in the corresponding class of 
left $A$-modules.

\medskip\noindent
{\bf Lemma 2.6} {\sl i) The $k$-dimension $\vert M_{(p,q)}\vert$ of a 
nonzero left $A$-module ${_AM}_{(p,q)}$ in the $(p,q)$ class obeys the 
inequality
$$
  \vert M_{(p,q)}\vert\geq \max\{ \vert A^L_p\vert, 
  \vert A^L_q\vert\} .\eqno(2.25)
$$

ii) The restriction of $A$ to the subalgebras $A^L_p$ and $A^R_q$ makes 
${_AM}_{(p,q)}$ a faithful left $A^L_p$- and $A^R_q$-modules, respectively.}
\smallskip

\noindent{\it Proof.} In the following first we prove that the left $A$-modules
$M_{(p_1,q_1)}$ and $N_{(p_2,q_2)}$ should obey the matching condition 
$q_1=p_2$ in order to get a nonzero product module 
$M_{(p_1,q_1)}\times N_{(p_2,q_2)}$. Then writing a left $A$-module 
$M_{(p,q)}$ as a product with the unit module and using this matching 
condition, the emerging tensor product space can be given as a sum 
of subspaces with respect to a basis of the corresponding simple submodule of 
the unit module. We will use Theorem 2.4 and Remark 2.5 to prove that 
$M_{(p,q)}$ is a faithful $A^L_p$- and $A^R_q$-module and then the 
estimation of the $k$-dimension of $M_{(p,q)}$ will follow.
    
Using property (1.12) of the separating idempotent of 
$A^L$ and the decomposition of the unit into primitive orthogonal 
idempotens in $Z^L$, one obtains
$$
  \Delta(\UN)=\sum_r S^{-1}(S(\UN^{(1)}))\otimes \UN^{(2)}z_r^L
             =\sum_r (z_r^R\otimes z_r^L)\Delta(\UN).\eqno(2.26)
$$
Therefore, for any two left $A$-modules $M,N$ within a certain class we have 
$$\eqalign{ 
  z_p^Lz_q^R\cdot M_{(p_1,q_1)}\times N_{(p_2,q_2)}
 &=z_p^Lz_q^R\cdot (\Delta(\UN)(M_{(p_1,q_1)}\otimes N_{(p_2,q_2)}))\cr
  &=\sum_r z_p^Lz_r^R\UN^{(1)}\cdot M_{(p_1,q_1)}\otimes 
   z_r^Lz_q^R\UN^{(2)}\cdot N_{(p_2,q_2)},\cr}\eqno(2.27)
$$
implying 
$$
  {_AM}_{(p_1,q_1)}\times {_AN}_{(p_2,q_2)}=\left\{\eqalign{
  &0,\qquad q_1\not= p_2,\cr
  &{_A(M\times N)}_{(p_1,q_2)},\qquad q_1=p_2.\cr}\right.\eqno(2.28)
$$
The separating idempotent of $A^L$ is a quasibasis for the counit due to
(1.13), hence, it has the expression 
$S(\UN^{(1)})\otimes \UN^{(2)}=\sum_i f_i\otimes e_i$, where
$\{ e_i\}_i,\{ f_i\}_i\subset A^L$ are dual bases with respect to the 
counit: $\varepsilon(e_if_j)=\delta_{i,j}$. Choosing a basis 
$\{ e_i\}_i=\cup_p\{ e_i\}_{i\in p}$ that respects the direct sum decomposition 
$A^L=\oplus_p A^L_p$, i.e. $\{ e_i\}_{i\in p}\subset A^L_p$, 
$f_i\otimes e_i\in A^L_p\otimes A^L_p, i\in p$ follows.
Since $X^L_M$ and $X^R_M$ defined in (2.3) are left $A$-module isomorphisms
$$\eqalign{
  \vert M_{(p,q)}\vert 
 &= \vert A^L\times M_{(p,q)}\vert  
  = \vert S(\UN^{(1)})\otimes \UN^{(2)}z^L_p\cdot M_{(p,q)}\vert\cr
 &= \vert \sum_{i\in p} f_i\otimes e_i\cdot M_{(p,q)}\vert
  =\sum_{i\in p}\vert e_i\cdot M_{(p,q)}\vert,\cr}\eqno(2.29a)
$$ 
$$\eqalign{
  \vert M_{(p,q)}\vert 
 &= \vert M_{(p,q)}\times A^L\vert  
  = \vert \UN^{(1)}z^R_q\cdot M_{(p,q)}\otimes \UN^{(2)}\vert\cr 
 &= \vert \sum_{i\in q} S^{-1}(f_i)\cdot M_{(p,q)}\otimes e_i\vert
  =\sum_{i\in q}\vert S^{-1}(f_i)\cdot M_{(p,q)}\vert\cr}\eqno(2.29b)
$$
for any left $A$-module $M_{(p,q)}$ in the $(p,q)$ class.
Hence, if we prove that $M_{(p,q)}$ is a faithful left $A^L_p$- and 
$A^R_q$-module, i.e. $x^L_p\cdot M_{(p,q)}$ and $x^R_q\cdot M_{(p,q)}$ 
are nonzero linear subspaces of $M_{(p,q)}$ for all non-zero elements 
$x^L_p\in A^L_p$ and $x^R_q\in A^R_q$, respectively, then we are done, 
because a nonzero linear subspace is at least one dimensional and 
$\vert A^R_q\vert =\vert S(A^R_q)\vert=\vert A^L_q\vert$ 
due to the invertibility of the antipode $S$. 

Let us suppose that $0\not= x^L_p\in A^L_p$ is in the annihilator ideal 
of ${_AM}_{(p,q)}$. Since in the simple module ${_AA}^L_p$ 
every non-zero element is cyclic, $A^L_p=\{ a\cdot x^L_p := 
a^{(1)}x^L_pS(a^{(2)})\,\vert\, a\in A\}$ should also be contained in the 
annihilator ideal of ${_AM}_{(p,q)}$.
But this contradicts the assumption that ${_AM}_{(p,q)}$ is a nonzero 
module in the $(p,q)$ class. Since the module 
$A^R_{qA}$ is simple (see Remark 2.5), one has $A^R_{qA}=\{ x^R_q\cdot a:= 
S(a^{(1)})x^R_qa^{(2)}\,\vert\, a\in A\}$ for any non-zero $x^R_q\in A^R_q$.
Hence, the assumption that a non-zero element of $A^R_q$ is in the annihilator 
ideal of ${_AM}_{(p,q)}$ leads to the contradiction as before.\qed

\medskip\noindent
{\bf Corollary 2.7} {\sl Given ${_AM}$ let ${_{A^L}M}$ and ${_{A^R}M}$ denote
the $A^L$- and $A^R$-module, respectively, defined by restriction of the 
$A$-module structure to these subalgebras. 
If $\,{\rm End}_{A^R} M=A^L\cdot\simeq A^L$ then ${_{A^L}M}$ and ${_{A^R}M}$ 
are free rank one $A^L$- and $A^R$-modules, respectively.}
\smallskip

\noindent{\it Proof.} Repeating the argument in [2, p.417], one obtains an 
upper bound for the $k$-dimension $\vert M\vert$ of the module $M$: 
being separable, $A^R$ is semisimple; hence, by the Wedderburn structure theorem
$A^R=\oplus_i A^R_i\simeq \oplus_i M_{n_i}(D_i)$, where $D_i$ is a division
algebra corresponding to the simple ideal $A^R_i$. Let $m_i$ denote the 
multiplicity of the simple $A^R_i$-submodules in the semisimple module
${_{A^R}M}$. Then ${\rm End}_{A^R} M\simeq\oplus_i M_{m_i}(D_i)$, which 
is isomorphic to $A^L$ by assumption. Hence, as a right action on $M$, it is 
antiisomorphic to $A^L$, i.e., isomorphic to $A^R$. This is possible only if 
there is a permutation $\sigma$ of simple ideals of $A^R$ such that 
$n_{\sigma(i)}=m_i$ and $D_{\sigma(i)}=D_i$. Therefore  
$M=\oplus_i{\rm Mat}(n_i\times\sigma(n_i),D_i)$ as an $A^R$-bimodule and  
$\vert M\vert=\sum_i \vert D_i\vert n_in_{\sigma(i)}$. The upper bound   
$\vert A^R\vert=\sum_i\vert D_i\vert n_i^2$ for $\vert M\vert$ follows from 
the Cauchy--Schwarz inequality. 

However, the $A^R$-bimodule structure of $M$ implies that $M$ is a faithful 
left $A^R$-module, hence, a faithful left $Z^R$-module. Therefore, the previous 
Lemma leads to the opposite estimation: $\vert M\vert\geq\vert A^R\vert$. Thus
$\sum_i \vert D_i\vert n_in_{\sigma(i)}=\vert M\vert
=\vert A^R\vert=\sum_i\vert D_i\vert n_i^2$, which is possible only if 
$n_{\sigma(i)}=n_i$.
But in this case ${_{A^R}M}$ and ${_{A^L}M}$ are isomorphic to the 
left regular $A^R$- and $A^L$-module, respectively, that is ${_AM}$ 
becomes a free rank one $A^R$- and $A^L$-module by restricting the 
$A$-action to these subalgebras.\qed 

\medskip\noindent
{\bf Lemma 2.8} {\sl If the isomorphism $A^L\simeq M\times\Cl{M}$ 
($A^L\simeq\Cr{M}\times M$) of left $A$-modules holds for a left $A$-module 
${_AM}$, where $A^L$ denotes the unit left $A$-module, then 
${\rm End}_{A^R} M=A^L\cdot\simeq A^L$ 
(${\rm End}_{A^L} M=A^R\cdot\simeq A^R$) as $k$-algebras.}
\smallskip

\noindent{\it Proof.} Since $\Cl{M}:={\rm Hom}_k(M,k)$ as a 
$k$-linear space, one can realize ${\rm End}\, {_kM}$ by 
$M\otimes\Cl{M}$ as $(\sum_a m_a\otimes\hat m_a)(m):=
\sum_a m_a\langle\hat m_a,m\rangle, m\in M$. The subalgebra 
${\rm End}\, {_{A^R}M}\subset {\rm End}\, {_kM}$ is given by 
$M\times\Cl{M}$: 
if $f=\sum_a m_a\otimes\hat m_a\in {\rm End}\, {_{A^R}M}$
then using (2.6) we get
$$\eqalign{
  f(m)
 &=f(\UN\cdot m)
  =f(\UN^{(1)}S(\UN^{(2)})\cdot m)
  =\UN^{(1)}\cdot f(S(\UN^{(2)})\cdot m)\cr
 &=\sum_a\UN^{(1)}\cdot m_a\langle\hat m_a,
   S(\UN^{(2)})\cdot m\rangle
  =\sum_a\UN^{(1)}\cdot m_a\langle\UN^{(2)}\cdot \hat m_a, m\rangle\cr
 &=(\sum_a\UN^{(1)}\cdot m_a\otimes\UN^{(2)}\cdot\hat m_a)(m),
   \quad m\in M,\cr}\eqno(2.30)
$$
that is ${\rm End}\, {_{A^R} M}\subset M\times\Cl{M}$. Choosing
$f=\sum_a m_a\otimes\hat m_a\in M\times\Cl{M}$ and $x^R\in
A^R$ 
$$\eqalign{
  f(x^R\cdot m)
 &\equiv(\sum_a\UN^{(1)}\cdot m_a\otimes\UN^{(2)}\cdot\hat m_a)
   (x^R\cdot m)
  =\sum_a\UN^{(1)}\cdot m_a\langle \UN^{(2)}\cdot\hat m_a,
   x^R\cdot m\rangle\cr
 &=\sum_a\UN^{(1)}\cdot m_a\langle S^{-1}(x^R)\UN^{(2)}\cdot\hat m_a,m\rangle
  =\sum_a x^R\UN^{(1)}\cdot m_a\langle\UN^{(2)}\cdot\hat m_a,m\rangle\cr
 &=x^R\cdot(\sum_a m_a\otimes\hat m_a)(m)=x^R\cdot f(m),
   \qquad m\in M,\cr}\eqno(2.31)
$$
leads to the opposite containment, hence ${\rm End}\, {_{A^R}M}=M\times\Cl{M}$.

The structure of the subalgebra $M\times\Cl{M}$ of $M\otimes\Cl{M}$ can be 
obtained as follows. The unit of ${\rm End}\, {_kM}$ is given by
$1_M=\sum_i m_i\otimes \hat m_i\in M\otimes\Cl{M}$, where $\{m_i\}_i\subset M$ 
and $\{\hat m_i\}_i\subset \Cl{M}$ are dual bases. $1_M$ is in 
$M\times\Cl{M}$ due to (2.7b). Let $U\colon A^L\to M\times\Cl{M}$ be 
the required left $A$-module isomorphism. Hence, there exists $u^L\in A^L$ such
that $1_M=U(u^L)$ and 
$$\eqalign{
  {\rm End}\, {_{A^R}M}&=M\times\Cl{M}=U(A^L)
  =U(A^L\cdot\UN)=A^L\cdot U(\UN)
  =:A^L\cdot(\sum_a m_a\otimes\hat m_a)\cr
 &=\sum_i A^L\cdot m_a\otimes\hat m_a
  =(A^L\cdot)\circ U(\UN).\cr}\eqno(2.32)
$$
using (1.4) in the sixth equality. Since the unit $A$-module ${_AA^L}$ becomes 
a free, hence faithful left $A^L$-module by restriction and since $U(\UN)\in
{\rm End}\, {_{A^R} M}$ is invertible due to 
$1_M=U(u^L)=(u^L\cdot)\circ U(\UN)$, 
(2.32) leads to ${\rm End}\, {_{A^R} M}=A^L\cdot\simeq A^L$.

The proof of the statement involving the right dual $\Cr{M}$ is similar.\qed

\bigskip\noindent
{\bf 3. Hopf modules in weak Hopf algebras}\medskip

Besides $A$-modules we need the notion of weak Hopf modules of a WBA $A$ 
[2, p.407]. First, a {\it left (right) $A$-comodule} is a pair $^AM\equiv
(M,\delta_L)$ ($M^A\equiv(M,\delta_R)$) consisting of a finite dimensional
$k$-linear space and a $k$-linear map $\delta_L\colon M\to A\otimes M$ 
($\delta_R\colon M\to M\otimes A$) called the coaction that obeys
$$\eqalign{
    (id_A\otimes\delta_L)\circ\delta_L&=(\Delta\otimes id_M)\circ\delta_L,\cr
    (\varepsilon\otimes id_M)\circ\delta_L&=id_M,\cr}\qquad\eqalign{
    (\delta_R\otimes id_A)\circ\delta_R&=(id_M\otimes\Delta)\circ\delta_R,\cr
    (id_M\otimes\varepsilon)\circ\delta_R&=id_M.\cr}\eqno(3.1)
$$
They incorporate only the coalgebra properties of $A$. In the following we will use the notations $\delta_L(m)\equiv m_{-1}\otimes m_0$ and 
$\delta_R(m)\equiv m_0\otimes m_1$. Lower and upper $A$-indices will indicate
$A$-modules and $A$-comodules, respectively.

The {\it weak Hopf modules} (WHM) ${M^A_A}, {_AM^A}, {^A_AM}, {^AM_A}$ 
of a WBA $A$ are $A$-modules and $A$-comodules simultaneously together 
with a compatibility condition restricting the comodule map to be an 
$A$-module map, e.g.
$$\eqalign{ 
  M^A_A&\equiv(M,\mu_R,\delta_R):\quad 
     (m\cdot a)_0\otimes (m\cdot a)_1=m_0\cdot a^{(1)}\otimes m_1a^{(2)},\cr
  _AM^A&\equiv(M,\mu_L,\delta_R):\quad 
     (a\cdot m)_0\otimes (a\cdot m)_1=a^{(1)}\cdot m_0\otimes a^{(2)}m_1,\cr}
  \qquad a\in A, m\in M.\eqno(3.2)
$$
As a consequence of these identities WHMs obey a kind of non-degeneracy 
property
$$\eqalign{
  m&=m_0\cdot\Pi^R(m_1),\quad m\in {M_A^A};\cr
  m&=m_0\cdot\bar\Pi^L(m_{-1}),\quad m\in {^AM_A};}\qquad\eqalign{  
  m&=\bar\Pi^R(m_1)\cdot m_0,\quad m\in {_AM^A};\cr
  m&=\Pi^L(m_{-1})\cdot m_0,\quad m\in {^A_AM}.\cr}\eqno(3.3)          
$$
We call ${_AM^A_A}, {^AM^A_A}, {^A_AM_A}, {^A_AM^A}, {_A^AM_A^A}$ 
{\it multiple weak Hopf modules} if they are pairwise WHMs of the WBA $A$ 
in the possible $A$-indices and if the different module or comodule maps 
commute, i.e., they are bimodules or bicomodules.
The {\it invariants} and {\it coinvariants} of left/right $A$-modules and 
left/right $A$-comodules, respectively, are defined to be
$$\eqalign{
  I({_AM})&:=\{ m\in M\,\vert\, a\cdot m=\pi^L(a)\cdot m, a\in A\},\cr
  I({M_A})&:=\{ m\in M\,\vert\, m\cdot a=m\cdot \pi^R(a), a\in A\},\cr}
  \ \eqalign{
  C({^AM})&:=\{ m\in M\,\vert\,\delta_L(m)\in A^R\otimes M\},\cr
  C({M^A})&:=\{ m\in M\,\vert\, \delta_R(m)\in M\otimes A^L\}.\cr}
  \eqno(3.4)
$$
For example, the left/right invariants and the left/right coinvariants
of the multiple weak Hopf module $_A^AA^A_A\equiv(A,\mu,\mu,\Delta,\Delta)$
of a WBA $A$ are the left/right integrals $I^{L/R}$ and the right/left 
subalgebras $A^{R/L}$, respectively. Dualizing left/right actions or coactions
of a WBA $A$ with the help of dual bases in $A$ and $\hat A$ with respect to 
the canonical pairing, one arrives at right/left coactions or actions of the 
dual WBA $\hat A$, respectively, e.g. 
$$\eqalignno{
  {^{\hat A}M}&:\qquad\hat\delta_L(m)
  :=\beta_i\otimes m\cdot b_i,\quad m\in M_A,&(3.5a) \cr
  {_{\hat A}M}&:\qquad \varphi\cdot m:=m_0\langle\varphi,m_1\rangle, 
  \quad m\in M^A,\varphi\in\hat A,&(3.5b)\cr}
$$
and the invariants (coinvariants) with respect to $A$ become
coinvariants (invariants) with respect to $\hat A$.   

If $A$ is not only a WBA, but also a WHA one can say more about the
invariants and coinvariants of (multiple) WHMs:

\medskip\noindent
{\bf Lemma 3.1} {\sl Let $A$ be a WHA. 
\item{i)} The coinvariants and the invariants 
of a WHM of $A$ can be equivalently characterized as 
$$\eqalign{ 
                        C({M_A^A})&=
   \{m\in M\,\vert\, \delta_R(m)=m\cdot\UN^{(1)}\otimes\UN^{(2)}\},\cr
                        C({_AM^A})&=
   \{m\in M\,\vert\, \delta_R(m)=\UN^{(1)}\cdot m\otimes\UN^{(2)}\},\cr
                        C({_A^AM})&=
   \{m\in M\,\vert\, \delta_L(m)=\UN^{(1)}\otimes\UN^{(2)}\cdot m\},\cr
                        C({^AM_A})&=
   \{m\in M\,\vert\, \delta_L(m)=\UN^{(1)}\otimes m\cdot\UN^{(2)}\},\cr}
  \eqno(3.6a)
$$$$\eqalign{
                        I({M_A^A})&=
   \{m\in M\,\vert\, m_0\cdot a\otimes m_1=m_0\otimes m_1S(a),a\in A\},\cr
                        I({_AM^A})&=
   \{m\in M\,\vert\, m_0\otimes am_1=S(a)\cdot m_0\otimes m_1,a\in A\},\cr
                        I({_A^AM})&=
   \{m\in M\,\vert\, m_{-1}\otimes a\cdot m_0=S(a)m_{-1}\otimes m_1,a\in A\},\cr
                        I({^AM_A})&=
   \{m\in M\,\vert\, m_{-1}a\otimes m_0=m_{-1}\otimes m_0\cdot S(a),a\in A\}.
   \cr}\eqno(3.6b)
$$ 
\item{ii)} The following maps define projections from WHMs onto their 
coinvariants and invariants, respectively
$$\eqalign{
            P^A(m)&:=m_0\cdot S(m_1),\ m\in {M_A^A},\cr 
            ^AP(m)&:=S(m_{-1})\cdot m_0,\ m\in {_A^AM},\cr}\quad\eqalign{ 
            \bar P^A(m)&:=S^{-1}(m_1)\cdot m_0,\ m\in {_AM^A},\cr
            ^A{\bar P}(m)&:=m_0\cdot S^{-1}(m_{-1}),\ m\in {^AM_A},\cr}
\eqno(3.7a)
$$$$\eqalign{
         P_A(m)&:=m_0\cdot R(m_1),\ m\in {M_A^A},\cr 
         _AP(m)&:=L(m_{-1})\cdot m_0,\ m\in {_A^AM},\cr}\quad\eqalign{ 
          _A{\bar P}(m)&:=\bar L(m_1)\cdot m_0,\ m\in {_AM^A},\cr
          \bar P_A(m)&:=m_0\cdot\bar R(m_{-1}),\ m\in {^AM_A},\cr}
\eqno(3.7b)
$$
where $S$ is the antipode and $R, \bar R, L, \bar L$ are the projection 
maps (1.17) to integrals in the WHA $A$.
\item{iii)} In case of the multiple WHMs ${_AM^A_A}$ and ${_A^AM_A}$ the 
coinvariants are left and right $A$-modules with respect to the induced 
left and right adjoint actions, respectively.}

\smallskip\noindent{\it Proof.} i) The characterization (3.6a) of coinvariants 
and the form (3.7a) of the projections onto them have been already proved in 
[2, p.409]. Concerning the invariants of $M^A_A$, first we note that the set 
given in (3.6b) is contained in the set of invariants defined in (3.4) since
$$\eqalign{
  m\cdot a&=(id\otimes\varepsilon)(m_0\cdot a^{(1)}\otimes m_1a^{(2)})
           =(id\otimes\varepsilon)(m_0\otimes m_1S(a^{(1)})a^{(2)})\cr
          &=(id\otimes\varepsilon)(m_0\otimes m_1\Pi^R(a))
           =(id\otimes\varepsilon)(\delta_R(m\cdot\Pi^R(a))
           =m\cdot\Pi^R(a),\cr}\eqno(3.8)
$$
for all $a\in A$. Using the third identity in (1.8) the opposite 
containment is as follows 
$$\eqalign{
  m_0\cdot a\otimes m_1&=m_0\cdot\UN^{(1)}a\otimes m_1\UN^{(2)}
           =m_0\cdot a^{(1)}\otimes m_1\Pi^L(a^{(2)})\cr
          &=m_0\cdot a^{(1)}\otimes m_1a^{(2)}S(a^{(3)})
           =(m\cdot a^{(1)})_0\otimes (m\cdot a^{(1)})_1S(a^{(2)})\cr
          &=(m\cdot \Pi^R(a^{(1)}))_0\otimes 
            (m\cdot \Pi^R(a^{(1)}))_1S(a^{(2)})
           =m_0\otimes m_1\Pi^R(a^{(1)})S(a^{(2)})\cr
          &=m_0\otimes m_1S(a),\qquad a\in A, m\in I({M^A_A}).\cr}\eqno(3.9)
$$
The cases of the other three WHMs can be proved similarly.

ii) The image of the map $P_A$ is in $I({M^A_A})$ due to the 
defining property (1.9) of the right integrals in $A$. Applying $P_A$ to
an invariant $m\in I({M^A_A})$ and using their characterization (3.6b) and 
the non-degeneracy property (3.3),
$$\eqalign{
  m_0\cdot R(m_1)&:=m_0\cdot ((m_1b_i)\Sl\hat S^2(\beta_i))
           =m_0\cdot ((m_1S(b_i))\Sl\hat S(\beta_i))\cr
          &=(m_0\cdot b_i)\cdot (m_1\Sl\hat S(\beta_i))
           =m_0\cdot b_im_1^{(2)}\langle m_1^{(1)},\hat S(\beta_i)\rangle\cr
          &=m_0\cdot S(m_1^{(1)})m_1^{(2)}
           =m_0\cdot \Pi^R(m_1)=m,\quad m\in I({M^A_A})\cr}\eqno(3.10)
$$
follows, that is $P_A$ is a projection onto the invariants of $M^A_A$. 
The cases of projections onto the invariants of the other three WHMs can 
be proved similarly.

iii) We have to show that the maps
$$\eqalign{
  \nu_L(a\otimes m)\equiv a\star m&:=a^{(1)}\cdot m\cdot S(a^{(2)}), \qquad
     a\in A,\, m\in C({_AM^A_A}),\cr
  \nu_R(m\otimes a)\equiv m\star a&:=S(a^{(1)})\cdot m\cdot a^{(2)}, \qquad
     a\in A,\, m\in C({_A^AM_A})\cr}\eqno(3.11)
$$
provide a left and a right $A$-module structure $(C({_AM^A_A}),\nu_L)$ 
and $(C({_A^AM_A}),\nu_R)$, respectively. The image of the map $\nu_L$ is 
in $C({_AM^A_A})$, because for all $a\in A$ and $m\in C({_AM^A_A})$ 
$$\eqalign{
   \delta_R(a^{(1)}\cdot m\cdot S(a^{(2)}))
  &=a^{(11)}\cdot (\UN^{(1)}\cdot m)\cdot S(a^{(2)})^{(1)}\otimes 
    a^{(12)}\UN^{(2)}S(a^{(2)})^{(2)}\cr 
  &=a^{(1)}\cdot m\cdot S(a^{(4)})\otimes a^{(2)}S(a^{(3)})\cr
  &=a^{(1)}\cdot m\cdot S(a^{(3)})\otimes\Pi^L(a^{(2)})\in M\otimes A^L.\cr}
  \eqno(3.12)
$$
The map $\nu_L$ is clearly a left $A$-action, i.e. 
$a\star(b\star m)=ab\star m$, for all $a,b\in A$ and $m\in C({_AM^A_A})$,
moreover, for all $m\in C({_AM^A_A})$
$$\eqalign{
   m&=\UN\cdot m=(\UN\cdot m)_0\cdot\Pi^R((\UN\cdot m)_1)
    =\UN^{(1)}\cdot m_0\cdot\Pi^R(\UN^{(2)}m_1)\cr
   &=\UN^{(1)}\cdot m\cdot\Pi^R(\UN^{(2)})
    =\UN^{(1)}\cdot m\cdot S(\UN^{(2)})=\UN\star m,\cr}
  \eqno(3.13)
$$
where we used the identities (3.3) and (3.2), the property (3.6a) of
the coinvariants and (1.10). 
The proof of the case $(C({_A^AM_A}),\nu_R)$ is similar.\qed

\smallskip
Extending the result of [2, p.410] concerning the structure of a WHM, the 
structure of a multiple WHM is given by the following

\medskip\noindent
{\bf Theorem 3.2} {\sl i) Let ${_AM^A_A}$ be a multiple weak Hopf module of 
the WHA $A$. Then ${_AM^A_A}$ is isomorphic as a multiple WHM to 
$_A(C(M)\times A^A_A)$, which as a left $A$-module is isomorphic to the 
product of the left $A$-modules $(C(M),\star)$ defined in the previous Lemma 
and the left regular module $_AA$, while the right 
$A$-module and right $A$-comodule structures are inherited from the 
WHM $A^A_A\equiv (A,\mu,\Delta)$. 

ii) In particular, $_A{\hat A}^A_A\equiv (\hat A,\mu_L,\mu_R,\delta_R)$
is a multiple WHM with structure maps
$$\eqalignno{
  \mu_L(a\otimes\varphi)\equiv a\cdot\varphi &:=\varphi\Sl S^{-1}(a),&(3.14a)\cr
  \mu_R(\varphi\otimes a)\equiv\varphi\cdot a&:=S(a)\Sr\varphi,
                  \qquad a\in A,\varphi\in\hat A,&(3.14b)\cr
  \delta_R(\varphi)\equiv \varphi_0\otimes\varphi_1
      &:=\beta_i\varphi\otimes b_i,&(3.14c)\cr}
$$
where $\{ b_i\}\subset A$ and $\{\beta_i\}\subset\hat A$ are dual bases 
with respect to the canonical pairing, therefore 
$$
  _A{\hat A}^A_A\simeq {_A(C(\hat A)\times A^A_A)}={_A(\hat I^L\times A^A_A)},
  \eqno(3.15)
$$
where $\hat I^L$ is the space of left integrals in the WHA $\hat A$.}

\noindent{\it Proof.} i) As a $k$-linear space ${_A(C(M)\times A^A_A)}\equiv 
(C(M)\times A, \mu_L,\mu_R,\delta_R)$ is 
(see (2.2a))
$$\eqalign{
  C(M)\times A&:=\UN^{(1)}\star C(M)\otimes\UN^{(2)}A
  =\UN^{(1)}\cdot C(M)\cdot S(\UN^{(2)})\otimes \UN^{(3)}A\cr
 &=C(M)\cdot \UN^{(1)}S(\UN^{(2)})\otimes \UN^{(3)}A
  =C(M)\cdot \Pi^L(\UN^{(1)})\otimes \UN^{(2)}A\cr
 &=C(M)\cdot S(\UN^{(1)})\otimes \UN^{(2)}A\cr}\eqno(3.16)
$$
due to the fact that 
$$
  x^R\cdot m=x^R\cdot (\UN\star m):=
             x^R\UN^{(1)}\cdot m\cdot S(\UN^{(2)})=m\cdot x^R;\quad 
        m\in C(M), x^R\in A^R,\eqno(3.17)
$$
which follows from the identities (3.13) and (1.12). 
One can easily check that the maps 
$$\eqalign{
  a\cdot (\sum_i n_i\otimes b_i)&:=\sum_i a^{(1)}\star n_i\otimes 
                                        a^{(2)}b_i,\cr
  (\sum_i n_i\otimes b_i)\cdot a&:=\sum_i n_i\otimes b_ia,\cr
  \delta_R(\sum_i n_i\otimes b_i)&:=\sum_i n_i\otimes b_i^{(1)}\otimes 
                                         b_i^{(2)},\cr}
  \qquad a\in A,\ \sum_i n_i\otimes b_i\in C(M)\times A.\eqno(3.18)
$$
provide $C(M)\times A$ with a  multiple WHM-structure. The $k$-linear maps 
$U\colon {_AM^A_A}\to C(M)\times A$ and $V\colon C(M)\times A\to {_AM^A_A}$
defined by
$$
  U(m):=m_0\cdot S(m_1)\otimes m_2,\qquad
  V(\sum_i n_i\otimes b_i):=\sum_i n_i\cdot b_i
  \eqno(3.19)
$$
are inverses of each other [2, p.411], i.e. $V\circ U=id_M$ and 
$U\circ V=id_{C(M)\times A}$. In order to prove that ${_AM^A_A}$
and ${_A(C(M)\times A^A_A)}$ are isomorphic as multiple WHMs as well,
we have to show that both $U$ and $V$ are left and right $A$-module and 
right $A$-comodule maps. We can restrict ourselves to the left $A$-module 
properties, because the two other properties were already shown in [2, p.411].
$$\eqalignno{
  U(a\cdot m)&=(a\cdot m)_0\cdot S((a\cdot m)_1)\otimes (a\cdot m)_2\cr
            &=a^{(1)}\cdot m_0\cdot S(m_1)S(a^{(2)})\otimes a^{(3)}m_2\cr
            &=a^{(1)}\star (m_0\cdot S(m_1))\otimes a^{(2)}m_2\cr
            &=a\cdot U(m),\quad a\in A,\, m\in M,&(3.20a)\cr\cr
  V(a\cdot (\sum_i n_i\otimes b_i))
           &=\sum_i V(a^{(1)}\star n_i\otimes a^{(2)}b_i)
            =\sum_i V(a^{(1)}\cdot n_i\cdot S(a^{(2)})\otimes a^{(3)}b_i)\cr
           &=\sum_i a^{(1)}\cdot n_i\cdot S(a^{(2)})a^{(3)}b_i
            =\sum_i a^{(1)}\cdot n_i\cdot \Pi^R(a^{(2)})b_i\cr
           &=\sum_i a^{(1)}\Pi^R(a^{(2)})\cdot n_i\cdot b_i
            =a\cdot (\sum_i n_i\cdot b_i)\cr
           &=a\cdot V(\sum_i n_i\otimes b_i),
   \quad a\in A,\ \sum_i n_i\otimes b_i\in C(M)\times A, &(3.20b)\cr}
$$
where we used (3.17) in the fifth equality of (3.20b). 

ii) The WHM structure $\hat A^A_A\equiv (\hat A,\mu_R,\delta_R)$ given by 
(3.14b and c) of the multiple WHM ${_A\hat A^A_A}$ has been shown in [2, p.409].
The map $\mu_L$ defined in (3.14a) is clearly a left $A$-module map on 
$\hat A$ that commutes with the given right $A$-module map $\mu_R$.
The right comodule map $\delta_R$ is also a left $A$-module map since
$$\eqalign{
  \delta_R(a\cdot\varphi)
 &:=\delta_R(\varphi\Sl S^{-1}(a))
  =\delta_R(\varphi\Sl \bar\Pi^L(a^{(2)})S^{-1}(a^{(1)}))\cr
 &=\delta_R(((\hat\UN\Sl\bar\Pi^L(a^{(2)}))\varphi)\Sl S^{-1}(a^{(1)}))
  =\delta_R(((\hat\UN\Sl a^{(2)})\varphi)\Sl S^{-1}(a^{(1)}))\cr
 &=\delta_R(((\hat\UN\Sl a^{(3)}S^{-1}(a^{(2)}))(\varphi\Sl S^{-1}(a^{(1)}))\cr
 &=\beta_i(\hat\UN\Sl \bar\Pi^R(a^{(2)}))(\varphi\Sl S^{-1}(a^{(1)}))
     \otimes b_i\cr
 &=(\beta_i\Sl \bar\Pi^R(a^{(2)}))(\varphi\Sl S^{-1}(a^{(1)}))
     \otimes b_i
  =((\beta_i\Sl a^{(2)})\varphi)\Sl S^{-1}(a^{(1)})\otimes b_i\cr
 &=(\beta_i\varphi)\Sl S^{-1}(a^{(1)})\otimes a^{(2)}b_i
   =a^{(1)}\cdot \varphi_0\otimes a^{(2)}\varphi_1,
  \qquad a\in A,\varphi\in\hat A,\cr}\eqno(3.21)
$$
where we used the identities (1.6) and (1.10--11).
Hence, the maps (3.14) provides $\hat A$ with a multiple WHM structure, and 
the statement (3.15) follows from the previously proved structure of a general 
multiple weak Hopf module. By dualizing the right $A$-coaction to left
$\hat A$-action as in (3.5b), the right coinvariants $C(\hat A^A)$ become 
the left invariants of the left regular module ${_{\hat A}{\hat A}}$, which 
is the space of left integrals $\hat I^L$ in $\hat A$.\qed

\medskip\noindent
{\bf Corollary 3.3} {\sl The left regular $A$-module ${_AA}$ is injective, 
i.e., $A$ is a quasi-Frobenius algebra.}

\smallskip\noindent
{\it Proof.} The inverse of the antipode provides the isomorphism 
of the right $A$-modules 
$$
 \hat S^{-1}\colon(\hat A_A,\Sl)\to (\hat A_A,\mu_R),\eqno(3.22)
$$
with right action $\mu_R$ given in (3.14b) and the structure theorem of 
multiple weak Hopf modules implies that $(\hat A_A,\mu_R)$ 
is isomorphic to a direct summand of the free right $A$-module
$\hat I^L\otimes A_A$. Therefore, $(\hat A_A,\Sl)$ is a projective right 
$A$-module, which implies the injectivity of its $k$-dual, that is of ${_AA}$. 
Hence, $A$ is a quasi-Frobenius algebra [5, p.414], which has been already 
established in [2, p.413].\qed

\medskip\noindent
{\bf Corollary 3.4} {\sl ${_AI^R}$ and $A^R_A$ are $A$-duals of each other.
${_A{\hat I^L}}$ is the right conjugate module ${_A\Cr{I^R}}$ of ${_AI^R}$ 
and they are the direct sum of simple submodules
$$\eqalignno{
  {_AI^R}&=\bigoplus_p {_AI^R_p},\qquad I^R_p:=z^R_pI^R,&(3.23a)\cr
  {_A{\hat I^L}}&=\bigoplus_p {_A{\hat I^L_p}},\qquad 
   \hat I^L_p:=z^R_p\star\hat I^L,&(3.23b)\cr}
$$
where $\{ z_p^R\}_p$ is the set of primitive orthogonal 
idempotents in $Z^R$.}

\smallskip\noindent 
{\it Proof.} Since the right integrals $I^R$ form a left ideal in $A$ 
and ${_AA}$ is injective by Corollary 3.3, it follows [5, p.392]
that every $\phi\in {\rm Hom}\, ({_AI^R},{_AA})$ can be extended to $\bar\phi
\in{\rm Hom}\, ({_AA},{_AA})$. But every such element $\bar\phi$ is given by 
a right multiplication of an element $a\in A$, hence any $\phi$ is given by the 
restriction of a right multiplication to $I^R$:
$$
  \phi(r)=\bar\phi(r)\equiv ra=r\Pi^R(a),\quad  r\in I^R.
$$
This establish that ${\rm Hom}\, ({_AI^R},{_AA})\simeq A^R$ as a $k$-linear 
space. The isomorphic right $A$-module structure, 
${\rm Hom}\, ({_AI^R},{_AA})_A\simeq A^R_A$, follows from
$$\eqalign{
  (\phi_x\cdot a)(r)&:=\phi_x(r)a =(rx)a=r\Pi^R(xa)=: r(x\cdot a)\cr
 &=\phi_{x\cdot a}(r),\quad x\in A^R, r\in I^R, a\in A.}\eqno(3.24)
$$
The proof of other duality relation is as follows. A map 
$f\in {\rm Hom}\, ({A^R_A},A_A)$ is just a left multiplication 
with the image $f(\UN)\in A$,
$$
  f(x)=f(\Pi^R(\UN x))=:f(\UN\cdot x)=f(\UN)x,\quad x\in A^R,\eqno(3.25)
$$
which should be a right integral, $f(\UN)\in I^R$, because of the module 
homomorphism property of $f$ and (3.25)
$$
  f(\UN)a=f(\UN\cdot a):=f(\Pi^R(\UN a))=f(\UN)\Pi^R(a).
  \eqno(3.26)
$$
The isomorphic left $A$-module strucure of $I^R$ and 
${\rm Hom}\, ({A^R_A},A_A)$ is immediate since it is given by left 
multiplication on the image $f(\UN)\in I^R$ in both cases.

Since in quasi-Frobenius algebras the $A$-duals of simple right $A$-modules 
are simple left $A$-modules [5, p.396], the direct sum decomposition
(3.23a) into simple submodules is induced by the corresponding 
decomposition (2.22) of $A^R_A$.

${_A\Cr{I^R}}={_A{\hat I^L}}$ follows since using (3.14a and b), (1.10) and 
(1.8), the left $A$-module structures of ${_A{\hat I^L}}$ and ${_AI^R}$ are 
related as required by (2.6), 
$$\eqalign{
  \langle a\star\lambda,r\rangle
 &:=\langle a^{(1)}\cdot\lambda\cdot S(a^{(2)}),r\rangle   
  =\langle S^2(a^{(2)})\Sr\lambda\Sl S^{-1}(a^{(1)}),r\rangle\cr   
 &=\langle\lambda, S^{-1}(a^{(1)})rS^2(a^{(2)})\rangle   
  =\langle\lambda, S^{-1}(a^{(1)})r\Pi^R(S^2(a^{(2)}))\rangle\cr   
 &=\langle\lambda, S^{-1}(a^{(1)})rS^2(\Pi^R(a^{(2)}))\rangle   
  =\langle\lambda, S^{-1}(a^{(1)})rS^3(\bar\Pi^L(a^{(2)}))\rangle\cr  
 &=\langle\lambda, S^{-1}(\UN^{(1)})S^{-1}(a)rS^3(\UN^{(2)})\rangle
  =\langle \UN\star\lambda, S^{-1}(a)r\rangle\cr
 &=\langle\lambda, S^{-1}(a)r\rangle,
  \quad a\in A,\lambda\in\hat I^L, r\in I^R,\cr}\eqno(3.27)
$$
and the restriction of the canonical pairing to these integrals is 
non-degenerate. Hence, ${_A\hat I^L}$ is also semisimple and the 
decomposition (3.23b) follows because $z_p^R$ is a central idempotent in $A$
and $I_p^R:=z_p^RI^R=I^Rz_p^R=I^RS^{-1}(z_p^R)=S^{-1}(z_p^R)I^R$.\qed

\medskip\noindent
{\bf Corollary 3.5} {\sl $\hat I^L$ becomes a free rank one left $A^R$- and
$A^L$-module by restricting its left $A$-module structure 
${_A\hat I^L}\equiv (_A\hat I^L,\star)$ to these canonical subalgebras.}  

\smallskip\noindent 
{\it Proof.} Due to Corollary 3.4 ${_A{\hat I^L}}$ is the right conjugate of 
${_AI^R}$, that is ${_AI^R}$ is the left conjugate of ${_A{\hat I^L}}$, 
because $\Cl{\Cr{M}}=M$ for any left $A$-module $M$.
Hence, if we prove the left $A$-module isomorphism
$A^L\simeq\hat I^L\times I^R$ then Lemma 2.8 and Corollary 2.7 lead
to the desired result.  

The restriction of the multiple WHM isomorphism 
${_A{\hat A}^A_A}\simeq {_A(\hat I^L\times A^A_A)}$ in (3.15) to the right 
invariants leads to the isomorphism ${_AI(\hat A_A)}\simeq 
{_AI(\hat I^L\times A_A)}$ of left $A$-modules, where the left $A$-module 
structure of the right invariants is inherited from that of the corresponding 
multiple WHM. 
In our case $I(\hat A_A)\equiv I(\hat A_A,\mu_R)=\hat A^L$ and 
$I(\hat I^L\times A_A)=\hat I^L\times I^R$. The latter equality can be seen 
by using the form (3.7b) of the projection $P_A$ to right invariants of the 
WHM $\hat I^L\times A^A_A$. To prove the former equality we note that the 
invariants of the right $A$-module $(\hat A_A,\mu_R)$ are the coinvariants 
of the dual left $\hat A$-comodule $({^{\hat A}{\hat A}},\hat\delta_L)$ given 
by (3.5a). Since in this case 
$\delta_L(\varphi)=\hat S(\varphi^{(2)})\otimes\varphi^{(1)}$, applying 
$\hat S^{-1}\otimes\varepsilon$ to the defining identity (3.4) of left 
coinvariants and using (1.10) one arrives at 
$\hat A^L=C({^{\hat A}{\hat A}},\hat\delta_L)=I(\hat A_A,\mu_R)$.    
Therefore, $\hat A^L\simeq\hat I^L\times I^R$ as left $A$-modules. 
However, $\hat A^L\simeq A^L$ also holds since the invertible map
$\hat S\circ\kappa_L\colon A^L\to\hat A^L$ with $\kappa_L$ in (1.5)
is an $A$-module map: 
$$\eqalign{
  (\hat S\circ\kappa_L)(a\cdot x^L)
 :&=(\hat S\circ\kappa_L)(\Pi^L(ax^L))
  =\hat S(\Pi^L(ax^L)\Sr\hat\UN)
  =\hat S(ax^L\Sr\hat\UN)\cr
  &=\hat S(a\Sr\kappa_L(x^L))=\hat S(\kappa_L(x^L))\Sl S^{-1}(a)
  =: a\cdot(\hat S\circ\kappa_L)(x^L),\cr}
  \eqno(3.28)
$$ 
where $a\in A$ and $x^L\in A^L$. Thus, $A^L\simeq\hat I^L\times I^R$ as 
left $A$-modules.\qed

\bigskip\noindent
{\bf 4. Existence of non-degenerate left integrals in weak Hopf algebras}
\medskip

Here we prove the generalization of the Larson--Sweedler theorem [10]. 

\smallskip\noindent
{\bf Theorem 4.1} {\sl A finite dimensional weak bialgebra $A$ over a field $k$ 
is a weak Hopf 
algebra iff there exists a non-degenerate left integral in $A$.}\smallskip

\noindent{\it Proof. Sufficiency.} A left integral $l\in A$ obeys the 
defining property $al=\Pi^L(a)l,a\in A$. Non-degeneracy means that the maps
$$
  \eqalign{R_l\colon\hat A&\to A\cr
           \varphi&\mapsto (\varphi\Sr l)\cr}\qquad
  \eqalign{L_l\colon\hat A&\to A\cr
           \varphi&\mapsto (l\Sl\varphi) \cr}
$$
are bijections. This implies that there exist $\lambda,\rho\in\hat A$ such that $l\Sl\rho\equiv L_l(\rho)=\UN=R_l(\lambda)\equiv\lambda\Sr l.$
Let us define the $k$-linear maps $S\colon A\to A$ and 
$\hat S\colon\hat A\to\hat A$ by
$$\eqalign{
  S(a)&:=(R_l\circ\hat L_\lambda)(a)\equiv(\lambda\Sl a)\Sr l
      = l^{(1)}\langle al^{(2)},\lambda\rangle,\cr
  \hat S(\psi)&:=(\hat R_\lambda\circ L_l)(\psi)
               \equiv(l\Sl \psi)\Sr \lambda
               =\lambda^{(1)}\langle \psi\lambda^{(2)}, l\rangle.\cr}\eqno(4.1)
$$
They are transposed to each other with respect to the canonical pairing and
$\hat S(\rho)=\lambda$. Now we prove that $\lambda$ ($\rho$) is a 
non-degenerate left (right) integral in $\hat A$ obeying 
$l\Sr \lambda=\hat\UN=l\Sr \rho$.

Since $R_l$ and $L_l$ are bijections the identities
$$\eqalign{
  R_l(\psi\lambda)&=(\psi\lambda)\Sr l
  =\psi\Sr(\lambda\Sr l)=\psi\Sr\UN
  =\hat\Pi^L(\psi)\Sr\UN
  =\hat\Pi^L(\psi)\Sr(\lambda\Sr l)\cr
  &=R_l(\hat\Pi^L(\psi)\lambda)\cr
  L_l(\rho\psi)&=l\Sl (\rho\psi)
  =(l\Sl\rho)\Sl \psi=\UN\Sl\psi
  =\UN\Sl\hat\Pi^R(\psi)
  =(l\Sl\rho)\Sl\hat\Pi^R(\psi)\cr
  &=L_l(\rho\hat\Pi^R(\psi))\cr}\eqno(4.2)$$
imply that $\lambda$ ($\rho$) is a left (right) integral in $\hat A$. Using the 
properties $l\Sl \rho=\UN=\lambda\Sr l$
$$
  \eqalign{
  \hat\Pi^R(l\Sr\rho)&=\hat\Pi^R(\rho^{(1)})\langle\rho^{(2)},l\rangle
  =\hat\UN^{(1)}\langle\rho\hat\UN^{(2)},l\rangle
  =\hat\UN^{(1)}\langle\hat\UN^{(2)},l\Sl\rho\rangle=\hat\UN,\cr
  \hat{\bar\Pi}^R(l\Sr\lambda)&=
   \hat{\bar\Pi}^R(\lambda^{(1)})\langle\lambda^{(2)},l\rangle
  =\hat\UN^{(1)}\langle\hat\UN^{(2)}\lambda,l\rangle
  =\hat\UN^{(1)}\langle\hat\UN^{(2)},\lambda\Sr l\rangle=\hat\UN,\cr}
  \eqno(4.3)
$$
which imply that $l\Sr\rho=\hat\UN=l\Sr\lambda$ since 
$l\Sr\rho, l\Sr\lambda\in\hat A^L$ and both of the 
$\hat A^L-A^R$ and $\hat A^L-A^L$ pairings are nondegenerate.
The proved properties of $\lambda,\rho\in\hat A$ allow us to construct the
inverse of the map $\hat S$:
$$
  \hat S^{-1}(\psi):=(\hat R_\rho\circ R_l)(\psi)
                   \equiv(\psi\Sr l)\Sr\rho
                   =\rho^{(1)}\langle\rho^{(2)}\psi,l\rangle.\eqno(4.4)
$$
Indeed, for all $\psi\in\hat A$ one obtains
$$\eqalign{
  (\hat S^{-1}\circ\hat S)(\psi)
  &:=\rho^{(1)}\langle\rho^{(2)}\lambda^{(1)},l\rangle
    \langle\psi\lambda^{(2)},l\rangle
   =\rho^{(1)}\langle\rho^{(2)}\hat\Pi^R(\psi^{(1)})\lambda^{(1)},l\rangle
    \langle\psi^{(2)}\lambda^{(2)},l\rangle\cr
  &=\rho^{(1)}\psi^{(1)}\langle\rho^{(2)}\psi^{(2)}\lambda^{(1)},l\rangle
    \langle\psi^{(3)}\lambda^{(2)},l\rangle\cr
  &=\rho^{(1)}\psi^{(1)}\langle\rho^{(2)}\hat\Pi^L(\psi^{(2)})\lambda^{(1)},
     l\rangle\langle\lambda^{(2)},l\rangle\cr
  &=\rho^{(1)}\psi^{(1)}\langle\rho^{(2)}\hat\Pi^L(\psi^{(2)}),l\rangle
   =\rho^{(1)}\psi\langle\rho^{(2)},l\rangle=\psi.\cr}\eqno(4.5)
$$
Therefore, the transposed map $S^{-1}:=(\hat S^{-1})^t\equiv L_l\circ\hat 
L_\rho$ is the inverse of $S$ and invertibility of $S:=R_l\circ \hat L_\lambda$ given in (4.1) implies that $\hat L_\lambda$ and $\hat L_\rho$ are invertible.
Hence, $\lambda$ and $\rho$ are non-degenerate left and right integrals in 
$\hat A$, respectively.

Since $\rho\in\hat A$ is a non-degenerate right integral, there exists
$r\in A$ such that $\rho\Sl r=\hat\UN$. In a similar way as before,
one proves that $r$ is a right integral obeying $r\Sl\rho=\UN$:
$$\eqalign{
  \hat L_\rho(ra)&=\rho\Sl (ra)=\hat\UN\Sl a
  =\hat\UN\Sl\Pi^R(a)=\hat L_\rho(r\Pi^R(a))\cr
  \Pi^R(r\Sl\rho)&=
   \langle\rho,r^{(1)}\rangle S(\bar\Pi^L(r^{(2)}))
  =\langle\rho,r\UN^{(1)}\rangle S(\UN^{(2)})
  =\langle\rho\Sl r,\UN^{(1)}\rangle S(\UN^{(2)})=\UN,\cr}\eqno(4.6)
$$
hence, $r\Sl\rho=\UN$ follows since $\rho$ is a right integral and the
$A^R-\hat A^R$ pairing is non-degenerate. But then $S=L_r\circ\hat R_\rho$ 
also holds (therefore, $r$ is non-degenerate), because 
$$
  \eqalign{(S^{-1}&\circ L_r\circ\hat R_\rho)(a)
  =(L_l\circ\hat L_\rho\circ L_r\circ\hat R_\rho)(a)
   =\langle r^{(1)}a,\rho\rangle\langle r^{(2)}l^{(1)},\rho\rangle l^{(2)}\cr
  &=\langle r^{(1)}a^{(1)},\rho\rangle\langle r^{(2)}\Pi^L(a^{(2)})l^{(1)},
    \rho\rangle l^{(2)}
   =\langle r^{(1)}a^{(1)},\rho\rangle\langle r^{(2)}a^{(2)}l^{(1)},\rho\rangle      a^{(3)}l^{(2)}\cr
  &=\langle r^{(1)},\rho\rangle\langle r^{(2)}\Pi^R(a^{(1)})l^{(1)},\rho\rangle      a^{(2)}l^{(2)}
   =\langle \Pi^R(a^{(1)})l^{(1)},\rho\rangle a^{(2)}l^{(2)}\cr
  &=\langle l^{(1)},\rho\rangle al^{(2)}=a.\cr}\eqno(4.7)
$$

Now, the defining properties (1.2) of the antipode are fulfilled 
for the map $S:=R_l\circ\hat L_\lambda=L_r\circ\hat R_\rho$, because for 
all $a\in A$ one has 
$$\eqalignno{
  a^{(1)}S(a^{(2)})
   &=a^{(1)}l^{(1)}\langle a^{(2)}l^{(2)},\lambda\rangle
   =\Pi^L(a)l^{(1)}\langle l^{(2)},\lambda\rangle
   =\Pi^L(a),&(4.8a)\cr
  S(a^{(1)})a^{(2)}
   &=r^{(2)}a^{(2)}\langle r^{(1)}a^{(1)},\rho\rangle
   =r^{(2)}\Pi^R(a)\langle r^{(1)},\rho\rangle
   =\Pi^R(a),&(4.8b)\cr
  S(a^{(1)})a^{(2)}S(a^{(3)})
   &=\Pi^R(a^{(1)})S(a^{(2)})
   =\Pi^R(a^{(1)})l^{(1)}\langle a^{(2)}l^{(2)},\lambda\rangle\cr
  &=l^{(1)}\langle al^{(2)},\lambda\rangle
   =S(a).&(4.8c)\cr}
$$

\noindent
{\it Necessity.} The statement follows from Lemma 2.6 and results 
(Theorem 3.16) in [2], but for completeness we give a full proof using 
the results of the previous chapter.

Applying the structure theorem of multiple WHMs to ${_A{\hat A}^A_A}$ 
given by (3.14) we get the isomorphism 
${_A{\hat A}^A_A}\simeq {_A(\hat I^L\times A^A_A)}$. Moreover, the restriction
of the left $A$-module structure of ${_A\hat I^L}\equiv (_A\hat I^L,\star)$ 
to the coideal subalgebra $A^R\subset A$ leads to a free $A^R$-module 
${_{A^R}\hat I^L}\equiv({_{A^R}\hat I^L},\star)$ with 
a single generator $\lambda_0\in\hat I^L$ due to Corollary 3.5. Hence, 
using the multiple WHM isomorphism $V\colon\hat I^L\times A\to \hat A$ 
given in (3.19) and the presence of the separating idempotent in 
$\hat I^L\times A:=\UN^{(1)}\star\hat I^L\otimes \UN^{(2)}\cdot A=
\hat I^L\cdot S(\UN^{(1)})\otimes \UN^{(2)}A$, which follows from (3.11),
(1.4) and from the property $\Delta(\UN)\in A^R\otimes A^L$,
one obtains
$$\eqalign{
  \hat A&=V(\hat I^L\times A)=V((A^R\star\lambda_0)\times A)
  =V((\lambda_0\cdot S(A^R))\times A)=V(\lambda_0\times S(A^R)A)\cr
  &=V(\lambda_0\times A):=\lambda_0\cdot A:=S(A)\Sr\lambda_0
   =A\Sr\lambda_0,\cr}\eqno(4.9)
$$
which implies the non-degeneracy of the left integral $\lambda_0$
in the weak Hopf algebra $\hat A$.\qed 

\medskip
Since a non-degenerate left integral in a WHA  
provides a non-degenerate associative bilinear form on the dual WHA:

\smallskip\noindent
{\bf Corollary 4.2} {\sl A finite dimensional weak Hopf algebra is a 
Frobenius algebra.} 

\bigskip\noindent
{\bf 5. Grouplike elements and invertible modules}\medskip

In this chapter first we define (left/right) grouplike elements in a WHA $A$. 
Then we give two equivalent descriptions of 
invertible $A$-modules in terms of the canonical coideal subalgebras in $A$ 
and in terms of left (right) grouplike elements in the dual WHA $\hat A$. 

The set of {\it grouplike elements} $G(H)$ in a Hopf algebra $H$ can be 
defined to be [17] 
$G(H):=\{ g\in H\vert \Delta(g)=g\otimes g, \varepsilon(g)\not= 0 \}$.
The grouplike elements are linearly independent, they obey the property
$S(g)g=\UN$ and they form a group. The generalization of this notion to 
a weak Hopf algebra $A$ 
$$
  G(A):=\{ g\in A\vert\Delta(\UN)(g\otimes g)=\Delta(g)
       =(g\otimes g)\Delta(\UN), gS(g)=\UN \}
$$ 
given in [2, p.433] seems to be too restrictive, hence we introduce 
slightly softened generalizations as well:

\smallskip\noindent{\bf Definition 5.1} {\sl The set of right/left grouplike
elements $G_{R/L}(A)$ in a weak Hopf algebra $A$ is defined to be
$$\eqalignno{
  G_R(A):=&\{ g\in A\vert (g\otimes g)\Delta(\UN)=\Delta(g)
         =\Delta(\UN)(g\otimes\Pi^L(g)^{-1}g); \Pi^{R/L}(g)\in A^{R/L}_*\} 
         \qquad &(5.1a)\cr
  G_L(A):=&\{ g\in A\vert (g\Pi^R(g)^{-1}\otimes g)\Delta(\UN)=\Delta(g)
         =\Delta(\UN)(g\otimes g); \Pi^{R/L}(g)\in A^{R/L}_*\} 
          \qquad &(5.1b)\cr}
$$
where $A^{R/L}_*$ denote the set of invertible elements in $A^{R/L}$.
The set of grouplike elements in $A$ is defined to be the intersection
$G(A):=G_R(A)\cap  G_L(A)$.}

\smallskip
Using the form (1.11) of the maps $\Pi^{R/L}$, the  
defining properties (5.1) lead to the relations
$$\eqalignno{
g\in G_R(A)&:\ \Pi^R(g)=S(g)\Pi^L(g)^{-1}g,\ \Pi^L(g)=gS(g)\Rightarrow 
               g\in A_*, \Pi^R(g)=\UN,&(5.2a)\cr 
g\in G_L(A)&:\ \Pi^L(g)=g\Pi^R(g)^{-1}S(g),\ \Pi^R(g)=S(g)g\Rightarrow 
               g\in A_*,\Pi^L(g)=\UN,&(5.2b)\cr}
$$
that is elements of $G_{R/L}(A)$ are themselves invertible.
Using (5.1--2) it is easy to show that $G_R(A)$ and $G_L(A)$, hence $G(A)$, 
too, are groups, $G_L(A)=S(G_R(A))$, and the definition of grouplike
elements $G(A)$ above is equivalent to that of given in [2]. For example, 
$gh\in G_R(A)$ if $g,h\in G_R(A)$ since
$$\eqalign{
  \Pi^R(gh)&=\Pi^R(\Pi^R(g)h)=\Pi^R(\UN h)=\UN\in A^R_*,\cr
  \Pi^L(gh)&=\Pi^L(g\Pi^L(h))=g^{(1)}\Pi^L(h)S(g^{(2)})
            =ghS(h)S(g)\in A^L_*,\cr}
$$$$\eqalign{
  \Delta(gh)&=(g\otimes g)\Delta(\UN)\Delta(h)=(g\otimes g)\Delta(h)
             =(gh\otimes gh)\Delta(\UN),\cr
  \Delta(gh)&=\Delta(g)\Delta(\UN)(h\otimes\Pi^L(h)^{-1}h)
             =\Delta(g)(h\otimes S(h^{-1}))\cr
            &=\Delta(\UN)(g\otimes S(g^{-1}))(h\otimes S(h^{-1}))
             =\Delta(\UN)(gh\otimes\Pi^L(gh)^{-1}gh).\cr}
$$
\noindent{\bf Corollary 5.2} {\sl The element $g\in A$ is right (left) 
grouplike iff $g$ is invertible and obeys the property
$\Delta(g)=(g\otimes g)\Delta(\UN)$ ($\Delta(g)=\Delta(\UN)(g\otimes g)$).}

\smallskip\noindent{\it Proof.} If $g\in A$ is right (left) grouplike
it is invertible due to the discussion above, while the required coproduct 
property follows by definition. Conversely, the relations
$$\eqalign{ 
  \UN&=g^{-1}g=g^{-1}({id}\otimes\varepsilon)(\Delta(g))
      =\UN^{(1)}\varepsilon(g\UN^{(2)})=:\Pi^R(g),\cr 
  \Pi^L(g)&=g^{(1)}S(g^{(2)})=g\UN^{(1)}S(\UN^{(2)})S(g)=gS(g),\cr}
$$
imply that $\Pi^{R/L}(g)\in A^{R/L}_*$. In conclusion, using (1.8) one derives 
$$
  \UN^{(1)}g\otimes\UN^{(2)}=g^{(1)}\otimes\Pi^L(g^{(2)})
  =g\UN^{(1)}\otimes\Pi^L(g\UN^{(2)})
  =g\UN^{(1)}\otimes g\UN^{(2)}S(g).
$$
Multiplying this identity by $\UN\otimes S(g^{-1})$ from the right and using
the form of $\Pi^L(g)$, one arrives the other coproduct property of 
a right grouplike element in (5.1a). The proof for left grouplike elements is 
similar.\qed 

We note that the set $G(A)$ in $G_R(A)$ can also be given by the subset 
of elements satisfying $\Pi^L(g)=\UN$ or by the subset of pointwise invariant 
elements with respect to $S^2$. For verification of the latter claim, we note 
that if $g=S^2(g)$ holds for $g\in G_R(A)$ then $\Pi^L(g)=
gS(g)=S^2(g)S(g)=S(\Pi^L(g))$, that is $\Pi^L(g)$, hence $\Pi^L(g)^{-1}$, 
too, are in $A^L\cap A^R\subset {\rm Center}\, A^L$. Using (5.2a), (5.1a) 
and these consequences one obtains
$$
  \UN=\Pi^R(g)=S(\UN^{(1)})S(g)g\UN^{(2)}
  \Rightarrow \Pi^L(g^{-1})^{-1}=\UN
  \Rightarrow g^{-1},g\in G(A).
$$ 

Now we turn to characterization of invertible modules of WHAs.

\smallskip\noindent{\bf Definition 5.3} {\sl An object $M$ of a  
monoidal category $({\cal L};\times,E)$ is invertible if there exists 
an inverse object $\bar M\in{\rm Obj}\,{\cal L}$ obeying 
$M\times\bar M\simeq E\simeq\bar M\times M$, where 
$E\in{\rm Obj}\, {\cal L}$ is the monoidal unit of the category and $\simeq$ 
denotes equivalence of objects in ${\cal L}$.\footnote{${}^1$}{\rm In case of 
symmetric or braided monoidal categories invertibility is defined by the 
condition $E\simeq M\times\bar M$ [7].}} 

\smallskip\noindent
{\bf Proposition 5.4} {\sl i) Let ${\cal L}$ be the autonomous 
monoidal category of finite dimensional left $A$-modules of a WHA 
$A$ given in Prop. 2.2. The module ${_AM}\in{\rm Obj}\,{\cal L}$ 
is invertible iff it becomes a free rank one left $A^L$- and $A^R$-module 
by restricting $A$ to the subalgebras $A^L$ and $A^R$, respectively.

ii) An invertible module ${_AM}\in{\rm Obj}\,{\cal L}$ is semisimple. Namely, 
it is the direct sum of simple submodules:
$$
  {_AM}=\oplus_p M_{(p,\tau_M(p))},\qquad M_{(p,\tau_M(p))}:=z_p^L\cdot M=
         z_{\tau_M(p)}^R\cdot M,
$$
where $\{ z_p^L\}_p\subset Z^L$ and $\{ z_p^R:=S(z_p^L)\}_p\subset Z^R$ are 
the sets of primitive orthogonal idempotents and $\tau_M$ is a permutation
on them.}

\smallskip\noindent{\it Proof.} i) First we show that ${_AM}$ is invertible iff
$$
  M\times\Cl{M}\simeq A^L\simeq\Cr{M}\times M\eqno(5.3)
$$
as left $A$-modules, where $A^L$ is the unit left $A$-module given in (2.1).

If (5.3) holds then, using the natural equivalences $X^L$ and $X^R$
given in (2.3), $\Cl{M}\simeq A^L\times\Cl{M}
\simeq\Cr{M}\times M\times\Cl{M}\simeq\Cr{M}\times A^L\simeq\Cr{M}$ follows, 
hence, $M$ is invertible. Conversely, let $\bar M$ be the inverse of $M$ 
and let $\sigma\colon\bar M\times M\to A^L$ and
$\tau\colon M\times\bar M\to A^L$ be the corresponding invertible arrows.
We will show that $\Cl{M}\simeq\bar M\simeq\Cr{M}$, which imply (5.3).
The arrow
$$
  \omega:=(X_{A^L}^L)^{-1}(\tau\times\tau)
  (1_M\times\sigma^{-1}\times 1_{\bar M})(X^R_M\times 1_{\bar M})\tau^{-1}
  \in{\rm End}\, {_AA^L}\eqno(5.4)
$$
is invertible; therefore, it is given by the action of an invertible element 
$z^L\in Z^L:=A^L\cap{\rm Center}\, A$ due to ${\rm End}\, {_AA^L}=Z^L\cdot$ 
[2, p.402]. Hence, if $z^L_N\colon N\to N$ denotes the arrow
given by the action $z^L\in Z^L$ for $N\in{\rm Obj}\,{\cal L}$ then
$\{ z^L_N\}_N$ is a natural automorphism of the identity functor on 
${\cal L}$ and $\omega=z^L_{A^L}$. Defining $\tilde\tau:=(z^L_{A^L})^{-1}\tau
\colon M\times\bar M\to A^L$, we have $\tilde\tau=\tau(z^L_{M\times\bar M})^{-1}
=\tau((z^L_M)^{-1}\times 1_{\bar M})$ due to naturality and (1.4). Therefore,
$(z^L_M)^{-1}\times 1_{\bar M}=\tau^{-1}\tilde\tau=\tau^{-1}\omega^{-1}\tau$,
which leads to $z^L_M=(X^L_M)^{-1}(\tau\times 1_M)(1_M\times\sigma^{-1})X^R_M$
due to the form (5.4) of $\omega$ and faithfulness of $-\times 1_{\bar M}$. 
Hence, using naturality and (1.4)
$$\eqalign{
  1_M&=(z^L_M)^{-1}z^L_M
      =(z^L_M)^{-1}(X^L_M)^{-1}(\tau\times 1_M)(1_M\times\sigma^{-1})X^R_M\cr
     &=(X^L_M)^{-1}(\tilde\tau\times 1_M)(1_M\times\sigma^{-1})X^R_M.\cr}
  \eqno(5.5a)
$$  
Then
$$
  (X^R_{\bar M})^{-1}(1_{\bar M}\times\tilde\tau)(\sigma^{-1}\times 1_{\bar M})
  X^L_{\bar M}=1_{\bar M}\eqno(5.5b)
$$  
also holds because of faithfulness of $1_M\times -$ and because of the identity
$$\eqalign{
  1_{A^L}&=\tilde\tau\tilde\tau^{-1}
  =\tilde\tau((X^L_M)^{-1}\times 1_{\bar M})
   (\tilde\tau\times 1_M\times 1_{\bar M})
   (1_M\times\sigma^{-1}\times 1_{\bar M})(X^R_M\times 1_{\bar M})
   \tilde\tau^{-1}\cr
 &=(X^L_{A^L})^{-1}(\tilde\tau\times\tilde\tau)
   (1_M\times\sigma^{-1}\times 1_{\bar M})(X^R_M\times 1_{\bar M})
   \tilde\tau^{-1}\cr
 &=\tilde\tau [1_M\times (X^R_{\bar M})^{-1}(1_{\bar M}\times\tilde\tau)
   (\sigma^{-1}\times 1_{\bar M})X^L_{\bar M}]
   \tilde\tau^{-1}.\cr}\eqno(5.6)
$$
Thus, using the right and left evaluation maps defined in (2.8) and 
(2.14), respectively,
$$\eqalignno{
  \Cl{\mu}:=
  (X^L_{\bar M})^{-1}(E^l_M\times 1_{\bar M})(1_{\CL{M}}\times\tilde\tau^{-1})
  X^R_{\CL{M}}&\colon\Cl{M}\to\bar M,&(5.7a)\cr
 \Cr{\mu}:=  
 (X^R_{\bar M})^{-1}(1_{\bar M}\times E^r_M)(\sigma^{-1}\times 1_{\CR{M}})
  X^L_{\CR{M}}&\colon\Cr{M}\to\bar M&(5.7b)\cr}
$$  
provide the equivalences $\Cl{M}\simeq\bar M\simeq\Cr{M}$ with the
inverse arrows
$$\eqalignno{
  {\Cl{\mu}}^{-1}=
  (X^L_{\CL{M}})^{-1}(\sigma\times 1_{\CL{M}})(1_{\bar M}\times C^l_M)
  X^R_{\bar M}&\colon\bar M\to\Cl{M},&(5.8a)\cr
  {\Cr{\mu}}^{-1}=
  (X^R_{\CR{M}})^{-1}(1_{\CR{M}}\times\tilde\tau)(C^r_M\times 1_{\bar M})
  X^L_{\bar M}&\colon\bar M\to\Cr{M}&(5.8b)\cr}
$$ 
due to the rigidity identities (2.10) and (2.15), respectively,
and due to (5.5a and b). For example, 
$$\eqalign{
  {\Cl{\mu}}^{-1}\Cl{\mu}
 &:=[(X^L_{\CL{M}})^{-1}(\sigma\times 1_{\CL{M}})(1_{\bar M}\times C^l_M)
     X^R_{\bar M}]
    [(X^L_{\bar M})^{-1}(E^l_M\times 1_{\bar M})(1_{\CL{M}}\times
   \tilde\tau^{-1})X^R_{\CL{M}}]\cr
 &=(X^L_{\CL{M}})^{-1}(E^l_M\times 1_{\Cl{M}})
   (1_{\Cl{M}}\times [(X^R_M)^{-1}(1_M\times\sigma)(\tilde\tau^{-1}\times 1_M)
    X^L_M]\times 1_{\Cl{M}})\cr
 &\phantom{=}
   (1_{\Cl{M}}\times C^l_M)X^R_{\Cl{M}}
  =(X^L_{\CL{M}})^{-1}(E^l_M\times 1_{\Cl{M}})
   (1_{\Cl{M}}\times C^l_M)X^R_{\Cl{M}}=1_{\Cl{M}},\cr}
$$
where we used the inverse of (5.5a) in the third equality and (2.10b) in the 
fourth one. 

Now we prove that (5.3) is fulfilled iff $M$ becomes a free rank one $A^L$- and 
$A^R$-module by restricting the left $A$-action 
to these subalgebras. If (5.3) holds then the statement follows from 
Lemma 2.8 and Corollary 2.7. Conversely, suppose that ${_AM}$ becomes a free 
$A^L$- and $A^R$-module with a single generator $m\in M$ by resricting the 
$A$-action to these subalgebras. The elements $\hat m_l$ and $\hat m_r$ of 
the $k$-dual $\hat M$ of $M$ defined by
$$
  \langle\hat m_l,x^R\cdot m\rangle_M:=\varepsilon(x^R),\qquad
  \langle\hat m_r,x^L\cdot m\rangle_M:=\varepsilon(x^L),\qquad
  x^{R/L}\in A^{R/L}\eqno(5.9)
$$
are $A^L$- and $A^R$-generators of $\Cl{M}$ and $\Cr{M}$, respectively, 
because the counit $\varepsilon$ is a nondegenerate functional on $A^R$ 
and on $A^L$. Moreover, 
choosing dual bases $\{ e_i\}_i,\{ f_i\}_i$ in $A^L$ with respect to the 
counit $\varepsilon$, the bases $\{ e_i\cdot \hat m_l\}_i\subset\Cl{M}, 
\{ S^{-1}(f_i)\cdot m\}_i\subset M$ and $\{ S^{-1}(f_i)\cdot\hat m_r\}_i
\subset\Cr{M}, \{ e_i\cdot m\}_i\subset M$ become dual to each other. Indeed,
for the dual $A^L$-bases we have 
$$
  \delta_{ij}=\varepsilon(e_if_j)=\varepsilon(f_jS^2(e_i))
  =\varepsilon(S(e_i)S^{-1}(f_j)).\eqno(5.10) 
$$
The third equality follows from the invariance of the counit with respect to
the antipode: $\varepsilon=\varepsilon\circ S$. The second is the 
consequence of the identities (1.14--15) claiming that $S^2$ is the Nakayama
automorphism $\theta_L\colon A^L\to A^L$ corresponding to the counit as a 
non-degenerate functional on $A^L$. Therefore
$$\eqalign{
  \delta_{ij}&=\varepsilon(S(e_i)S^{-1}(f_j))
  :=\langle\hat m_l,S(e_i)S^{-1}(f_j)\cdot m\rangle_M
 =\langle e_i\cdot\hat m_l,S^{-1}(f_j)\cdot m\rangle_M\, ,\cr
  \delta_{ij}&=\varepsilon(S^{-2}(f_j)e_i)
  :=\langle\hat m_r,S^{-2}(f_j)e_i\cdot m\rangle_M
 =\langle S^{-1}(f_j)\cdot\hat m_r,e_i\cdot m\rangle_M\, 
  .\cr}\eqno(5.11)
$$
Thus, we can prove that the left and right coevaluation 
maps $C^l_M\colon A^L\to M\times\Cl{M}$ and $C^r_M\colon A^L\to\Cr{M}\times M$ 
defined in (2.8) and (2.14) are invertible, that is (5.3) holds: using that
$\UN^{(1)}\otimes\UN^{(2)}=S^{-1}(f_i)\otimes e_i$ (summation suppressed),
rank one $A^L$- and $A^R$-freeness of $M$ in the fourth equalities, 
respectively, and (2.13a) in the sixth equality of (5.12b), one obtains 
$$\eqalignno{
  M\times\Cl{M}&:=\UN^{(1)}\cdot M\otimes \UN^{(2)}\cdot\Cl{M}
  =\UN^{(1)}\cdot M\otimes \UN^{(2)}A^L\cdot\hat m_l\cr
 &=\UN^{(1)}S^{-1}(A^L)\cdot M\otimes \UN^{(2)}\cdot\hat m_l
  =\UN^{(1)}A^L\cdot m\otimes \UN^{(2)}\cdot\hat m_l&(5.12a)\cr
 &=A^LS^{-1}(f_i)\cdot m\otimes e_i\cdot\hat m_l
  =C^l_M(A^L),\cr
  \Cr{M}\times M&:=\UN^{(1)}\cdot \Cr{M}\otimes\UN^{(2)}\cdot M
  =\UN^{(1)}A^R\cdot\hat m_r\otimes \UN^{(2)}\cdot M\cr
 &=\UN^{(1)}\cdot\hat m_r\otimes \UN^{(2)}S(A^R)\cdot M
  =\UN^{(1)}\cdot\hat m_r\otimes \UN^{(2)}A^R\cdot m&(5.12b)\cr
 &=S^{-1}(f_i)\cdot\hat m_r\otimes A^Re_i\cdot m
  =A^LS^{-1}(f_i)\cdot\hat m_r\otimes e_i\cdot m
  =C^r_M(A^L),\cr}
$$
i.e. $C^l_M$ and $C^r_M$ are surjective. Injectivity of $C^l_M$ and $C^r_M$
follow from the faithfulness of $M$ as a left $A^L$- and $A^R$-module,
respectively.  

ii) From (5.3) and Lemma 2.8 we can deduce that
$$
  {\rm End}\,{_AM}\subset {\rm End}\,{_{A^R}M}\cap
  {\rm End}\,{_{A^L}M}
  =({\rm Center}\, A^L)\cdot 
  =({\rm Center}\, A^R)\cdot\ .\eqno(5.13)
$$
Let $m\in M$ be a free $A^L$-generator. The action by an element $x^L\in 
{\rm Center}\, A^L$ on $M$ is an element of ${\rm End}\,{_AM}$ only if
$$
  \Pi^L(a)x^L\cdot m=a^{(1)}S(a^{(2)})x^L\cdot m
  =a^{(1)}x^LS(a^{(2)})\cdot m=\Pi^L(ax^L)\cdot m,\qquad a\in A,\eqno(5.14)
$$
i.e. only if $\Pi^L(a)x^L=\Pi^L(ax^L)$ for all $a\in A$. However, this 
relation implies that $x^L\in{\rm Center}\, 
A$:
$$\eqalign{
  S(a)x^L&=S(a^{(1)})\Pi^L(a^{(2)})x^L=S(a^{(1)})\Pi^L(a^{(2)}x^L)
  =S(a^{(1)})a^{(2)}x^LS(a^{(3)})\cr
 &=\Pi^R(a^{(1)})x^LS(a^{(2)})
  =x^LS(a),\quad a\in A.\cr}\eqno(5.15)
$$
Therefore, $x^L\in A^L\cap {\rm Center}\, A =:Z^L$, that is 
${\rm End}\,{_AM}\subset Z^L\cdot $. The opposite containment is trivial.
The proof of the relation ${\rm End}\,{_AM}= Z^R\cdot $ is similar. 
Hence, the direct summands of ${_AM}$ in the statement ii) are indecomposable
submodules. Since ${_AM}$ is a free rank one $A^L$- and $A^R$-module 
due to i), $\tau_M$ is a permutation and the $k$-dimensions 
of the indecomposable submodules $M_{(p,\tau_M(p))}$ saturate the lower bound 
(2.25) given in Lemma 2.6. Therefore, $M_{(p,\tau_M(p))}$ is simple since it
cannot contain a non-trivial submodule.\qed  

\smallskip Now we turn to the characterization of invertible $\hat A$-modules
in terms of right (left) grouplike elements in the WHA $A$. First, we give 
the connection between (right/left) grouplike elements in $A$ and invertible 
submodules of $({_{\hat A}A},\Sr)$:

\medskip\noindent{\bf Lemma 5.5} {\sl Let $A$ be a WHA and let $F_a:=
(\hat A\Sr a,\Sr)$ denote the cyclic left $\hat A$-submodule of
${_{\hat A}A}:=({_{\hat A}A},\Sr)$ generated by $a\in A$.
 
\item{i)} $g\in A$ is (right/left) grouplike iff $g$ is an 
element of an invertible submodule ${_{\hat A}F}$ of ${_{\hat A}A}$ and 
$g$ obeys the normalization conditions ($\Pi^{R/L}(g)=\UN$) 
$\Pi^R(g)=\UN=\Pi^L(g)$.

\item{ii)} The cyclic submodules $F_g,F_h\subset{_{\hat A}A}$ generated by 
(right/left) grouplike elements are in the same module isomorphism class 
iff $gh^{-1}\in A^T$. 

\item{iii)} In any module isomorphism class of invertible 
submodules of ${_{\hat A}A}$, there is a submodule which contains a right 
(left) grouplike element.}

\smallskip\noindent{\it Proof.} i) Let $g\in G_{R/L}(A)$ or $g\in G(A)$. 
Clearly, $F_g$ is a submodule of ${_{\hat A}A}$ that contains $g$
satisfying the required normalization conditions. According to Prop. 5.4 i) 
invertibility of $F_g$ follows if $F_g$ becomes a free $\hat A^L$- and 
$\hat A^R$-module with the single generator $g$ by restricting the 
$\hat A$-action to these subalgebras. If $g\in G_R(A)$ then 
the identities (1.6--7) and (5.1--2a) lead to the relations
$$\eqalignno{
  \varphi\Sr g&=\UN^{(1)}g\langle\varphi,\UN^{(2)}S(g)^{-1}\rangle
               =\UN^{(1)}g\langle S(g)^{-1}\Sr\varphi,\Pi^L(\UN^{(2)})\rangle\cr
              &=(\hat\Pi^L(S(g)^{-1}\Sr\varphi)\Sr\UN)g
               =\hat\Pi^L(S(g)^{-1}\Sr\varphi)\Sr g,
         \qquad \varphi\in\hat A,&(5.16a)\cr
  \varphi\Sr g&=g\UN^{(1)}\langle\varphi,g\UN^{(2)}\rangle
               =g\UN^{(1)}\langle\varphi\Sl g,\bar\Pi^L(\UN^{(2)})\rangle
               =g\UN^{(1)}\langle\hat{\bar\Pi}^R(\varphi\Sl g),
                 \UN^{(2)})\rangle\cr
              &=g(\hat{\bar\Pi}^R(\varphi\Sl g)\Sr\UN)
               =\hat{\bar\Pi}^R(\varphi\Sl g)\Sr g,
         \qquad \varphi\in\hat A.&(5.16b)\cr}
$$
They imply that $F_g\subset(\hat A^R\Sr g)\cap(\hat A^L\Sr g)=gA^R\cap A^Rg$.
Moreover, if $0=\varphi^L\Sr g=(\varphi^L\Sr\UN)g$ or 
$0=\varphi^R\Sr g=g(\varphi^R\Sr\UN)=g(S(\varphi^R)\Sr\UN)$ for certain 
$\varphi^{L/R}\in\hat A^{L/R}$ then $\varphi^{L/R}=0$, because $g$ is 
invertible and the maps $\hat\kappa_L$ in (1.5) and
the antipode $S$ are bijections. Therefore, $F_g$ is a 
free rank one $\hat A^R$- and $\hat A^L$-module for any $g\in G_R(A)$, hence
for any $g\in G(A)\subset G_R(A)$, too. 
The case of $g\in G_L(A)$ can be proved similarly. 

Conversely, let ${_{\hat A}F}$ be an invertible submodule of ${_{\hat A}A}$.
Then $F$ is a right coideal in $A$ and a free left $\hat A^{L}$- and 
$\hat A^R$-module with a single generator $f\in F$. Thus, one can define 
two projections $\Phi^L_f\colon\hat A\to\hat A^L$ and 
$\bar\Phi^R_f\colon\hat A\to\hat A^R$ by requiring
$$
  \Phi^L_f(\varphi)\Sr f:=\varphi\Sr f,\qquad
  \bar\Phi^R_f(\varphi)\Sr f:=\varphi\Sr f,\eqno(5.17)
$$
for $\varphi\in\hat A$. They are left $\hat A^L$- and $\hat A^R$-module maps,
respectively. Since $F$ is a right coideal in $A$, defining $\hat f_l$ and 
$\hat f_r$ in the $k$-dual $\hat F$ of $F$ like in (5.9) by
$$\eqalign{ 
  \langle\hat f_r,\varphi^L\Sr f\rangle_F
  &=\langle\hat f_r,f^{(1)}\rangle_F\langle f^{(2)},\varphi^L\rangle
   \equiv\langle f\Sl\hat f_r,\varphi^L\rangle
 :=\hat\varepsilon(\varphi^L),\qquad \varphi^L\in\hat A^L,\cr
  \langle\hat f_l,\varphi^R\Sr f\rangle_F
  &=\langle\hat f_l,f^{(1)}\rangle_F\langle f^{(2)},\varphi^R\rangle
   \equiv\langle f\Sl\hat f_l,\varphi^R\rangle
 :=\hat\varepsilon(\varphi^R),\qquad\varphi^R\in\hat A^R,\cr}
  \eqno(5.18)
$$
we have $\Pi^L(f\Sl\hat f_r)=\UN=\Pi^R(f\Sl\hat f_l)$ and
$$\eqalign{
  \Phi^L_f(\varphi)
  &=\hat S(\hat\UN^{(1)})\langle\hat f_r,\hat\UN^{(2)}\Phi^L_f(\varphi)\Sr 
     f\rangle_F
   =\hat S(\hat\UN^{(1)})\langle\hat f_r,\hat\UN^{(2)}\varphi\Sr f\rangle_F\cr
  &=\hat S(\hat\UN^{(1)})\langle f\Sl\hat f_r,\hat\UN^{(2)}\varphi\rangle
   =\hat\Pi^L(\varphi^{(1)})\langle f\Sl\hat f_r,\varphi^{(2)}\rangle,\cr
  \bar\Phi^R_f(\varphi)
  &=\langle\hat f_l,\hat\UN^{(1)} \bar\Phi^R_f(\varphi)\Sr 
    f\rangle_F\hat S^{-1}(\hat\UN^{(2)})
   =\langle\hat f_l,\hat\UN^{(1)}\varphi\Sr 
    f\rangle_F\hat S^{-1}(\hat\UN^{(2)})\cr
  &=\langle f\Sl\hat f_l,\hat\UN^{(1)}\varphi\rangle\hat S^{-1}(\hat\UN^{(2)})
   =\langle f\Sl\hat f_l,\varphi^{(1)}\rangle\hat{\bar\Pi}^R(\varphi^{(2)})
  \cr}\eqno(5.19)
$$
using (5.17--18), (1.11) and upper right and lower left eqs. in (1.8).
Thus, using (1.7--8)
$$\eqalignno{
  \varphi\Sr f&=:\Phi^L_f(\varphi)\Sr f
  =f^{(1)}\langle \hat\Pi^L(\varphi^{(1)}),f^{(2)}\rangle
   \langle f\Sl\hat f_r,\varphi^{(2)}\rangle\cr
 &=f^{(1)}\langle\varphi,\Pi^L(f^{(2)})(f\Sl\hat f_r)\rangle
  =\UN^{(1)}f\langle\varphi,\UN^{(2)}(f\Sl\hat f_r)\rangle,&(5.20a)\cr
  \varphi\Sr f&=:\bar\Phi^R_f(\varphi)\Sr f
  =f^{(1)}\langle\varphi^{(1)},f\Sl\hat f_l\rangle
   \langle\hat{\bar\Pi}^R(\varphi^{(2)}),f^{(2)}\rangle\cr
 &=f^{(1)}\langle\varphi,(f\Sl\hat f_l)\bar\Pi^L(f^{(2)})\rangle
  =f\UN^{(1)}\langle\varphi,(f\Sl\hat f_l)\UN^{(2)}\rangle,&(5.20b)\cr}
$$
for all $\varphi\in\hat A$, which imply
$$
  \UN^{(1)}f\otimes\UN^{(2)}(f\Sl\hat f_r)
  =f^{(1)}\otimes f^{(2)}
  =f\UN^{(1)}\otimes (f\Sl\hat f_l)\UN^{(2)}.\eqno(5.21)
$$
Applying the counit $\varepsilon$ to the first tensor factor we obtain
$$
  \Pi^L(f)(f\Sl\hat f_r)=f=(f\Sl\hat f_l)\bar\Pi^L(f).\eqno(5.22)
$$
Let $g\in F$ such that $\Pi^R(g)=\UN$. 
Then $g=\varphi^R\Sr f=f(\varphi^R\Sr\UN)=:fx^R$ with $x^R\in A^R$ for some 
$\varphi^R\in\hat A^R$ due to $\hat A^R$-freeness of $F$ and 
(1.6). Thus, $\UN=\Pi^R(g)=\Pi^R(f)x^R$, that is $x^R$ is invertible. This 
implies that $g$ is also an $\hat A^{L/R}$-generator of $F$, hence,
(5.21--22) hold for $f=g\in F$, too. 
Since $\UN=\Pi^R(g)=S(g)(g\Sl\hat g_r)$ by assumption and due to the first 
equality of (5.21), $S(g)$, hence $g$, too, is invertible. Since 
$\bar\Pi^L(g)=S^{-1}(\Pi^R(g))=\UN$ due to (1.11), the second equality of 
(5.22) implies that $g\Sl\hat g_l=g$. Hence, the second equality of (5.21) 
together with invertibility of $g$ implies that $g\in G_R(A)$ due to 
Corollary 5.2. The cases $g\in G_L(A), G(A)$ can be proved similarly.

ii) First we note that for $g,h\in G_R(A)$ ($G_L(A),G(A)$) the invertible 
left $\hat A$-modules $F_{gh}$ and $F_g\times F_h$ are isomorphic, 
because the maps
$$\eqalign{
  U\colon F_g\times F_h&\to F_{gh}\cr
           m\otimes n&\mapsto mn\cr}\qquad
  \eqalign{
  V\colon F_{gh}&\to F_g\times F_h\cr
           m&\mapsto \hat\UN^{(1)}\Sr mh^{-1}\otimes\hat\UN^{(2)}\Sr h\cr}
  \eqno(5.23)
$$
are left $\hat A$-module maps, which are inverses of each other. Hence, 
it is enough to prove that $F_g\simeq F_{\UN}$ as left $\hat A$-modules 
for $g\in G_R(A)$ ($G_L(A),G(A)$) iff $g\in A^T$.

Let $g\in G^T_R(A):=G_R(A)\cap A^T$. Then $(A^T)^\perp:=\{\varphi\in\hat A
\vert\langle\varphi, A^T\rangle=0\}\subset\hat A$ is an ideal contained 
in the annihilator
ideal of both of the left $\hat A$-modules $F_{\UN}$ and $F_g$, because 
$F_{\UN}, F_g\subset A^T$ and $A^T$ is a subcoalgebra of $A$. Therefore 
$F_{\UN}$ and $F_g$ are also left modules with respect to the factor 
algebra $\hat A/(A^T)^\perp$ and the isomorphism of the modules $F_{\UN}$ and 
$F_g$ with respect to this factor algebra ensures their isomorphism as 
$\hat A$-modules. The factor algebra $\hat A/(A^T)^\perp$ is isomorphic to the 
dual WHA $\widehat{A^T}$ of $A^T$ as an algebra, which is isomorphic to a 
direct sum of simple matrix algebras, $\widehat{A^T}\simeq \oplus_\alpha
M_{n_\alpha}(Z_\alpha)$, due to Lemma 2.3. The $Z_\alpha$s are separable 
field extensions of the ground field $k$ determined by the ideal decomposition 
$Z=\oplus_\alpha Z_\alpha$ of $Z\equiv A^L\cap A^R$ and the dimensions
obey $n_\alpha={\rm dim}_{Z_\alpha} A^L_\alpha$. Hence, $F_{\UN}$ and $F_g$ 
are isomorphic $\widehat{A^T}$-modules if the multiplicities of simple 
submodules corresponding to the Wedderburn components of $\widehat{A^T}$ 
in their direct sum decompositions are equal. In order to prove this, 
first we note that the primitive idempotents $\{ z_\alpha\}_\alpha\subset Z$ 
are central in $A^T$, hence they are in the hypercenter $H$ of $A^T$ and 
they are related to the primitive central idempotents 
$\{\hat e_\alpha\}_\alpha$ of $\widehat{A^T}$ as 
$$
  \hat e_\alpha\Sr\UN=z_\alpha=\UN\Sl\hat e_\alpha\eqno(5.24)
$$
due to (1.6) and the remarks after it. 
Hence, $\hat e_\alpha\Sr g=(\hat e_\alpha\Sr\UN)g=z_\alpha g$ and
$F_{\UN}$ and $F_g$ are faithful left $\widehat{A^T}$-modules, because $\UN$
and $g$ are invertible. Therefore, the multiplicity corresponding to a 
Wedderburn component of $\widehat{A^T}$ is at least one in both of the modules
$F_{\UN}$ and $F_g$. Then the identity 
$$
  \vert F_\UN\vert =\vert F_g\vert=\vert\hat A^R\vert=\vert A^L\vert
  =\sum_\alpha\vert Z_\alpha\vert {\dim}_{Z_\alpha} A^L_\alpha
  =\sum_\alpha\vert Z_\alpha\vert n_\alpha\eqno(5.25)
$$
for $k$-dimensions coming from the $\hat A^R$-freeness of invertible 
$\hat A$-modules and from the algebra structure of $\widehat{A^T}$ 
ensures that these multiplicities are equal to one,
that is $F_{\UN}$ and $F_g$ are isomorphic $\widehat{A^T}\simeq\hat A/
(A^T)^\perp$, hence isomorphic $\hat A$-modules.  
    
Conversely, let $g\in G_R(A)$ be such that there exists an isomorphism
$U\colon F_{\UN}\to F_g$ between the invertible left 
$\hat A$-modules $F_{\UN}$ and $F_g$. Using that $U$ is an $\hat A$-module 
map, we have 
$$\eqalign{
  U(\UN)^{(1)}\langle\varphi,U(\UN)^{(2)}\rangle&=\varphi\Sr U(\UN)
  =U(\varphi\Sr\UN)=U(\hat\Pi^L(\varphi)\Sr\UN)\cr
 &=\hat\Pi^L(\varphi)\Sr U(\UN)
  =U(\UN)^{(1)}\langle\hat\Pi^L(\varphi),U(\UN)^{(2)}\rangle\cr
 &=U(\UN)^{(1)}\langle\varphi,\Pi^L(U(\UN)^{(2)})\rangle
  =\UN^{(1)}U(\UN)\langle\varphi,\UN^{(2)}\rangle, \qquad\varphi\in\hat A,\cr}
  \eqno(5.26)
$$
that is $\Delta(U(\UN))=\UN^{(1)}U(\UN)\otimes\UN^{(2)}$, which ensures that 
$U(\UN)\in A^L$. 
Moreover, $U(\UN)$ is an $\hat A^{L/R}$ generator of $F_g$, because 
it is the image of the $\hat A^{L/R}$ generator $\UN\in F_{\UN}$. 
Hence, there exists an invertible element $\varphi^L\in\hat A^L$ such that 
$$
  g=\varphi^L\Sr U(\UN)=(\varphi^L\Sr\UN)U(\UN)\in A^RA^L=A^T.\eqno(5.27)
$$
The case of (left) grouplike elements can be proved similarly. 

iii) Let $f$ be an $\hat A^{L/R}$-generator of the invertible submodule
$F_f\subset{_{\hat A}A}$. If there is no right grouplike element in 
$F_f=\hat A^L\Sr f=A^Rf$, that is, due to i), there is no such element 
$g$ in $F_f$ that obeys $\Pi^R(g)=\UN$, let us define $g:=f\Sl\hat f_l\in A$ 
with $\hat f_l$ given in (5.18). 
Then $F_g:=\hat A\Sr(f\Sl\hat f_l)=(\hat A\Sr f)\Sl\hat f_l=A^Rf\Sl\hat f_l
=A^R(f\Sl\hat f_l)=A^Rg$ due to (1.4) and the maps
$$\eqalignno{
  \Sl\hat f_l\colon F_f&\to F_g,&(5.28a)\cr
  x_Rf&\mapsto x_Rf\Sl\hat f_l=x_R(f\Sl\hat f_l)=x_Rg\cr  
  \Sl(\bar\Pi^L(f)\Sr\hat\UN)\colon F_g&\to F_f&(5.28b)\cr
  x_Rg&\mapsto x_Rg\Sl(\bar\Pi^L(f)\Sr\hat\UN)=x_Rg\bar\Pi^L(f)=x_Rf\cr}
$$
where $x_R\in A^R$, commute with the left Sweedler action, 
i.e. they are left $\hat A$-module maps.
They are also inverses of each other due to (5.22), which property has been 
already indicated in (5.28b). Therefore, $F_g$ and $F_f$ are equivalent 
submodules of ${_{\hat A}A}$, that is $F_g$ is also invertible. Since 
$\Pi^R(g):=\Pi^R(f\Sl\hat f_l)=\UN$ due to (5.18) and due to the 
non-degeneracy of the $A^R-\hat A^R$ pairing, $g$ is a right grouplike element 
due to i). The proof is similar for left grouplike elements: one has to 
define $g:=f\Sl\hat f_r$ with $\hat f_r$ given in (5.18) to get $g\in G_L(A)$
in the submodule $F_g$ isomorphic to $F_f$.\qed

\medskip\noindent{\bf Corollary 5.6} {\sl The elements of 
$G^T_{R/L}(A):=G_{R/L}(A)\cap A^T$ are of the form 
$g_LS(g_L^{-1})\in G^T_R(A)$ and $g_LS^{-1}(g_L^{-1})\in G^T_L(A)$. They are
in $G^T(A):=G(A)\cap A^T$ iff $g_L=S^2(g_L)$. 

\noindent
$G_{R/L}^T(A)$ and $G^T(A)$ is a normal subgroup in $G_{R/L}(A)$ and $G(A)$,
respectively.}

\smallskip\noindent{\it Proof.} An element $g\in G^T_R(A)$ has the product 
form $g=g_Lg_R$ due to (5.27) with $g_L:=U(\UN)\in A^L$ and $g_R:=
\varphi^L\Sr\UN\in A^R$. Since $g$ is invertible, $g_L$ and $g_R$ are 
invertible. Using property (5.2a) one obtains $\UN=\Pi^R(g)\equiv\Pi^R(g_Lg_R)
=g_RS(g_L)$. The other cases follow since 
$G^T_L(A)=S(G^T_R(A))$ and since $G^T(A)=G^T_R(A)\cap G^T_L(A)$. 

Since $\hat A\Sr g\Sl\hat A=gA^T=A^Tg$ for $g\in G_{R/L}(A)$ ($G(A)$) 
due to (5.1), $gA^Tg^{-1}=A^T$ follows. Therefore, $G_{R/L}^T(A)$ and 
$G^T(A)$ are normal subgroups.\qed

\medskip\noindent{\bf Proposition 5.7} {\sl Every invertible left 
$\hat A$-module is isomorphic to a cyclic submodule of 
$({_{\hat A}A},\Sr)$ generated by an element in $G_R(A)$ 
($G_L(A)$). The isomorphism classes of invertible left $\hat A$-modules  
are in one-to-one correspondence with elements of the (finite) factor group 
$G_R(A)/G_R^T(A)$ ($G_L(A)/G_L^T(A)$).} 

\smallskip\noindent{\it Proof.} Due to Prop. 5.4 ii) an invertible left 
$\hat A$-module $M$ is a direct sum of inequivalent simple submodules: 
$M=\oplus_p\hat z^L_p\cdot M=:\oplus_p M_p$, where $\{\hat z^L_p\}_p$ is the
set of primitive orthogonal idempotents in $\hat Z^L$. 
Since $\hat A$ is a quasi-Frobenius algebra, see Corollary 3.3, the simple 
submodules $M_p$ are isomorphic to left ideals in $\hat A$ [5, p.401]. 
Since they are inequivalent for different $p$, the invertible module 
${_{\hat A}M}$ itself is isomorphic to a left ideal in $\hat A$. Due to 
Corollary 4.2 $\hat A$ is a Frobenius algebra, hence, the isomorphism 
${_{\hat A}\hat A}\simeq ({_{\hat A}A},\Sr)$ of left regular modules 
holds [5, p.413]. Thus, ${_{\hat A}M}$ is isomorphic to an invertible submodule of $({_{\hat A}A},\Sr)$, that is to a cyclic submodule $F_g$ with $g\in G_R(A)$ ($g\in G_L(A)$) by Lemma 5.5 iii). Due to Lemma 5.5 ii) the isomorphism classes of cyclic submodules $F_g, g\in G_{R/L}(A)$ are given by the elements of the 
factor group $G_{R/L}(A)/G^T_{R/L}(A)$.  

Since a finite dimensional $k$-algebra has a finite number of inequivalent 
simple modules, there is only a finite number of inequivalent semisimple modules 
with a given $k$-dimension. Therefore, the factor groups
$G_{R/L}(A)/G_{R/L}^T(A)$ are finite groups.\qed
\medskip

In consideration of Prop. 5.7 we can formulate why the notion of grouplike 
elements in a WHA is too restrictive: one cannot always associate a grouplike 
element in $A$ to an invertible module of the dual WHA $\hat A$. We 
formulate this claim as follows:

\smallskip\noindent{\bf Proposition 5.8} {\sl
Let $t_L\in A^L_*$ denote the element that relates the counit and the 
reduced trace as non-degenerate functionals on the separable algebra $A^L$: 
$\varepsilon(\cdot)={\rm tr}\, (\cdot\, t_L)$. The coset $gG^T_R(A)\subset 
G_R(A)$ for $g\in G_R(A)$ contains a grouplike element iff there exists 
$x_L\in A^L_*$ such that 
$$
  gt_Lg^{-1}=x_Lt_Lx_L^{-1}.\eqno(5.29)
$$
In general, $G(A)/G^T(A)$ is a proper subgroup of $G_R(A)/G^T_R(A)$.}

\smallskip\noindent{\it Proof.} The adjoint action by $g\in G_R(A)$ 
on $A$ gives rise to 
algebra automorphisms of $A^L$ and $A^R$, because (5.1--2a) imply
that $\Pi^{R/L}(gy_{R/L}g^{-1})=gy_{R/L}g^{-1}$ for $y_{R/L}\in A^{R/L}$. 
Using the invariance of the reduced trace with respect to algebra 
automorphisms and the WBA identity 
$\varepsilon(abc)=\varepsilon(\Pi^R(a)b\Pi^L(c)); a,b,c\in A$, 
which follows from (1.1b) and (1.3), one obtains
$$\eqalign{
  \varepsilon(y_LgS(g))&=\varepsilon(\Pi^R(g^{-1})y_L\Pi^L(g))
  =\varepsilon(g^{-1}y_Lg)={\rm tr}\,(g^{-1}y_Lgt_L)\cr
  &={\rm tr}\,(y_Lgt_Lg^{-1})=\varepsilon(y_Lgt_Lg^{-1}t_L^{-1}),
   \quad y_L\in A^L,\cr}\eqno(5.30)
$$
i.e. $gS(g)=gt_Lg^{-1}t_L^{-1}$ due to non-degeneracy of the counit 
on $A^L$. Therefore, for all $g\in G_R(A)$ we have
$$
  S(g)=t_Lg^{-1}t_L^{-1},\qquad S^2(g)=tgt^{-1},\quad t:=t_LS(t_L^{-1}).
  \eqno(5.31)
$$
The element $t_L$ implements the Nakayama automorphism $\theta_\varepsilon=S^2$ of $\varepsilon$ on $A^L$: $\theta_\varepsilon={\rm Ad}\, t_L$. Hence, 
$t:=t_LS(t_L^{-1})\in A^T$ implements $S^2$ on $A^T$ and due to (5.31) on 
the subcoalgebras $gA^T$ of $A$, $g\in G_R(A)$ as well. In addition, 
$t\in G^T(A)$ due to Corollary 5.6.

Hence, if for a given $g\in G_R(A)$ there exists $x_L\equiv x_L(g)\in A^L_*$ 
such that $gt_Lg^{-1}=x_Lt_Lx_L^{-1}$, then 
$gS(g)=gt_Lg^{-1}t_L^{-1}=x_Lt_Lx_L^{-1}t_L^{-1}=x_LS^2(x_L^{-1})$ due to 
(5.31). Therefore, $h:=x_L^{-1}S(x_L)g\in G_R(A)$ is a grouplike element in 
the coset $gG^T_R(A)$ because $\Pi^L(h)=\UN$.

Conversely, if $h$ is a grouplike element in the coset $gG^T_R(A)\subset 
G_R(A)$ then $h=x_LS(x_L^{-1})g$ for some $x_L\in A^L_*$ due to 
Corollary 5.6. Therefore, using (5.31) 
$$
  \UN=\Pi^L(h)=x_LgS(g)S^2(x_L^{-1})=x_Lgt_Lg^{-1}t_L^{-1}S^2(x_L^{-1})
     =x_Lgt_Lg^{-1}x_L^{-1}t_L^{-1}.\eqno(5.32)
$$

For the second statement of the proposition first we note that the inclusion
$gG^T(A)\subset gG^T_R(A)$ for $g\in G(A)$ induces the inclusion
$G(A)/G^T(A)\subset G_R(A)/G^T_R(A)$ of the factor groups. To show that 
this inclusion is proper in general an example will suffice.

Let the WHA $A$ over the rational field ${\bf Q}$ be given as follows.
Let $A^L$ be a full matrix algebra $M_m({\bf Q}(\sqrt{2})), m>1$, where 
${\bf Q}(\sqrt{2})$ denotes the (separable) field extension of ${\bf Q}$ 
by $\sqrt{2}$.
Let the counit $\varepsilon$ as a non-degenerate index $\UN$ functional on 
the separable algebra $A^L$ be given with the help of the reduced trace: 
$\varepsilon(\cdot):={\rm tr}\,(\cdot t_L)$, where 
$t_L\in A^L_*$ satisfying ${\rm tr}\,(t_L^{-1})=1$. Let $A^T$ be the WHA
of the form $A^L\otimes A^{Lop}=: A^L\otimes A^R$ given in the Appendix 
of [2]. Let $A$ as an algebra over ${\bf Q}$ be given by the crossed 
product $A:=A^T\cros Z_2$, where $Z_2=\{e,g\}$ is the cyclic group of order 
two and the action of the non-trivial element $g\in Z_2$ on $A^L$ ($A^R$)
is the outer automorphism that changes the sign of the central element 
$z_L=\sqrt{2}\cdot\UN$ of $A^L$ ($z_R=\sqrt{2}\cdot\UN\in A^R$). 
Now it is a straightforward calculation that one extends the WHA structure 
of $A^T$ to $A:=A^T\cros Z_2$ by defining
$$\eqalign{
  \tilde\varepsilon(g^nx)&:=\epsilon(x),\cr
  \tilde\Delta(g^nx)&:=(g^n\otimes g^n)\Delta(x),\cr
  \tilde S(g^nx)&:=S(x)t_Lg^nt_L^{-1},\cr}\eqno(5.33)
$$
where $x\in A^T$ and $n=0,1$.

Due to Corollary 5.2 $g\in A$ becomes a right grouplike element for
any possible choice of $t_L$, i.e. $G_R(A)/G^T_R(A)\simeq Z_2$. 
However, if $t_L\in A^L_*$ is such that the prescribed 
outer automorphism on $A^L$ induced by $g$ is not inner on $t_L$, 
that is (5.29) is not fulfilled,
there is no grouplike element in the coset $gG^T_R(A)\subset G_R(A)$; 
thus, $G(A)/G^T(A)\simeq \{ e\}$.\qed

\medskip\noindent{\bf Corollary 5.9} {\sl $G_R(A)/G^T_R(A)=G(A)/G^T(A)$ if
$A^L$ is central simple or if $S^2_{\vert A^L}=id_{\vert A^L}$. In the latter 
case even $G_{R/L}(A)=G(A)$ holds.}

\smallskip\noindent{\it Proof.} If $A^L$ is central simple (5.29) is 
fulfilled by definition. In the other case $t_L$ 
is central in $A^L$, $G^T_R(A)=G^T(A)$ and 
(5.29) reads as $gt_Lg^{-1}=t_L, g\in G_R(A)$. Due to (5.31) 
$S(g)g=t_Lg^{-1}t_L^{-1}g$, and it is a central element in $A^L$.
Therefore, $\UN=\Pi^R(g)=S(\UN^{(1)})S(g)g\UN^{(2)}=S(g)g$ 
due to (5.1--2a), which proves the claim.\qed 

\bigskip\noindent
{\bf 6. Distinguished (left/right) grouplike elements, Radford formula 
and the order of the antipode}\medskip

After defining distinguished (left/right) grouplike elements and deriving some 
basic properties of them we prove the generalization of the Radford formula: 
the fourth power of the antipode in a WHA can be expressed in terms of 
distinguished left (right) grouplike elements like in the finite dimensional 
Hopf case [15]. Using this result we derive a finiteness type claim about 
the order of the antipode in a WHA and prove that the double of a 
WHA is unimodular.

We note that the Radford formula was proved in [13] for WHAs in the  
case when the square of the antipode is the identity mapping on 
$A^L$.\footnote{${}^1$}{For WHAs based on certain separable, 
but not strongly separable [9] algebra $A^L$ the property 
$S^2_{\vert A^L}\not=id_{\vert A^L}$, i.e. the non-triviality of the 
Nakayama automorphism corresponding to the counit as a non-degenerate 
functional $\varepsilon\colon A^L\to k$, is not only a possibility, but the
only possibility because $\varepsilon$ should be an index $\UN$ 
functional on $A^L$. For example, if $A^L=M_2({\bf Z}_2)$, that is a two 
by two matrix algebra over the finite field ${\bf Z}_2$, the reduced trace 
${\rm tr}$ on $A^L$ is non-degenerate but it has index $0$. The two 
non-degenerate index $\UN$ functional on $A^L$ have the form ${\rm tr}\, 
(\cdot\, t_L)$ with $t_L^{\pm 1}=\left(\matrix{1&1\cr1&0\cr}\right)$ and lead
to $S^2_{\vert A^L}={\rm Ad}\, t_L\not=id_{\vert A^L}$.} 
For such WHAs the sets of various grouplike elements
coincide, see Corollary 5.9. 

Before turning to the definition of (left/right) distinguished grouplike 
elements in a WHA let us examine the connection between integrals in dual 
pairs $A, \hat A$ of WHAs. 

The pair $(l,\lambda)\in I^L\times\hat I^L\subset A\times\hat A$ 
($(r,\rho)\in I^R\times\hat I^R$) is called a {\it dual pair of left 
(right) integrals} if they are non-degenerate and if they obey one of the 
equivalent relations $l\Sr\lambda=\hat\UN, \lambda\Sr l=\UN$ 
($r\Sl\rho=\UN, \rho\Sl r=\hat\UN$). 
Due to Theorem 4.1 such pairs exist in any dual pair of WHAs. 
$({_A{\hat I^L}},\star)$ is an invertible $A$-module due to Corollary 3.5 and
Prop. 5.4 i). Since this module is the right conjugate of the 
module ${_AI^R}$ due to Corollary 3.4, ${_AI^R}$ is also an invertible 
left $A$-module due to (5.7--8b). Hence, it is a free rank one 
left $A^{L/R}$-module due to Prop. 5.4 i).
An element $r$ is a free $A^L$ ($A^R$) generator in ${_AI^R}$ iff $r$ is 
a non-degenerate right integral, thus non-degenerate right 
integrals $r,r'\in I^R$ are related by an element 
$x_L\in A^L_*$ ($x_R\in A^R_*$): $r'=x_Lr$ ($r'=x_Rr$). The corresponding 
statement holds for non-degenerate right integrals in $\hat I^R$ by duality. 
Hence dual pairs of right integrals, $(r_1,\rho_1)$ and $(r_2,\rho_2)$, 
are related by a `common' invertible element $x_L\in A^L_*$ ($x_R\in A^R_*$):
$$
  (r_2,\rho_2)=(x_Lr_1,(\hat\UN\Sl x_L^{-1})\rho_1)
  =(x_Rr_1,(S^{-2}(x_R^{-1})\Sr\hat\UN)\rho_1).\eqno(6.1)
$$
Let us consider the element $s_R:=\rho\Sr r\in A$ constructed from the elements
of a dual pair $(r,\rho)$ of right integrals. Since $r$ is a non-degenerate 
functional on $\hat A$ and since $\rho$ is a free $\hat A^{L/R}$-generator 
of the left $\hat A$-module ${_{\hat A}{\hat I^R}}$, $s_R$ becomes a free 
left $\hat A^{L/R}$-generator of the cyclic left $\hat A$-module 
$(\hat A\Sr s_R,\Sr)$, i.e. it is an invertible $\hat A$-submodule 
in $(A,\Sr)$. Moreover, using (1.8)
$$\eqalign{
  \Pi^R(s_R)&:=\Pi^R(\rho\Sr r)=\Pi^R(r^{(1)})\langle r^{(2)},\rho\rangle
  =\UN^{(1)}\langle r\UN^{(2)},\rho\rangle\cr
 &=\UN^{(1)}\langle\UN^{(2)},\rho\Sl r\rangle
  =\UN^{(1)}\langle\UN^{(2)},\hat\UN\rangle=\UN,\cr}\eqno(6.2)
$$
that is $s_R$ is a right grouplike element in $A$ due to Lemma 5.5 i).
If $(r_i,\rho_i); i=1,2$ are dual pairs of right integrals the corresponding
right grouplike elements differ by a right grouplike element in $A^T$ due to
(6.1), (1.5--6) and Corollary 5.6:
$$
  \rho_2\Sr r_2=(\hat\UN\Sl x_L^{-1})\rho_1\Sr x_Lr_1
  =x_LS(x_L^{-1})(\rho_1\Sr r_1),\qquad x_L\in A^L_*.\eqno(6.3)
$$
However, it is not known to us whether the coset $G_R^T(A)s_R$ in 
$G_R(A)$ is special enough in order to contain always a grouplike element.
But we note that if $s_R:=\rho\Sr r$ is grouplike, 
i.e. $\Pi^L(s_R)=\UN$ also holds, 
then $\sigma_R:=r\Sr\rho\in G(\hat A)$ already follows: by duality $\sigma_R$ 
is a free $A^{L/R}$-generator in the cyclic left $A$-module 
$(A\Sr\sigma_R,\Sr)$ with the property $\hat\Pi^R(\sigma_R)=\hat\UN$ and
$$\eqalign{
   \hat\Pi^L(\sigma_R)&:=\hat\Pi^L(r\Sr\rho)
  =\hat\Pi^L(\rho^{(1)})\langle\rho^{(2)},r\rangle
  =\hat S(\hat\UN^{(1)})\langle\hat\UN^{(2)}\rho,r\rangle\cr 
 &=\hat S(\Pi^L(\rho\Sr r)\Sr\hat\UN)=:\hat S(\Pi^L(s_R)\Sr\hat\UN)
  =\hat S(\UN\Sr\hat\UN)=\hat\UN,\cr}
  \eqno(6.4)
$$
that is $\sigma_R$ is grouplike by Lemma 5.5 i).

Similarly, a dual pair $(l,\lambda)$ of left integrals leads to left 
grouplike elements: $s_L:=l\Sl\lambda\in G_L(A)$ and $\sigma_L:=\lambda\Sl l\in
G_L(\hat A)$. $s_L$ is grouplike iff $\sigma_L$ is grouplike, because 
$\Pi^R(s_L)$ and $\hat\Pi^R(\sigma_L)$ obey a relation analogous to (6.4):
$$
  \Pi^R(s_L)=S(\UN\Sl\hat\Pi^R(\sigma_L))=S(\UN\Sl\sigma_L)
            =\Pi^R(\UN\Sl\sigma_L).\eqno(6.5)
$$
These considerations lead to the following 
  
\medskip\noindent
{\bf Definition 6.1} {\sl Let $(l,\lambda)$ ($(r,\rho)$) be dual pair of left
(right) integrals in a dual pair $A,\hat A$ of WHAs. 
The elements $s_L:=l\Sl\lambda$ and $\sigma_L:=\lambda\Sl l$ ($s_R:=\rho\Sr r$ 
and $\sigma_R:=r\Sr\rho$) are called distinguished left (right) grouplike 
elements in $A$ and $\hat A$, respectively. 

\noindent
A dual pair of left (right) integrals is called 
a distinguished pair of left (right) integrals if $s_L$ ($s_R$) is not only 
left (right) grouplike but also grouplike. In this case $s_L$ ($s_R$) is
called distinguished grouplike element.}

\smallskip 
Let us introduce some notations we use in the forthcoming Lemma.
Using properties (5.1--2) it is easy to see that a left/right grouplike 
element $\gamma_{L/R}\in G_{L/R}(\hat A)$ gives rise to a projection 
$\Pi^{L/R}_{\gamma_{L/R}}\colon A\to A^{L/R}$ by defining  
$$
  \Pi^L_{\gamma_L}(a):=\Pi^L(\gamma_L\Sr a), \qquad 
  \Pi^R_{\gamma_R}(a):=\Pi^R(a\Sl\gamma_R), 
  \quad a\in A.\eqno(6.6)
$$
The invertible right/left $A$-module structures of left/right integrals 
in $A$ can be made explicit by using these projections and distinguished 
left/right grouplike elements $\sigma_{L/R}$ connected to the dual pair
$(l,\lambda)/(r,\rho)$ of left/right integrals:
$$
  la=l\Pi^R_{\hat S(\sigma_L^{-1})}(a),\qquad 
  ar=\Pi^L_{\hat S(\sigma_R^{-1})}(a)r,
  \quad a\in A.\eqno(6.7)
$$
For example, the first relation can be proved by using (5.1--2b), (1.6) and 
the non-degeneracy of $\lambda$: 
$$\eqalign{
  \lambda\Sl la&=\sigma_L\Sl a
                =\langle\hat S(\sigma_L^{-1})\hat\UN^{(1)},a\rangle
                 \sigma_L\hat\UN^{(2)}
                =\sigma_L(\hat\UN\Sl\Pi^R(a\Sl\hat S(\sigma_L^{-1})))\cr
               &=\sigma_L\Sl\Pi^R(a\Sl\hat S(\sigma_L^{-1}))
                =\lambda\Sl l\Pi^R_{\hat S(\sigma_L^{-1})}(a).\cr}
$$
\smallskip\noindent
{\bf Lemma 6.2} {\sl Let $C_b:=AbA\subset A$ be 
the cyclic ideal with the generator $b=b(\gamma,\delta)\in A$ 
characterized by a left and a right grouplike element $\gamma\in G_L(\hat A)$
and $\delta\in G_R(\hat A)$, respectively, through the property
$$
  abc=\Pi^L_\gamma(a)b\Pi^R_\delta(c),\qquad a,c\in A,\eqno(6.8)
$$
where the projections $\Pi^L_\gamma$ and $\Pi^R_\delta$ are defined in
(6.6). The left/right Sweedler actions by left/right grouplike elements 
in $\hat A$ provide isomorphisms between such types of cyclic ideals 
as (possibly non-unital) rings. The image $\tilde b$ of the generator 
$b=b(\gamma,\delta)$ obeys the characterization property
$$
   \beta_L\Sr b\Sl \beta_R=:\tilde b
  =\tilde b(\hat S(\beta_R)\gamma\beta_L^{-1},
   \beta_R^{-1}\delta\hat S(\beta_L)),
   \qquad \beta_{L/R}\in G_{L/R}(\hat A).\eqno(6.9)
$$}  
\noindent{\it Proof.} First, we note that the set of such cyclic ideals 
is non-empty: $l\in I^L$ from a dual pair 
$(\l,\lambda)$ of left integrals is a generator with characterization property 
$l=l(\hat\UN,\hat S(\sigma_L^{-1}))$ due to (1.9) and (6.7), where 
$\sigma_L:=l\Sl\lambda$ is the corresponding distinguished left 
grouplike element.

Since left (right) Sweedler actions by left (right) grouplike elements in 
$\hat A$ provide algebra automorphisms of $A$, the isomorphism of the 
corresponding 
cyclic ideals as rings follows. The only open question is the characterization 
property (6.9) of the image $\tilde b$ of the generator $b=b(\gamma,\delta)$. 
Using properties (6.1--2b) of left grouplike elements, characterization 
property (6.8) of the generator $b$, coproduct properties (1.4) of elements in 
$A^{L/R}$ and properties (1.7) of the projections $\Pi^{L/R}$ and 
$\hat\Pi^{L/R}$, one derives
$$\eqalignno{
  a(\beta_L&\Sr b)=\beta_L\Sr(\beta_L^{-1}\Sr a)b
  =\beta_L\Sr\Pi^L_\gamma(\beta_L^{-1}\Sr a)b
  =\Pi^L_\gamma(\beta_L^{-1}\Sr a)(\beta_L\Sr b)\cr
 &=:\Pi^L(\gamma\Sr(\beta_L^{-1}\Sr a))(\beta_L\Sr b)
  =\Pi^L_{\gamma\beta_L^{-1}}(a)(\beta_L\Sr b),&(6.10a)\cr
  (\beta_L&\Sr b)c=\beta_L\Sr b(\beta_L^{-1}\Sr c)
  =\beta_L\Sr b\Pi^R_\delta(\beta_L^{-1}\Sr c)\cr
 &=(\beta_L\Sr b)(\beta_L\Sr\Pi^R_\delta(\beta_L^{-1}\Sr c))
  =(\beta_L\Sr b)\UN^{(1)}\langle\beta_L,\Pi^R_\delta(\beta_L^{-1}\Sr c)
   \UN^{(2)}\rangle\cr
 &=(\beta_L\Sr b)\UN^{(1)}
   \langle\hat\UN^{(1)}\beta_L\otimes\hat\UN^{(2)}\beta_L,
   \Pi^R(\beta_L^{-1}\Sr c\Sl\delta)\otimes\UN^{(2)}\rangle\cr
 &=(\beta_L\Sr b)\UN^{(1)}
   \langle\hat\Pi^R(\hat\UN^{(1)}\beta_L)\otimes\hat\Pi^L(\hat\UN^{(2)}\beta_L),
   \beta_L^{-1}\Sr c\Sl\delta\otimes\UN^{(2)}\rangle\cr
 &=(\beta_L\Sr b)\UN^{(1)}
   \langle\hat S(\beta_L)\hat\UN^{(1)}\beta_L\otimes \hat\UN^{(2)},
   \beta_L^{-1}\Sr c\Sl\delta\otimes\UN^{(2)}\rangle\cr
 &=(\beta_L\Sr b)\UN^{(1)}
   \langle\Delta(\hat \UN),c\Sl\delta\hat S(\beta_L))\otimes\UN^{(2)}\rangle
  =(\beta_L\Sr b)\UN^{(1)}\varepsilon((c\Sl\delta\hat S(\beta_L))\UN^{(2)})\cr
 &=(\beta_L\Sr b)\Pi^R(c\Sl\delta\hat S(\beta_L))
  =(\beta_L\Sr b)\Pi^R_{\delta\hat S(\beta_L)}(c).&(6.10b)\cr}
$$
The change of the characterization property of the generator $b$ 
due to right Sweedler actions
$b\Sl\beta_R, \beta_R\in G_R(\hat A)$ can be proved similarly.\qed

\medskip\noindent {\bf Corollary 6.3} {\sl Distinguished left grouplike 
elements in $\hat A$ fall into a central element of the factor group 
$G_L(\hat A)/G^T_L(\hat A)$. There exists a two-sided non-degenerate 
integral in $A$ iff distinguished left grouplike elements in $\hat A$
fall into the unit element of this factor group.}

\smallskip\noindent{\it Proof.} For any $\beta\in G_L(\hat A)$ the map 
$B_\beta(a):=\beta\Sr a\Sl\hat S^{-1}(\beta), a\in A$
defines an algebra automorphism of $A$, which maps the space $I^L$ of 
left integrals into itself due to the previous Lemma. The image 
$\tilde l:= B_\beta(l)$ of a non-degenerate left integral 
$l=l(\hat\UN, \hat S(\sigma_L^{-1}))$ is a non-degenerate left integral 
having the characterization property $\tilde l=
\tilde l(\hat\UN,\hat S^{-1}(\beta^{-1})\hat S(\sigma_L^{-1})\hat S(\beta))$
due to (6.9). Hence, the distinguished left grouplike element $\tilde\sigma_L$
corresponding to $\tilde l$ is given by
$$
  \tilde\sigma_L=\hat S^{-2}(\beta)\sigma_L\beta^{-1}
                =:\varphi\beta\sigma_L\beta^{-1},\quad
  S^{-2}(\beta)\beta^{-1}=:\varphi=
               \hat S^{-1}(\varphi_L^{-1})\varphi_L\in G^T_L(\hat A),
  \eqno(6.11)
$$ 
with $\varphi_L=\hat S^{-1}(\hat\Pi^R(\beta^{-1}))\in\hat A^L_*$ due to 
the form (5.2b) of $\hat\Pi^R(\beta^{-1})$. However, distinguished 
left grouplike elements differ by elements in $G^T_L(\hat A)$, in analogy
with the case (6.3) of distinguished right grouplike elements. 
Hence, for the $G^T_L(\hat A)$-cosets (6.11) implies the relation
$[\sigma_L]=[\tilde\sigma_L]=[\varphi][\beta][\sigma_L][\beta]^{-1}
=[\beta][\sigma_L][\beta]^{-1}$, that is $[\sigma_L]$ is central in the factor
group $G_L(\hat A)/G^T_L(\hat A)$.

If the non-degenerate left integral $l\in I^L$ is also a right integral then
we have the relation $\Pi^R_{\hat S(\sigma_L^{-1})}=\Pi^R$ due to (6.7) and 
(1.9). Hence, $\sigma_L=\hat\UN$ since 
$$
  \langle\hat\UN, a\rangle=\langle\hat\UN,\Pi^R(a)\rangle    
         =\langle\hat\UN,\Pi^R_{\hat S(\sigma_L^{-1})}(a)\rangle   
         =\langle\hat\UN,a\Sl\hat S(\sigma_L^{-1})\rangle
         =\langle\hat S(\sigma_L^{-1}),a\rangle, \quad a\in A,
$$    
using (6.6) and (1.7). Conversely, if $[\sigma_L]$ is the unit element 
of the factor group then there exists a dual pair $(l,\lambda)$ of 
left integrals with distinguished left grouplike element $\sigma_L=\hat\UN$
due to a relation analogous with (6.3). Therefore, 
$\Pi^R_{\hat S(\sigma_L^{-1})}=\Pi^R$ and (6.7) implies that $l$ is a
(non-degenerate) two-sided integral.\qed

\medskip\noindent
{\bf Theorem 6.4} {\sl Let $A,\hat A$ be a dual pair of WHAs and let 
$(s_L,\sigma_L)$ be the pair of distinguished left grouplike elements
corresponding to a dual pair $(l,\lambda)$ of left integrals in 
$A\times\hat A$. The Nakayama automorphism 
$\theta_\lambda:=\hat R_\lambda^{-1}\circ \hat L_\lambda\colon A\to A$ 
corresponding to the non-degenerate functional $\lambda\colon A\to k$ can 
be written as
$$
  \theta_\lambda(a)=\sigma_L^{-1}\Sr S^2(a)
                   =s_L^{-1}S^{-2}(a)s_L\Sl\hat S^{-1}(\sigma_L),
  \qquad a\in A.\eqno(6.12)
$$
The fourth power of the antipode $S$ of $A$ can be written as:
$$
  S^4(a)=\sigma_L\Sr s_L^{-1}as_L\Sl\hat S^{-1}(\sigma_L),\qquad a\in A.
  \eqno(6.13)
$$
The order of the antipode is finite up to an inner automorphism by a 
grouplike element in the trivial subalgebra $A^T$.}
\smallskip
   
\noindent{\it Proof.} In the sufficiency proof of Theorem 4.1 we have seen 
that the antipode and its inverse can be given with the help of pairs of 
non-degenerate integrals $l/r\in I^{L/R}, \lambda / \rho\in\hat I^{L/R}$
$$\eqalignno{
  S(a)&=(R_l\circ\hat L_\lambda)(a):=(\lambda\Sl a)\Sr l,\qquad\quad 
        \lambda\Sr l=\UN, l\Sr\lambda=\hat\UN,&(6.14a)\cr
  S^{-1}(a)&=(L_l\circ\hat L_\rho)(a):=l\Sl(\rho\Sl a),\qquad\quad 
        l\Sl\rho=\UN, l\Sr\rho=\hat\UN&(6.14b)\cr
  S^{-1}(a)&=(R_r\circ\hat R_\lambda)(a):=(a\Sr\lambda)\Sr r,\qquad\quad 
        \lambda\Sl r=\hat\UN, \lambda\Sr r=\UN,&(6.14c)\cr
  S(a)&=(L_r\circ\hat R_\rho)(a):=r\Sl(a\Sr\rho),\qquad\quad 
        r\Sl\rho=\UN, \rho\Sl r=\hat\UN.&(6.14d)\cr}
$$
Choosing a dual pair $(l,\lambda)$ of left integrals, we rewrite the 
antipode relations (6.14b--d) in terms of $(l,\lambda)$ and the corresponding 
pair $(s,\sigma)\equiv(s_L,\sigma_L)$ of distinguished left grouplike 
elements. We note that the second relations 
between the members of integral pairs given in (6.14a--d) are consequences of 
the first ones (see the proof of Theorem 4.1), hence, it is enough to ensure 
only these ones.

For (6.14b) the new member of the required pair $(l,\rho)$ of integrals is 
given by $\rho:=\hat S^{-1}(\lambda)=(\lambda\Sl s)\hat\Pi^R(\sigma)^{-1}$. 
Indeed, $\rho$ is a non-degenerate right integral and 
$\lambda=\hat S(\rho)=(l\Sl\rho)\Sr\lambda$ implies the relation 
$l\Sl\rho=\UN$ due to injectivity of $\hat R_\lambda$. Moreover,
using property (1.16) of left integrals 
$$\eqalign{
  (\lambda\Sl s)\hat\Pi^R(\sigma)^{-1}
 &:=(\lambda\Sl (l\Sl\lambda))\hat\Pi^R(\sigma)^{-1}
  =\langle\lambda\lambda^{(1)},l\rangle\lambda^{(2)}\hat\Pi^R(\sigma)^{-1}\cr
 &=\langle\lambda^{(1)},l\rangle \hat S^{-1}(\lambda)\lambda^{(2)}
   \hat\Pi^R(\sigma)^{-1}
  =\hat S^{-1}(\lambda)\sigma\hat\Pi^R(\sigma)^{-1}\cr
 &=\hat S^{-1}(\lambda)\hat\Pi^R(\sigma)\hat\Pi^R(\sigma)^{-1}
  =\hat S^{-1}(\lambda)=\rho.\cr}\eqno(6.15)
$$ 
Hence, interchanging the role of $A$ and $\hat A$, the new member of 
integrals for (6.14c) is given by $r:=S^{-1}(l)=(l\Sl\sigma)\Pi^R(s)^{-1}$.
For (6.14d) the pair is given by 
$(r:=S^{-1}(l),\rho:=\hat S(\lambda)=s\Sr\lambda)$,
because $\rho=\hat S(\lambda)=(l\Sl\lambda)\Sr\lambda=s\Sr\lambda$ and 
$r=S^{-1}(l)$ are non-degenerate right integrals and 
$r\Sl\rho=S^{-1}(l)\Sl\hat S(\lambda)=S^{-1}(\lambda\Sr l)=\UN$.   
Therefore, we can rewrite (6.14b and c) as
$$\eqalignno{
 S^{-1}(a)&=l\Sl(\rho\Sl a)=l\Sl((\lambda\Sl s)\hat\Pi^R(\sigma)^{-1}\Sl a)
  =l\Sl(\lambda\Sl sa)\hat\Pi^R(\sigma)^{-1}\cr
 &=[l\Sl(\lambda\Sl sa)](\UN\Sl\hat\Pi^R(\sigma)^{-1})
  =(L_l\circ\hat L_\lambda)(sa)(\UN\Sl\hat\Pi^R(\sigma)^{-1}),&(6.16b)\cr
 S^{-1}(\sigma&\Sr a)=S^{-1}(a)\Sl\hat S(\sigma)
  =((a\Sr\lambda)\Sr r)\Sl\hat S(\sigma)\cr
 &=(a\Sr\lambda)\Sr (l\Sl\sigma)\Pi^R(s)^{-1}\Sl\hat S(\sigma)
  =(a\Sr\lambda)\Sr (l\Sl\sigma\hat S(\sigma))\Pi^R(s)^{-1}\cr
 &=(a\Sr\lambda)\Sr (l\Sl\hat\Pi^R(\sigma^{-1})^{-1})\Pi^R(s)^{-1}\cr
 &=(a\Sr\lambda)\Sr l(\UN\Sl\hat\Pi^R(\sigma^{-1})^{-1})\Pi^R(s)^{-1}\cr
 &=(a\Sr\lambda)\Sr l\Pi^R((\UN\Sl\sigma\hat S(\sigma))
   \Sl\hat S(\sigma^{-1}))\Pi^R(s)^{-1}\cr
 &=(a\Sr\lambda)\Sr l\Pi^R(\UN\Sl\sigma)\Pi^R(s)^{-1}
  =(a\Sr\lambda)\Sr l\Pi^R(s)\Pi^R(s)^{-1}\cr
 &=(a\Sr\lambda)\Sr l
  =(R_l\circ\hat R_\lambda)(a),&(6.16c)\cr}
$$
using relations (1.4) and (1.6) for elements in $A^R$, the identity 
$\sigma\hat S(\sigma)=\hat\Pi^R(\sigma^{-1})^{-1}$
following from (5.2b), the right $A$-module property (6.7) of left integrals 
and the relation (6.5). Finally, using property (1.16) of left integrals,
(6.14d) can be rewritten as
$$\eqalignno{
 S(a)&=r\Sl(a\Sr\rho)
  =(l\Sl\sigma)\Pi^R(s)^{-1}\Sl(a\Sr(s\Sr\lambda))\cr
 &=[(l\Sl\sigma)\Sl(as\Sr\lambda)]\Pi^R(s)^{-1}
  =[l\Sl\sigma\lambda^{(1)}\langle\lambda^{(2)},as\rangle]\Pi^R(s)^{-1}\cr
 &=[l\Sl\lambda^{(1)}\langle\hat S^{-1}(\sigma)\lambda^{(2)},as\rangle]
   \Pi^R(s)^{-1}
  =[l\Sl((as\Sl\hat S^{-1}(\sigma))\Sr\lambda)]\Pi^R(s)^{-1}\cr
 &=(L_l\circ\hat R_\lambda)(as\Sl\hat S^{-1}(\sigma))\Pi^R(s)^{-1}.&(6.16d)\cr}
$$
Therefore using (6.14a), (6.16b--d), the algebra isomorphism property of 
the map $\hat\kappa_R$ given in (1.5), the relation (6.5) and the form (5.2b) 
of $\Pi^R(s)$ we get 
$$\eqalignno{
  (R_l\circ \hat L_\lambda)(a)
 &=S(a)=S^{-1}(\sigma\Sr(\sigma^{-1}\Sr S^2(a)))
  =(R_l\circ\hat R_\lambda)(\sigma^{-1}\Sr S^2(a)),&(6.17a)\cr
  (L_l\circ\hat L_\lambda)(a)
 &=S^{-1}(s^{-1}a)(\UN\Sl\hat\Pi^R(\sigma))
  =S^{-1}(\Pi^R(s)s^{-1}a)=S^{-1}(a)s\cr
 &=S[s^{-1}S^{-2}(a)]\Pi^R(s)
  =(L_l\circ\hat R_\lambda)(s^{-1}S^{-2}(a)s\Sl\hat S^{-1}(\sigma)).
  &(6.17b)\cr}
$$
Due to injectivity of $R_l$ and $L_l$ (6.17a and b) lead to connections
between $\hat R_\lambda$ and $\hat L_\lambda$ that imply (6.12). The equality
of these two different forms of the Nakayama automorphism 
$\theta_\lambda$ gives rise to the Radford formula (6.13).

Since left (right) Sweedler actions by left (right) grouplike elements are 
algebra automorphisms, iterating the Radford formula $m$ times one arrives at
$$
  S^{4m}(a)=S^{4m}(s^{-1})\ldots S^4(s^{-1})
            (\sigma^{m}\Sr a\Sl\hat S^{-1}(\sigma^m))
            S^4(s)\ldots S^{4m}(s),\qquad a\in A.
  \eqno(6.18)
$$
For $g\in G_L(A)$ the relation 
$S^2(g)=S(\Pi^R(g^{-1})^{-1})\Pi^R(g^{-1})g\in G_L^T(A)g$ holds due to (5.2b)
and Corollary 5.6. Hence, $S^{2n}(g)\in G_L^T(A)g$ is for any integer
$n$. Since the factor group $G_L(A)/G_L^T(A)$ is finite due to 
Prop. 5.7, there exists an integer $m$ and $x\equiv S(x_R)x_R^{-1}\in G_L^T(A), 
\varphi\equiv\hat S(\varphi_R)\varphi_R^{-1}\in G_L^T(\hat A)$ with 
$x_R\in A^R_*,\varphi_R\in\hat A^R_*$ such that (6.18) reads as
$$\eqalign{
  S^{4m}(a)&=x^{-1}(\varphi\Sr a\Sl\hat S^{-1}(\varphi))x\cr
           &=x^{-1}(\hat S(\varphi_R)\Sr\UN)
             (\UN\Sl\hat S^{-1}(\varphi_R^{-1})) a
             (\varphi_R^{-1}\Sr\UN)(\UN\Sl\varphi_R)x\cr
           &=x^{-1}S^{-1}(\UN\Sl\varphi_R)
             (\UN\Sl\varphi_R^{-1}) a
             S^{-1}(\UN\Sl\varphi_R^{-1})(\UN\Sl\varphi_R)x\cr
           &=x^{-1}S^{-1}(\UN\Sl\varphi_R)
             (\UN\Sl\varphi_R)^{-1} a
             S^{-1}(\UN\Sl\varphi_R)^{-1}(\UN\Sl\varphi_R)x\cr
           &=S(y_R^{-1})y_R a S(y_R)y_R^{-1},\qquad a\in A\cr}
  \eqno(6.19)
$$
where we used the identities (1.6) and the notation
$y_R:=x_RS^{-1}(\UN\Sl\varphi_R)\in A^R_*$. Due to (6.19) $S^{4m}$ is an inner
algebra automorphism of $A$ by an element $y:=S(y_R)y_R^{-1}\in G_L^T(A)$.
However, $S^{4m}$ is also a coalgebra automorphism of $A$, which requires $y$
to be a grouplike element. Indeed, using the coproduct property (1.4) and 
separability identities (1.12) for $A^L$ and $A^R$, one derives the relation 
$$\eqalign{
  \Delta(a)&=(S(y_R)\otimes y_R^{-1})\Delta(S(y_R^{-1})y_R a S(y_R)y_R^{-1}) 
             (S(y_R^{-1})\otimes y_R)\cr
           &=(S(y_R)\otimes y_R^{-1})\Delta(S^{4m}(a))
             (S(y_R^{-1})\otimes y_R)\cr
           &=(S(y_R)\otimes y_R^{-1})(S^{4m}\otimes S^{4m})(\Delta(a))
             (S(y_R^{-1})\otimes y_R)\cr
           &=(y_R\otimes S(y_R^{-1}))\Delta(a)
             (y_R^{-1}\otimes S(y_R))
            =(y_RS^2(y_R^{-1})\otimes\UN)\Delta(a)
             (y_Ry_R^{-1}\otimes\UN)\cr
           &=(y_RS^2(y_R^{-1})\otimes\UN)\Delta(a),\qquad a\in A,\cr}
  \eqno(6.20)
$$
which leads to the equality $\UN=y_RS^2(y_R^{-1})$ by applying the counit to 
the second tensor factor. Hence, $y=S(y_R)y_R^{-1}$ is not only in 
$G_L^T(A)$ but also in $G^T(A)$ due to Corollary 5.6, which together with
(6.19) proves the last claim in the theorem.\qed 

\smallskip
The Radford formula [15] was used in [16] to prove unimodularity of
the Drinfeld double $\cD(H)$ of a Hopf algebra $H$. 
In the case of the double $\cD(A)$ of a WHA $A$ [3] the same result
holds:
 
\medskip\noindent
{\bf Corollary 6.5} {\sl The double $\cD(A)$ of a WHA $A$ is unimodular,
i.e. there exists a non-degenerate two-sided integral in $\cD(A)$. Namely,
if $(l,\lambda)$ is a dual pair of left integrals in $A\times\hat A$ then 
$\cD(l\otimes\hat S(\lambda))$ is a two-sided non-degenerate 
integral in $\cD(A)$.}
\smallskip
   
\noindent{\it Proof.} The double $\cD(A)$ of a WHA $A$ [3] 
is the $k$-linear space of the tensor product of $A$ and $\hat A$ 
over the subalgebras $A^L\simeq\hat A^R$ and $A^R\simeq\hat A^L$ 
$$
  \cD(A)\ni\cD(ax_Lx_R\otimes\varphi)
   =\cD(a\otimes (x_L\Sr\hat\UN)(\hat\UN\Sl x_R)\varphi),\quad
   a\in A,\varphi\in\hat A, x_{L/R}\in A^{L/R}\eqno(6.21)
$$ 
together with the WHA structure maps
$$\eqalign{
  \cD(a\otimes\varphi)\cD(b\otimes\psi)&:=\cD(ab^{(2)}\otimes\varphi^{(2)}\psi)
  \langle a^{(1)},\hat S^{-1}(\varphi^{(3)})\rangle 
  \langle a^{(3)},\varphi^{(1)}\rangle,\cr
  \varepsilon_{\cD}(\cD(a\otimes\varphi))&:=\varepsilon(a(\varphi\Sr\UN))
               =\hat\varepsilon((\hat\UN\Sl a)\varphi),\cr 
  \Delta_\cD(\cD(a\otimes\varphi))&:=\cD(a^{(1)}\otimes\varphi^{(2)})\otimes
             \cD(a^{(2)}\otimes\varphi^{(1)}),\cr
   S_\cD(\cD(a\otimes\varphi))&:=\cD(\UN\otimes\hat S^{-1}(\varphi))
                                \cD(S(a)\otimes\hat\UN).\cr}\eqno(6.22)
$$
Let $(l,\lambda)$ be a dual pair of left integrals in $A\times\hat A$ with
the corresponding pair $(s,\sigma)\equiv(s_L,\sigma_L)$ of distinguished left 
grouplike elements. The expression (6.12) of the Nakayama automorphism 
$\theta_l$ corresponding to $l$ implies that 
$\Delta^{op}(l)=l^{(1)}\otimes S^2(l^{(2)}s_L^{-1})$.
Hence, using properties (5.1--2b) of a left grouplike element
$$\eqalignno{
  l^{(2)}\otimes &l^{(3)}s^{-1}S^{-1}(l^{(1)})
  =l^{(1)}\otimes l^{(2)}s^{-1}S(l^{(3)}s^{-1})
  =l^{(1)}\otimes l^{(2)}\UN^{(1)}s^{-1}S(s^{-1})S(\UN^{(2)})S(l^{(3)})\cr
 &=l^{(1)}\otimes l^{(2)}\Pi^L(s^{-1})S(l^{(3)})
  =l^{(1)}\otimes \Pi^L(l^{(2)}).&(6.23)\cr}
$$
Moreover, the first relation in (6.7) implies that for $\lambda\in\hat I^L$ 
and $\varphi\in\hat A$ 
$$
  \varphi\hat S(\lambda)
  =\hat S(\lambda\hat\Pi^R_{S(s^{-1})}(\hat S^{-1}(\varphi)))
  =(\hat S\circ\hat\Pi^R\circ\hat S^{-1})(s^{-1}\Sr\varphi)\hat S(\lambda)
  =\hat\Pi^L_{s^{-1}}(\varphi)\hat S(\lambda),
   \eqno(6.24)
$$
where we used (1.10) in the third equality.
The relations (6.23--24) and the properties (6.21--22) of $\cD(A)$ together
with the identities (1.6) and (1.11) lead to
$$\eqalignno{
  \cD(a&\otimes\varphi)\cD(l\otimes\hat S(\lambda))
  :=\cD(al^{(2)}\otimes\varphi^{(2)}\hat S(\lambda))
   \langle l^{(1)},\hat S^{-1}(\varphi^{(3)})\rangle
   \langle l^{(3)},\varphi^{(1)}\rangle\cr
 &=\cD(al^{(2)}\otimes\hat\Pi^L(s^{-1}\Sr\varphi^{(2)})\hat S(\lambda))
   \langle S^{-1}(l^{(1)}),\varphi^{(3)}\rangle
   \langle l^{(3)},\varphi^{(1)}\rangle\cr
 &=\cD(al^{(2)}(\hat\Pi^L(\varphi^{(2)})\Sr\UN)\otimes\hat S(\lambda))
   \langle l^{(3)},\varphi^{(1)}\rangle
   \langle s^{-1}S^{-1}(l^{(1)}),\varphi^{(3)}\rangle\cr
 &=\cD(a(\hat{\bar\Pi}^R(\varphi^{(2)})\Sr l^{(2)})\otimes\hat S(\lambda))
   \langle l^{(3)},\varphi^{(1)}\rangle
   \langle s^{-1}S^{-1}(l^{(1)}),\varphi^{(3)}\rangle\cr
 &=\cD(al^{(2)}\otimes\hat S(\lambda))
   \langle l^{(3)},\hat{\bar\Pi}^R(\varphi^{(2)})\varphi^{(1)}\rangle
   \langle s^{-1}S^{-1}(l^{(1)}),\varphi^{(3)}\rangle&(6.25)\cr
 &=\cD(al^{(2)}\otimes\hat S(\lambda))
   \langle l^{(3)}s^{-1}S^{-1}(l^{(1)}),\varphi\rangle\cr
 &=\cD(al^{(1)}\otimes\hat S(\lambda))
   \langle \Pi^L(l^{(2)}),\varphi\rangle
  =\cD(a(\hat\Pi^L(\varphi)\Sr l)\otimes\hat S(\lambda))\cr
 &=\cD(aS(\hat\Pi^L(\varphi)\Sr\UN)l\otimes\hat S(\lambda))
  =\cD(a(\UN\Sl\hat S^{-1}(\hat\Pi^L(\varphi)))l\otimes\hat S(\lambda))\cr
 &=\cD(\Pi^L(a(\UN\Sl\hat{\bar\Pi}^R(\varphi)))l\otimes\hat S(\lambda))
  =\cD(\Pi^L(a(\UN\Sl\hat{\bar\Pi}^R(\varphi)))\otimes\hat\UN)
    \cD(l\otimes\hat S(\lambda))\cr
 &=\Pi^L_\cD(\cD(a\otimes\varphi))\cD(l\otimes\hat S(\lambda)),
   \qquad a\in A,\varphi\in\hat A,\cr}
$$
that is $\cD(l\otimes\hat S(\lambda))$ is a left integral in $\cD(A)$.
A similar computation shows that it is also a right integral. 

Now, we prove that $\cD(l\otimes\hat S(\lambda))$ is a non-degenerate 
functional on the dual $\hat\cD(A)$ of $\cD(A)$. The WHA $\hat\cD(A)$
[3] is the $k$-linear space of the tensor product of $\hat A$ and $A$ 
over the subalgebras $\hat A^R\simeq A^L$ and $\hat A^L\simeq A^R$ 
$$
  \hat\cD(A)\ni\hat\cD(\varphi\otimes x_LaS^{-1}(x_R))
  =\hat\cD(\hat S^{-1}(\hat\UN\Sl x_R)\varphi(x_L\Sr\hat\UN)\otimes a),
  \eqno(6.26)
$$ 
where $\varphi\in\hat A,a\in A, x_{L/R}\in A^{L/R}$. The WHA 
structure maps of $\hat\cD(A)$ are transposed to that of $\cD(A)$ with 
respect to the non-degenerate pairing
$$
  \langle\hat\cD(\varphi\otimes a),\cD(b\otimes\psi)\rangle 
 :=\langle\varphi\otimes a, P(b\otimes\psi)\rangle
  =\langle\hat P(\varphi\otimes a), b\otimes\psi\rangle,\quad
  a,b\in A,\ \varphi,\psi\in\hat A,\eqno(6.27)
$$  
where $P\colon A\otimes\hat A\to  A\otimes\hat A$ and 
$\hat P\colon\hat A\otimes A\to\hat A\otimes A$ are $k$-linear projections 
given with the help of separating idempotents of $A^L$ and $A^R$
$$\eqalign{
  P(b\otimes\psi)&:=b\UN^{(1)}S(\UN^{(1')})\otimes
                 (\UN^{(2')}\Sr\hat\UN)(\hat\UN\Sl S(\UN^{(2)}))\psi,\cr
  \hat P(\varphi\otimes a)&:=
                (\UN^{(1')}\Sr\hat\UN)\varphi(\UN^{(1)}\Sr\hat\UN)\otimes
  \UN^{(2)}a\UN^{(2')}.\cr}\eqno(6.28)
$$
Clearly, $P(A\otimes\hat A)$ and $\cD(A)$ ($\hat P(\hat A\otimes A)$ and 
$\hat\cD(A)$) are isomorphic $k$-linear spaces and 
$\cD(P(b\otimes\psi))=\cD(b\otimes\psi)$
($\hat\cD(\hat P(\varphi\otimes a))=\hat\cD(\varphi\otimes a)$) also holds
due to (6.21) (or (6.26)). Since $\cD(A)$ is finite dimensional, the two-sided 
integral $\cD(l\otimes\hat S(\lambda))$ is non-degenerate if the $k$-linear map
$R_{\cD(l\otimes\hat S(\lambda))}\colon\hat\cD(A)\to\cD(A)$ is injective, 
that is if $0=\hat\cD(\varphi\otimes a)\Sr\cD(l\otimes\hat S(\lambda))$ implies
$0=\hat\cD(\varphi\otimes a)$. Using the mentioned isomorphisms of the 
$k$-linear spaces, the definition (6.28) of the projection $P$, 
the form of the coproduct in $\cD(A)$ and the identity 
$(\hat\UN\Sl\Pi^R_{\hat S(\sigma^{-1})}(S(\UN^{(1)})))
\hat\Pi^L_{s^{-1}}(\UN^{(2)}\Sr\hat\UN)=\hat\UN$
we prove later on one computes
$$\eqalignno{
  &\empty P(l^{(1)}\otimes\hat S(\lambda)^{(2)})
      \langle\hat\cD(\varphi\otimes a),
      \cD(l^{(2)}\otimes\hat S(\lambda)^{(1)})\rangle:=\cr
   &:=l^{(1)}\UN^{(1)}S(\UN^{(1')})\otimes
      (\UN^{(2')}\Sr\hat\UN)(\hat\UN\Sl S(\UN^{(2)}))\hat S(\lambda)^{(2)}
      \langle\hat\cD(\varphi\otimes a),
      \cD(l^{(2)}\otimes\hat S(\lambda)^{(1)})\rangle\cr
   &=l^{(1)}S(\UN^{(1')})\otimes
      (\UN^{(2')}\Sr\hat\UN)\hat S(\lambda)^{(2)}
      \langle\hat\cD(\varphi\otimes a),
      \cD(l^{(2)}S(\UN^{(1)})\otimes
      \hat S(\hat\UN\Sl S(\UN^{(2)}))\hat S(\lambda)^{(1)})\rangle\cr  
   &=l^{(1)}S(\UN^{(1')})\otimes
      (\UN^{(2')}\Sr\hat\UN)\hat S(\lambda)^{(2)}
      \langle\hat\cD(\varphi\otimes a),
      \cD(l^{(2)}S(\UN^{(1)})\UN^{(2)}\otimes
      \hat S(\lambda)^{(1)})\rangle\cr  
   &=[lS(\UN^{(1)})]^{(1)}\otimes
      [(\UN^{(2)}\Sr\hat\UN)\hat S(\lambda)]^{(2)}
      \langle\hat\cD(\varphi\otimes a),
      \cD([lS(\UN^{(1)})]^{(2)}\otimes
      [(\UN^{(2)}\Sr\hat\UN)\hat S(\lambda)]^{(1)})\rangle\cr    
   &=l^{(1)}\otimes\hat S(\lambda)^{(2)}
      \langle\hat\cD(\varphi\otimes a),
      \cD(l^{(2)}\Pi^R_{\hat S(\sigma^{-1})}(S(\UN^{(1)}))\otimes
      \hat\Pi^L_{s^{-1}}(\UN^{(2)}\Sr\hat\UN)\hat S(\lambda)^{(1)})\rangle\cr      &=l^{(1)}\otimes\hat S(\lambda)^{(2)}
      \langle\hat\cD(\varphi\otimes a),
      \cD(l^{(2)}\otimes(\hat\UN\Sl\Pi^R_{\hat S(\sigma^{-1})}(S(\UN^{(1)})))
      \hat\Pi^L_{s^{-1}}(\UN^{(2)}\Sr\hat\UN)\hat S(\lambda)^{(1)})\rangle\cr      &=l^{(1)}\otimes\hat S(\lambda)^{(2)}
      \langle\hat\cD(\varphi\otimes a),
      \cD(l^{(2)}\otimes\hat S(\lambda)^{(1)})\rangle\cr    
   &=l^{(1)}\otimes\hat S(\lambda)^{(2)}
      \langle\hat P(\varphi\otimes a),
       l^{(2)}\otimes\hat S(\lambda)^{(1)}\rangle
    =(R_l\otimes\hat L_{\hat S(\lambda)})(\hat P(\varphi\otimes a)),&(6.29)\cr}
$$
where we used (1.16) in the second equality, (1.4) in the fourth and 
(6.7) and (1.10) in the fifth one. The $k$-linear map 
$R_l\otimes\hat L_{\hat S(\lambda)}
\colon\hat A\otimes A\to A\otimes\hat A$ is injective due to the 
non-degeneracy of the integrals $l$ and $\lambda$. Hence, (6.29) implies that 
$\hat P(\varphi\otimes a)$, or equivalently $\hat\cD(\varphi\otimes a)$, 
should be zero if the left hand side of (6.29), or equivalently 
$\hat\cD(\varphi\otimes a)\Sr\cD(l\otimes\hat S(\lambda))$, is zero. 
Finally, the proof of the identity we used in (6.29) is as follows:
$$\eqalign{
  (\hat\UN&\Sl\Pi^R_{\hat S(\sigma^{-1})}(S(\UN^{(1)})))
  \hat\Pi^L_{s^{-1}}(\UN^{(2)}\Sr\hat\UN)
  =(\hat\UN\Sl\Pi^R(S(\UN^{(1)})\Sl\hat S(\sigma^{-1})))
  \hat\Pi^L(s^{-1}\UN^{(2)}\Sr\hat\UN)\cr
 &=(\hat\UN\Sl S(S(\UN^{(1)})\Sl\hat S(\sigma^{-1})))
   \hat S(\hat\UN^{(1)})\langle s^{-1}\UN^{(2)},\hat\UN^{(2)}\rangle\cr
 &=\hat S(\hat\UN^{(1)})\langle s^{-1}\UN^{(2)},
   \hat\UN^{(2)}(\hat\UN\Sl(\sigma^{-1}\Sr S^2(\UN^{(1)})))\rangle\cr
 &=\hat S(\hat\UN^{(1)})\langle s^{-1}\UN^{(2)},
   \hat\UN^{(2)}\Sl(\sigma^{-1}\Sr S^2(\UN^{(1)}))\rangle\cr
 &=\hat S(\hat\UN^{(1)})\langle(\sigma^{-1}\Sr S^2(\UN^{(1)})) 
   s^{-1}\UN^{(2)},\hat\UN^{(2)}\rangle
  =\hat S(\hat\UN^{(1)})\langle S^2(\UN^{(1)}) 
   (\sigma\Sr s^{-1}\UN^{(2)}),\hat\UN^{(2)}\sigma^{-1}\rangle\cr
 &=\hat S(\sigma^{-1}\hat\UN^{(1)}\sigma)\langle S^2(\UN^{(1)}) 
   (\sigma\Sr s^{-1}\UN^{(2)}),\hat S^{-1}(\sigma)\hat\UN^{(2)}\rangle\cr
 &=\hat S(\sigma^{-1}\hat\UN^{(1)}\sigma)\langle S^2(\UN^{(1)}) 
   (\sigma\Sr s^{-1}\UN^{(2)})\Sl\hat S^{-1}(\sigma),\hat\UN^{(2)}
   \rangle\cr
 &=\hat S(\sigma^{-1}\hat\UN^{(1)}\sigma)\langle S^2(\UN^{(1)}) 
   (\sigma\Sr s^{-1}\UN^{(2)}ss^{-1}\Sl\hat S^{-1}(\sigma)),
   \hat\UN^{(2)}\rangle\cr
 &=\hat S(\sigma^{-1}\hat\UN^{(1)}\sigma)\langle S^2(\UN^{(1)}) 
   S^4(\UN^{(2)})S^4(s^{-1}),\hat\UN^{(2)}\rangle\cr
 &=\hat S(\sigma^{-1}\hat\UN^{(1)}\sigma)\langle\Pi^L(S^2(\UN^{(1)}) 
   S^4(\UN^{(2)})S^4(s^{-1})),\hat\UN^{(2)}\rangle\cr
 &=\hat S(\sigma^{-1}\hat\UN^{(1)}\sigma)\langle S^4(\UN^{(2)}) 
   S^3(\UN^{(1)}),\hat\UN^{(2)}\rangle
  =\hat S(\sigma^{-1}\hat\UN^{(1)}\sigma)\langle\UN,\hat\UN^{(2)}\rangle
  =\hat\UN,\cr}
$$
where we used (1.10) in the second equality, (1.12) and the anticoalgebra map
property of the antipode in the third one, algebra automorphism properties 
of Sweedler actions by left grouplike elements in the sixth, (5.1--2b) in the 
seventh and twelfth, the Radford formula in the tenth and (1.7) in the 
eleventh equality.\qed

\bigskip\noindent
{\bf Appendix A}\medskip

Here we give examples of finite dimensional WHAs of $A=A^L\otimes A^R\equiv
B\otimes B^{op}$ type, where $B$ is a separable $k$-algebra equipped with a 
non-degenerate functional $E\colon B\to k$ of index $\UN$ (see Appendix of 
[2]), having antipode of infinite order. 

Let $B=M_n({\bf R})$, i.e. a full matrix algebra over the real field, and let
${\rm tr}\colon B\to {\bf R}$ denote the trace functional with 
${\rm tr}(\UN)=n$. Any invertible element $t\in B$ with ${\rm tr}(t^{-1})=1$ 
defines a non-degenerate functional $E\colon B\to {\bf R}$ by
$$
  E(x):={\rm tr}(tx),\qquad x\in B,\eqno(A.1)
$$
which has index $\UN$. Indeed, if $\{ e_{ab}\}$ is a set of matrix units
then $\{ e_i\}_i\equiv\{ t^{-1}e_{ab}\}_{(a,b)}$ and  
$\{ f_i\}_i\equiv\{ e_{ba}\}_{(a,b)}$ are dual ${\bf R}$-bases of $B$ with
respect to $E$, $E(e_if_j)=\delta_{ij}$, and the index of $E$ is
$$
  {\rm Ind}\, E:=\sum_if_ie_i=\sum_{a,b}e_{ba}t^{-1}e_{ab}={\rm tr}(t^{-1})
  \sum_b e_{bb}=\UN.\eqno(A.2)
$$
The Nakayama automorphism $\theta$ of $E$ defined by $E(xy)=:E(y\theta(x)); x,y
\in B$ is inner,
$$
  \theta(x)=txt^{-1},\qquad x\in B\eqno(A.3)
$$
due to the form (A.1) of $E$. We can construct the WHA $B\otimes B^{op}$ 
[2]: it is the ${\bf R}$-linear space $B\otimes B$ with structure maps
$$\eqalign{
  (x_1\otimes x_2)(y_1\otimes y_2)&:=x_1y_1\otimes y_2x_2,\cr
  \Delta(x_1\otimes x_2)&:=\sum_i (x_1\otimes f_i)\otimes(e_i\otimes x_2),\cr
  \varepsilon(x_1\otimes x_2)&:=E(x_1x_2),\cr
  S(x_1\otimes x_2)&:=x_2\otimes\theta(x_1).\cr}\eqno(A.4)
$$
Clearly, $S^2=\theta\otimes\theta$, therefore the form (A.3) of the Nakayama
automorphism $\theta$ shows that the order of the antipode $S$ is finite
iff $t^m\in{\rm Center}\, B$ for a certain positive integer $m$. However,
this is not the case for a generic invertible $t\in B=M_n({\bf R})$ with 
${\rm tr}(t^{-1})$. 

Although the order of the antipode is not finite in the generic case, 
already $S^2$ is an inner automorphism by a grouplike element in the trivial 
subalgebra $A^T$, which, in this case, is equal to $A$ itself. Indeed 
$$
  S^{2}(x\otimes y)=\theta(x)\otimes\theta(y)
  =(t\otimes t^{-1})(x\otimes y)(t^{-1}\otimes t),\eqno(A.5)
$$
and $t\otimes t^{-1}=(t\otimes \UN)S(t^{-1}\otimes \UN)$ with 
$t\otimes \UN=S^{2}(t\otimes \UN)\in A^L$. Therefore $t\otimes t^{-1}$ is
a grouplike element in the trivial subalgebra $A^T$ by Corollary 5.6.
 
\bigskip\noindent
{\bf Appendix B}\medskip

Here we give the generalization of the cyclic module [4] 
$A^\natural_{(\sigma,s)}$ for weak Hopf algebras having a modular pair 
$(\sigma,s)$ in involution. The details will be published elsewhere.

Let $A$ be a weak Hopf algebra. The pair $(\sigma,s)\in G(\hat A)\times G(A)$
of grouplike elements is called a {\it modular pair} for $A$ if
$$
  \sigma\Sr s=s=s\Sl\sigma,\qquad s\Sr\sigma=\sigma=\sigma\Sl s.\eqno(B.1)
$$
They form a {\it modular pair in involution} if they implement the square
of the antipode
$$
  S^2(a)=\sigma\Sr sas^{-1}\Sl\sigma^{-1},\quad a\in A.\eqno(B.2)
$$
Clearly, a modular pair (in involution) is a self-dual notion for WHAs. 

The identity (B.2) is a kind of square root of the Radford formula, hence, 
modular pairs in involution do not exist for arbitrary WHAs. However,
there is a wide class of WHAs having such a pair. For example, in a weak 
Hopf $C^*$-algebra $A$ there is a canonical grouplike element $g\in A$ 
implementing $S^2$ on $A$ [2], hence $(\hat\UN, g)$ is a modular pair 
in involution for $A$. Another example is as follows: let $A$ be 
a WHA over $k$ and let the WHA $A_G:=\langle A^T,G_R(A)\rangle$ be the 
subWHA of $A$ generated by the trivial subWHA $A^T$ and by (a subgroup of) 
the right grouplike elements $G_R(A)$ in $A$. Then $(\hat\UN,t)$ with $t\in
G^T(A)$ defined in (5.31) is a modular pair in involution for $A_G$, 
because $t$ implements $S^2$ for $A^T$ and $G_R(A)$ due to (5.31). 

\smallskip\noindent
{\bf Proposition} {\sl Let $A$ be a WHA over the field $k$ and 
$(\sigma,s)\in G(\hat A)\times G(A)$
be a modular pair in involution. Let the cochains $C^n_{(\sigma,s)}(A),\ 
n\geq 0$ be defined by the $n$-fold product of the left regular module
${_AA}$, i.e. the $k$-linear spaces
$$\eqalign{
  C^0_{(\sigma,s)}(A)&:=A^L\cr
  C^n_{(\sigma,s)}(A)&:=A\times A\times\dots\times A\equiv
                       \Delta^{n-1}(\UN)\cdot(A\otimes A\otimes\dots\otimes A).
  \cr}\eqno(B.3)
$$
The face operators $\delta^{(n)}_i\colon C^{n-1}_{(\sigma,s)}(A)\to
C^n_{(\sigma,s)}(A),\ 0\leq i\leq n$ are 
$$\eqalign{
  \delta_0^{(1)}(x_L)&:=\bar\Pi^R(x_L),\cr
  \delta_1^{(1)}(x_L)&:=x_Ls,\cr
  \delta_0^{(n)}(a_1\otimes\dots\otimes a_{n-1})&:=
     \UN^{(1)}\otimes\UN^{(2)}a_1\otimes a_2\otimes\dots\otimes a_{n-1},
     \quad 1 < n,\cr
  \delta_i^{(n)}(a_1\otimes\dots\otimes a_{n-1})&:=
     a_1\otimes a_2\otimes\dots\otimes\Delta(a_i)\otimes\dots\otimes a_{n-1},
     \quad 1\leq i< n, 1 < n,\cr
  \delta_n^{(n)}(a_1\otimes\dots\otimes a_{n-1})&:=
     a_1\otimes\dots\otimes a_{n-2}\otimes\UN^{(1)}a_{n-1}\otimes\UN^{(2)}s,
     \quad 1 < n,;\cr}\eqno(B.4)
$$
the degeneracy operators $\sigma^{(n)}_i\colon C^{n+1}_{(\sigma,s)}(A)\to
C^n_{(\sigma,s)}(A),\ 0\leq i\leq n$ are
$$\eqalign{
  \sigma_0^{(0)}(a)&:=\Pi^L(a),\cr
  \sigma_i^{(n)}(a_1\otimes\dots\otimes a_{n+1})&:=
     a_1\otimes\dots\otimes\Pi^L(a_{i+1})a_{i+2}\otimes\dots\otimes a_{n+1},
     \quad 0\leq i< n, 0 < n,\cr
  \sigma_n^{(n)}(a_1\otimes\dots\otimes a_{n+1})&:=
     a_1\otimes\dots\otimes a_{n-1}\otimes\bar\Pi^R(a_{n+1})a_n,
     \quad  0 < n,\cr}\eqno(B.5)
$$
and the cyclic operators $\tau_{(n)}\colon C^n_{(\sigma,s)}(A)\to
C^n_{(\sigma,s)}(A)$ are given by
$$\eqalign{
  \tau_{(0)}(x_L)&:=x_L,\cr
  \tau_{(n)}(a_1\otimes\dots\otimes a_n)&:=
     \Delta^{(n-1)}(S(a_1\Sl\sigma))\cdot(a_2\otimes\dots\otimes 
     a_n\otimes s),
     \quad n\geq 1.\cr}\eqno(B.6)
$$
With the definitions (B.3--6) $A^\natural_{(\sigma,s)}\equiv
\{C^n_{(\sigma,s)}(A)\}_{n\geq 0}$ becomes a $\Lambda$-module, where
$\Lambda$ is the cyclic category.}

\bigskip\noindent
{\bf References}\medskip

\item{[1]} G. B\"ohm and K. Szlach\'anyi, A coassociative $C^*$-quantum group
with nonintegral dimensions, {\it Lett. Math. Phys.} {\bf 35} (1996), 437-456.
\smallskip

\item{[2]} G. B\"ohm, F. Nill and K. Szlach\'anyi, Weak Hopf algebras I. 
Integral theory and $C^*$-structure, {\it J. Algebra} {\bf 221} (1999), 
385-438.
\smallskip

\item{[3]} G. B\"ohm, ``Weak $C^*$-Hopf algebras and their application  
to spin models,'' Ph.D thesis, Budapest, 1997.
\item{} G. B\"ohm, Doi--Hopf modules over weak Hopf algebras, {\it Comm. 
Algebra} {\bf 28} (2000), 4687-4698.
\smallskip

\item{[4]} A. Connes and H. Moscovici, Cyclic cohomology and Hopf algebras,
{\it Lett. Math. Phys.} {\bf 48} (1999), 97-108.
\item{} A. Connes and H. Moscovici, Differentiable cyclic cohomology and Hopf 
algebraic structures in transverse geometry, {\it Essays on geometry and 
related topics}, Monogr. Enseign. Math. {\bf 38} (2001), 217-255. 
\smallskip

\item{[5]} C.W. Curtis and I. Reiner, ``Representation Theory of Finite 
Groups and Associative Algebras,'' John Wiley \& Sons, 1962.
\smallskip

\item{[6]} C.W. Curtis and I. Reiner, ``Methods of Representation Theory,'' 
John Wiley \& Sons, 1990.
\smallskip

\item{[7]} P. Deligne and J.S. Milne, Tannakian categories, {\it in} Lecture 
Notes in Math. 900, p. 101, Springer, Berlin/New York, 1982.
\smallskip

\item{[8]} J. Fuchs, A. Ganchev and P. Vecserny\'es, Rational Hopf Algebras: 
Polynomial equations, gauge fixing, and low-dimensional
examples, {\it Internat. J. Modern Phys.} {\bf A10} (1995), 3431-3476. 
\smallskip

\item{[9]} L. Kadison, ``New Examples of Frobenius Extensions,''  
University Lecture Series, Vol. 14, Am. Math. Soc., Providence, 1999. 
\smallskip

\item{[10]} R.G. Larson and M.E. Sweedler, An associative orthogonal 
bilinear form for Hopf algebras, {\it Amer. J. Math.} {\bf 91} (1969), 75-93.
\smallskip

\item{[11]} G. Mack and V. Schomerus, Quasi Hopf quantum symmetry in 
quantum theory, {\it Nucl. Phys.} {\bf B370} (1992), 185. 
\smallskip

\item{[12]} D. Nikshych and L. Vainerman, Finite quantum grupoids and their 
applications, {\it Math. Sci. Res. Inst. Publ.} {\bf 43} (2002), 211-262.
\smallskip

\item{[13]} D. Nikshych, On the structure  of weak Hopf algebras, 
{\it Adv. Math.} {\bf 170} (2002), 257-286. 
\smallskip

\item{[14]} R.S. Pierce, ``Associative Algebras,'' Graduate Texts in 
Mathematics, Vol.88, Springer, Berlin/New York, 1982.
\smallskip

\item{[15]} D.E. Radford, The order of the antipode of a finite dimensional
Hopf algebra is finite, {\it Amer. J. Math.} {\bf 98} (1976), 333-355.
\smallskip

\item{[16]} D.E. Radford, Minimal quasitriangular Hopf algebras, 
{\it J. Algebra} {\bf 157} (1993), 285-315.
\smallskip

\item{[17]} M.E. Sweedler, ``Hopf Algebras,'' Benjamin, Elmsford, NY, 1969.
\smallskip

\item{[18]}  K. Szlach\'anyi, Weak Hopf algebras, {\it in} ``Operator Algebras 
and Quantum Field Theory'', (S. Doplicher, R. Longo, J.E. Roberts, and 
L. Zsid\'o, Eds.) International Press, 1996.
\smallskip

\item{[19]} P. Vecserny\'es, On the quantum symmetry of the chiral Ising 
model, {\it Nucl. Phys.} {\bf B415} (1994), 557-588. 
\smallskip

\item{[20]} Y. Watatani, Index for $C^*$-subalgebras, {\it Memoirs Amer. Math. 
Soc.} {\bf 424} (1990). 
\smallskip

\item{[21]} D.N. Yetter, Framed tangles and a theorem of Deligne on braided 
deformations of Tannakian categories, {\it Contemp. Math.} {\bf 134} (1992), 
325-349. 
\smallskip

\bye